\newcommand{\xbD}{\Delta}

\newcommand{\xbO}{\Omega}
\newcommand{\xbP}{\Pi}

\newcommand{\xba}{\alpha}
\newcommand{\xbb}{\beta}

\newcommand{\xbe}{\in}
\newcommand{\xbf}{\phi}
\newcommand{\xbg}{\gamma}

\newcommand{\xbj}{\vartheta}
\newcommand{\xbk}{\kappa}

\newcommand{\xbm}{\mu}

\newcommand{\xbo}{\omega}

\newcommand{\xbq}{\psi}

\newcommand{\xbs}{\sigma}
\newcommand{\xbt}{\tau}

\newcommand{\xCK}{\times}

\newcommand{\xCN}{\neg}
\newcommand{\xCQ}{\emptyset}

\newcommand{\xCf}{\hspace{0.1em}}

\newcommand{\xcA}{\forall}

\newcommand{\xcC}{\not\subseteq}

\newcommand{\xcE}{\exists}

\newcommand{\xcH}{\not\Rightarrow}
\newcommand{\xcI}{\not\Leftarrow}

\newcommand{\xcL}{\not\vdash}
\newcommand{\xcM}{\not\models}
\newcommand{\xcN}{\hspace{0.2em}\not\sim\hspace{-0.9em}\mid\hspace{0.8em}}

\newcommand{\xcS}{\bigcap}
\newcommand{\xcT}{\bot}
\newcommand{\xcU}{\bigwedge}
\newcommand{\xcV}{\bigcup}

\newcommand{\xca}{\infty}
\newcommand{\xcb}{\subset}
\newcommand{\xcc}{\subseteq}
\newcommand{\xcd}{\supseteq}
\newcommand{\xce}{\not\in}
\newcommand{\xcf}{\supset}
\newcommand{\xcg}{\geq}
\newcommand{\xch}{\Rightarrow}
\newcommand{\xci}{\Leftarrow}
\newcommand{\xcj}{\Leftrightarrow}
\newcommand{\xck}{\leq}
\newcommand{\xcl}{\vdash}
\newcommand{\xcm}{\models}
\newcommand{\xcn}{\hspace{0.2em}\sim\hspace{-0.9em}\mid\hspace{0.58em}}

\newcommand{\xco}{\vee}
\newcommand{\xcp}{\rightarrow}

\newcommand{\xcr}{\leftrightarrow}
\newcommand{\xcs}{\cap}
\newcommand{\xcu}{\wedge}
\newcommand{\xcv}{\cup}

\newcommand{\xcz}{\Box}

\newcommand{\xDH}{\item }

\newcommand{\xDN}{\ominus}

\newcommand{\xdC}{\mbox{\boldmath$C$}}
\newcommand{\xdD}{\mbox{\boldmath$D$}}

\newcommand{\xdf}{{\cal F}}

\newcommand{\xdi}{{\cal I}}

\newcommand{\xdl}{{\cal L}}
\newcommand{\xdm}{{\cal M}}
\newcommand{\xdn}{{\cal N}}

\newcommand{\xdp}{{\cal P}}

\newcommand{\xdu}{{\cal U}}

\newcommand{\xdx}{{\cal X}}
\newcommand{\xdy}{{\cal Y}}
\newcommand{\xdz}{{\cal Z}}

\newcommand{\xEH}{ & }
\newcommand{\xEI}{\begin{itemize}}
\newcommand{\xEJ}{\end{itemize}}
\newcommand{\xEP}{ \\ }

\newcommand{\xEd}{\neq}
\newcommand{\xEh}{\begin{enumerate}}
\newcommand{\xEj}{\end{enumerate}}

\newcommand{\xeA}{\nabla}

\newcommand{\xeb}{\prec}
\newcommand{\xec}{\preceq}
\newcommand{\xed}{\succeq}
\newcommand{\xee}{\succ}

\newcommand{\xex}{\lceil}

\newcommand{\xFO}{\parallel}

\newcommand{\xfA}{\mid}
\newcommand{\xfB}{\uparrow}

\newcommand{\xfb}{\downarrow}

\newcommand{\Xl}{\ldots}

\newcommand{\ol}{\overline}

\newcommand{\wt}{\overbrace}

\newcommand{\bl}{\begin{lemma} \rm}
\newcommand{\el}{\end{lemma}}
\newcommand{\br}{\begin{remark} \rm}
\newcommand{\er}{\end{remark}}
\newcommand{\be}{\begin{example} \rm}
\newcommand{\ee}{\end{example}}
\newcommand{\bco}{\begin{corollary} \rm}
\newcommand{\eco}{\end{corollary}}
\newcommand{\bc}{\begin{claim} \rm}
\newcommand{\ec}{\end{claim}}
\newcommand{\bfa}{\begin{fact} \rm}
\newcommand{\efa}{\end{fact}}
\newcommand{\bp}{\begin{proposition} \rm}
\newcommand{\ep}{\end{proposition}}
\newcommand{\bd}{\begin{definition} \rm}
\newcommand{\ed}{\end{definition}}
\newcommand{\bcs}{\begin{construction} \rm}
\newcommand{\ecs}{\end{construction}}
\newcommand{\bcd}{\begin{condition} \rm}
\newcommand{\ecd}{\end{condition}}
\newcommand{\bt}{\begin{theorem} \rm}
\newcommand{\et}{\end{theorem}}
\newcommand{\bn}{\begin{notation} \rm}
\newcommand{\en}{\end{notation}}
\newcommand{\bfi}{\begin{bild} \rm}
\newcommand{\efi}{\end{bild}}
\newcommand{\bsta}{\begin{statement} \rm}
\newcommand{\esta}{\end{statement}}
\newcommand{\bcom}{\begin{comment} \rm}
\newcommand{\ecom}{\end{comment}}
\newcommand{\bdia}{\begin{diagram} \rm}
\newcommand{\edia}{\end{diagram}}

\newcommand{\bfc}{\begin{figure}[htb] \begin{center}}
\newcommand{\efc}{\end{center} \end{figure}}

\documentstyle[12pt]{article}
\title{
Roadmap for preferential logics
}

\author{Dov M Gabbay
\thanks{
Dov.Gabbay@kcl.ac.uk, www.dcs.kcl.ac.uk/staff/dg
} \\
King's College, London
\thanks{
Department of Computer Science, King's College London, Strand,
London WC2R 2LS, UK
} \\ \\
Karl Schlechta
\thanks{
ks@cmi.univ-mrs.fr, karl.schlechta@web.de, http://www.cmi.univ-mrs.fr/ $\sim$ ks
} \\
Laboratoire d'Informatique Fondamentale de Marseille
\thanks{
UMR 6166, CNRS and Universit\'{e} de Provence,
Address: CMI, 39, rue Joliot-Curie, F-13453 Marseille Cedex 13, France
}
}

\begin{document}

\newtheorem{lemma}{Lemma}[section]
\newtheorem{theorem}[lemma]{Theorem}
\newtheorem{proposition}[lemma]{Proposition}
\newtheorem{corollary}[lemma]{Corollary}
\newtheorem{claim}[lemma]{Claim}
\newtheorem{fact}[lemma]{Fact}
\newtheorem{remark}[lemma]{Remark}
\newtheorem{definition}{Definition}[section]
\newtheorem{construction}{Construction}[section]
\newtheorem{condition}{Condition}[section]
\newtheorem{example}{Example}[section]
\newtheorem{notation}{Notation}[section]
\newtheorem{bild}{Figure}[section]
\newtheorem{comment}{Comment}[section]
\newtheorem{statement}{Statement}[section]
\newtheorem{diagram}{Diagram}[section]

\maketitle

\renewcommand{\labelenumi}
  {(\arabic{enumi})}
\renewcommand{\labelenumii}
  {(\arabic{enumi}.\arabic{enumii})}
\renewcommand{\labelenumiii}
  {(\arabic{enumi}.\arabic{enumii}.\arabic{enumiii})}
\renewcommand{\labelenumiv}
  {(\arabic{enumi}.\arabic{enumii}.\arabic{enumiii}.\arabic{enumiv})}

\setcounter{secnumdepth}{3}
\setcounter{tocdepth}{3}

\tableofcontents

%
%
%

\section{
Introduction
}

\subsection{
Purpose of the paper
}

The purpose of these pages is to give the reader a systematic overview of
logical and algebraic rules used in nonmonotonic and related logics.
We try to give orientation in a multitude of sometimes quite
similar rules, and in translating the different versions to each other.

The emphasis is on systematisation, and we will not go into deeper
completeness proofs.
\subsection{
Organisation of the paper
}

The article is built around several tables.

They show

 \xEh

 \xDH connections between semantical and proof theoretical rules, but
also their (sometimes subtle) differences,

 \xDH connections between different semantical rules, but again also their
(sometimes subtle) differences.

 \xEj

Further tables summarize

 \xEh

 \xDH representation results for preferential structures,

 \xDH connections between the different concepts of AGM revision,

 \xDH results for distance based theory revision,

 \xDH connections between filters, the notion of size, and nonmonotonic
logic.

 \xEj

The last table is probably the most innovative part of the paper, and it
led to the introduction of perhaps new rules (variants of the $ \xCf (OR)$
rule).

Yet, as tables go, the emphasis is more on systematisation than on
novelty.

The ``subtle'' part of the comparisons and differences concerns mostly
domain closure problems:

 \xEI

 \xDH is a domain closed under finite union?

 \xDH does the operator preserve definability, i.e. is $f(M(T))=M(T' )$
for some
$T' $ - where $T,T' $ are sets of formulas, and $M(T)$ is the set of
classical
models of $T?$

 \xDH is the complement of $M(T)$ again some $M(T' )?$ etc.

 \xEJ

Thus, as a good roadmap should, the article points out easy ways to go
from
A to $B,$ but also puts up warning signs where there are problems ahead.
\subsection{
Summary of the tables
}

 \xEh

 \xDH Tables about rules for nonmonotonic logics

 \xEI

 \xDH Definition of the rules, Definition \ref{Definition Log-Cond}

 \xDH Connections between the different semantical rules,
Fact \ref{Fact Mu-Base}

 \xDH Translations between the logical and the semantical variants,
Proposition \ref{Proposition Alg-Log}

 \xEJ

 \xDH Summary of representation by preferential structures,
Table \ref{Table Pref-Representation-With-Ref}

 \xDH Tables about rules for theory revision

 \xEI

 \xDH AGM revision, definitions of logical and semantical versions,
Definition \ref{Definition AGM}

 \xDH Interdefinability of AGM concepts,
Proposition \ref{Proposition AGM-Equiv}

 \xDH Definition of rules for distance based revision,
Conditions \ref{Conditions TR-Dist}

 \xDH Translation between semantical and logical versions,
Proposition \ref{Proposition TR-Alg-Log}

 \xEJ

 \xDH Tables concerning size and coherence

 \xEI

 \xDH Definitions of (weak) filters, ideals, and coherence conditions,
Definition \ref{Definition Filter}

 \xDH Correspondence between coherence conditions and semantical rules
for nonmonotonic logics, Fact \ref{Fact Ref-Class-Mu}

 \xEJ

 \xEj

\section{
Generalities
}

\index{Definition Alg-Base}

\bd

$\hspace{0.01em}$


\label{Definition Alg-Base}

We use $ \xdp $ to denote the power set operator,
$ \xbP \{X_{i}:i \xbe I\}$ $:=$ $\{g:$ $g:I \xcp \xcV \{X_{i}:i \xbe I\},$
$ \xcA i \xbe I.g(i) \xbe X_{i}\}$ is the general cartesian
product, $card(X)$ shall denote the cardinality of $X,$ and $V$ the
set-theoretic
universe we work in - the class of all sets. Given a set of pairs $ \xdx
,$ and a
set $X,$ we denote by $ \xdx \xex X:=\{<x,i> \xbe \xdx:x \xbe X\}.$ When
the context is clear, we
will sometime simply write $X$ for $ \xdx \xex X.$

$A \xcc B$ will denote that $ \xCf A$ is a subset of $B$ or equal to $B,$
and $A \xcb B$ that $ \xCf A$ is
a proper subset of $B,$ likewise for $A \xcd B$ and $A \xcf B.$

Given some fixed set $U$ we work in, and $X \xcc U,$ then $ \xdC (X):=U-X$
.

If $ \xdy \xcc \xdp (X)$ for some
$X,$ we say that $ \xdy $ satisfies

$( \xcs )$ iff it is closed under finite intersections,

$( \xcS )$ iff it is closed under arbitrary intersections,

$( \xcv )$ iff it is closed under finite unions,

$( \xcV )$ iff it is closed under arbitrary unions,

$( \xdC )$ iff it is closed under complementation.

We will sometimes write $A=B \xFO C$ for: $A=B,$ or $A=C,$ or $A=B \xcv
C.$

We make ample and tacit use of the Axiom of Choice.
\index{Definition Rel-Base}

\ed

\bd

$\hspace{0.01em}$


\label{Definition Rel-Base}

$ \xeb^{*}$ will denote the transitive closure of the relation $ \xeb.$
If a relation $<,$
$ \xeb,$ or similar is given, $a \xcT b$ will express that a and $b$ are
$<-$ (or $ \xeb -)$
incomparable - context will tell. Given any relation $<,$ $ \xck $ will
stand for
$<$ or $=,$ conversely, given $ \xck,$ $<$ will stand for $ \xck,$ but
not $=,$ similarly
for $ \xeb $ etc.
\index{Definition Tree-Base}

\ed

\bd

$\hspace{0.01em}$


\label{Definition Tree-Base}

A child (or successor) of an element $x$ in a tree $t$ will be a direct
child in $t.$
A child of a child, etc. will be called an indirect child. Trees will be
supposed to grow downwards, so the root is the top element.
\index{Definition Seq-Base}

\ed

\bd

$\hspace{0.01em}$


\label{Definition Seq-Base}

A subsequence $ \xbs_{i}:i \xbe I \xcc \xbm $ of a sequence $ \xbs_{i}:i
\xbe \xbm $ is called cofinal, iff
for all $i \xbe \xbm $ there is $i' \xbe I$ $i \xck i'.$

Given two sequences $ \xbs_{i}$ and $ \xbt_{i}$ of the same length, then
their Hamming distance
is the quantity of $i$ where they differ.
\newpage

\section{
Logical rules
}

\index{Definition Log-Base}

\ed

\bd

$\hspace{0.01em}$


\label{Definition Log-Base}

We work here in a classical propositional language $ \xdl,$ a theory $T$
will be an
arbitrary set of formulas. Formulas will often be named $ \xbf,$ $ \xbq
,$ etc., theories
$T,$ $S,$ etc.

$v( \xdl )$ will be the set of propositional variables of $ \xdl.$

$M_{ \xdl }$ will be the set of (classical) models of $ \xdl,$ $M(T)$ or
$M_{T}$
is the set of models of $T,$ likewise $M( \xbf )$ for a formula $ \xbf.$

$ \xdD_{ \xdl }:=\{M(T):$ $T$ a theory in $ \xdl \},$ the set of definable
model sets.

Note that, in classical propositional logic, $ \xCQ,M_{ \xdl } \xbe
\xdD_{ \xdl },$ $ \xdD_{ \xdl }$ contains
singletons, is closed under arbitrary intersections and finite unions.

An operation $f: \xdy \xcp \xdp (M_{ \xdl })$ for $ \xdy \xcc \xdp (M_{
\xdl })$ is called definability
preserving, $ \xCf (dp)$ or $( \xbm dp)$ in short, iff for all $X \xbe
\xdD_{ \xdl } \xcs \xdy $ $f(X) \xbe \xdD_{ \xdl }.$

We will also use $( \xbm dp)$ for binary functions $f: \xdy \xCK \xdy \xcp
\xdp (M_{ \xdl })$ - as needed
for theory revision - with the obvious meaning.

$ \xcl $ will be classical derivability, and

$ \ol{T}:=\{ \xbf:T \xcl \xbf \},$ the closure of $T$ under $ \xcl.$

$Con(.)$ will stand for classical consistency, so $Con( \xbf )$ will mean
that
$ \xbf $ is clasical consistent, likewise for $Con(T).$ $Con(T,T' )$ will
stand for
$Con(T \xcv T' ),$ etc.

Given a consequence relation $ \xcn,$ we define

$ \ol{ \ol{T} }:=\{ \xbf:T \xcn \xbf \}.$

(There is no fear of confusion with $ \ol{T},$ as it just is not useful to
close
twice under classical logic.)

$T \xco T':=\{ \xbf \xco \xbf ': \xbf \xbe T, \xbf ' \xbe T' \}.$

If $X \xcc M_{ \xdl },$ then $Th(X):=\{ \xbf:X \xcm \xbf \},$ likewise
for $Th(m),$ $m \xbe M_{ \xdl }.$
\index{Fact Log-Base}

\ed

We recollect and note:

\bfa

$\hspace{0.01em}$


\label{Fact Log-Base}

Let $ \xdl $ be a fixed propositional language, $ \xdD_{ \xdl } \xcc X,$ $
\xbm:X \xcp \xdp (M_{ \xdl }),$ for a $ \xdl -$theory $T$
$ \ol{ \ol{T} }:=Th( \xbm (M_{T})),$ let $T,$ $T' $ be arbitrary theories,
then:

(1) $ \xbm (M_{T}) \xcc M_{ \ol{ \ol{T} }}$,

(2) $M_{T} \xcv M_{T' }=M_{T \xco T' }$ and $M_{T \xcv T' }=M_{T} \xcs
M_{T' }$,

(3) $ \xbm (M_{T})= \xCQ $ $ \xcr $ $ \xcT \xbe \ol{ \ol{T} }$.

If $ \xbm $ is definability preserving or $ \xbm (M_{T})$ is finite, then
the following also hold:

(4) $ \xbm (M_{T})=M_{ \ol{ \ol{T} }}$,

(5) $T' \xcl \ol{ \ol{T} }$ $ \xcr $ $M_{T' } \xcc \xbm (M_{T}),$

(6) $ \xbm (M_{T})=M_{T' }$ $ \xcr $ $ \ol{T' }= \ol{ \ol{T} }.$
$ \xcz $
\\[3ex]
\index{Fact Th-Union}

\efa

\bfa

$\hspace{0.01em}$


\label{Fact Th-Union}

Let $A,B \xcc M_{ \xdl }.$

Then $Th(A \xcv B)$ $=$ $Th(A) \xcs Th(B).$
\index{Fact Th-Union Proof}

\efa

\subparagraph{
Proof
}

$\hspace{0.01em}$


$ \xbf \xbe Th(A \xcv B)$ $ \xcj $ $A \xcv B \xcm \xbf $ $ \xcj $ $A \xcm
\xbf $ and $B \xcm \xbf $ $ \xcj $ $ \xbf \xbe Th(A)$ and $ \xbf \xbe
Th(B).$

$ \xcz $
\\[3ex]
\index{Fact Log-Form}

\bfa

$\hspace{0.01em}$


\label{Fact Log-Form}

Let $X \xcc M_{ \xdl },$ $ \xbf, \xbq $ formulas.

(1) $X \xcs M( \xbf ) \xcm \xbq $ iff $X \xcm \xbf \xcp \xbq.$

(2) $X \xcs M( \xbf ) \xcm \xbq $ iff $M(Th(X)) \xcs M( \xbf ) \xcm \xbq
.$

(3) $Th(X \xcs M( \xbf ))= \ol{Th(X) \xcv \{ \xbf \}}$

(4) $X \xcs M( \xbf )= \xCQ $ $ \xcj $ $M(Th(X)) \xcs M( \xbf )= \xCQ $

(5) $Th(M(T) \xcs M(T' ))= \ol{T \xcv T' }.$
\index{Fact Log-Form Proof}

\efa

\subparagraph{
Proof
}

$\hspace{0.01em}$


(1) `` $ \xch $ '': $X=(X \xcs M( \xbf )) \xcv (X \xcs M( \xCN \xbf )).$ In
both parts holds $ \xCN \xbf \xco \xbq,$ so
$X \xcm \xbf \xcp \xbq.$ `` $ \xci $ '': Trivial.

(2) $X \xcs M( \xbf ) \xcm \xbq $ (by (1)) iff $X \xcm \xbf \xcp \xbq $
iff $M(Th(X)) \xcm \xbf \xcp \xbq $ iff (again by (1))
$M(Th(X)) \xcs M( \xbf ) \xcm \xbq.$

(3) $ \xbq \xbe Th(X \xcs M( \xbf ))$ $ \xcj $ $X \xcs M( \xbf ) \xcm \xbq
$ $ \xcj_{(2)}$ $M(Th(X) \xcv \{ \xbf \})=M(Th(X)) \xcs M( \xbf ) \xcm
\xbq $ $ \xcj $
$Th(X) \xcv \{ \xbf \} \xcl \xbq.$

(4) $X \xcs M( \xbf )= \xCQ $ $ \xcj $ $X \xcm \xCN \xbf $ $ \xcj_{(1)}$
$M(Th(X)) \xcm \xCN \xbf $ $ \xcj $ $M(Th(X)) \xcs M( \xbf )= \xCQ.$

(5) $M(T) \xcs M(T' )=M(T \xcv T' ).$

$ \xcz $
\\[3ex]
\index{Fact Dp-Base}

\bfa

$\hspace{0.01em}$


\label{Fact Dp-Base}

If $X=M(T),$ then $M(Th(X))=X.$
\index{Fact Dp-Base Proof}

\efa

\subparagraph{
Proof
}

$\hspace{0.01em}$


$X \xcc M(Th(X))$ is trivial. $Th(M(T))= \ol{T}$ is trivial by classical
soundness and
completeness. So $M(Th(M(T))=M( \ol{T})=M(T)=X.$ $ \xcz $
\\[3ex]
\index{Example Not-Def}

\be

$\hspace{0.01em}$


\label{Example Not-Def}

If $v( \xdl )$ is infinite, and $m$ any model for $ \xdl,$ then $M:=M_{
\xdl }-\{m\}$ is not definable
by any theory $T.$ (Proof: Suppose it were, and let $ \xbf $ hold in $M,$
but not in $m,$ so in $m$ $ \xCN \xbf $ holds, but as $ \xbf $ is finite,
there is a model $m' $ in
$M$ which coincides on all propositional variables of $ \xbf $ with $m,$
so in $m' $ $ \xCN \xbf $
holds, too, a contradiction.) Thus, in the infinite case, $ \xdp (M_{ \xdl
}) \xEd \xdD_{ \xdl }.$

(There is also a simple cardinality argument, which shows that almost no
model sets are definable, but it is not constructive and thus less
instructive
than above argument. We give it nonetheless: Let $ \xbk:=card(v( \xdl
)).$ Then
there are $ \xbk $ many formulas, so $2^{ \xbk }$ many theories, and thus
$2^{ \xbk }$ many
definable model sets. But there are $2^{ \xbk }$ many models, so $(2^{
\xbk })^{ \xbk }$ many model
sets.)

$ \xcz $
\\[3ex]
\index{Definition Def-Clos}

\ee

\bd

$\hspace{0.01em}$


\label{Definition Def-Clos}

Let $ \xdy \xcc \xdp (Z)$ be given and closed under arbitrary
intersections.

For $A \xcc Z,$ let $ \wt{A}$ $:=$ $ \xcS \{X \xbe \xdy:A \xcc X\}.$

Intuitively, $Z$ is the set of all models for $ \xdl,$ $ \xdy $ is $
\xdD_{ \xdl }$, and $ \wt{A}=M(Th(A)),$
this is the intended application. Note that then $ \wt{ \xCQ }= \xCQ.$
\index{Fact Def-Clos}

\ed

\bfa

$\hspace{0.01em}$


\label{Fact Def-Clos}

(1) If $ \xdy \xcc \xdp (Z)$ is closed under arbitrary intersections and
finite unions,
$Z \xbe \xdy,$ $X,Y \xcc Z,$ then the following hold:

$(Cl \xcv )$ $ \wt{X \xcv Y}$ $=$ $ \wt{X} \xcv \wt{Y}$

$(Cl \xcs )$ $ \wt{X \xcs Y} \xcc \wt{X} \xcs \wt{Y},$ but usually not
conversely,

$ \xCf (Cl-)$ $ \wt{A}- \wt{B} \xcc \wt{A-B},$

$(Cl=)$ $X=Y$ $ \xcp $ $ \wt{X}= \wt{Y},$ but not conversely,

$(Cl \xcc 1)$ $ \wt{X} \xcc Y$ $ \xcp $ $X \xcc Y,$ but not conversely,

$(Cl \xcc 2)$ $X \xcc \wt{Y}$ $ \xcp $ $ \wt{X} \xcc \wt{Y}.$

(2) If, in addition, $X \xbe \xdy $ and $ \xdC X:=Z-X \xbe \xdy,$ then
the following two properties
hold, too:

$(Cl \xcs +)$ $ \wt{A} \xcs X= \wt{A \xcs X},$

$(Cl-+)$ $ \wt{A}-X= \wt{A-X}.$

(3) In the intended application, i.e. $ \wt{A}=M(Th(A)),$ the following
hold:

(3.1) $Th(X)$ $=$ $Th( \wt{X}),$

(3.2) Even if $A= \wt{A},$ $B= \wt{B},$ it is not necessarily true that $
\wt{A-B} \xcc \wt{A}- \wt{B}.$
\index{Fact Def-Clos Proof}

\efa

\subparagraph{
Proof:
}

$\hspace{0.01em}$


\label{Section Proof:}

$(Cl=),$ $(Cl \xcc 1),$ $(Cl \xcc 2),$ (3.1) are trivial.

$(Cl \xcv )$ Let $ \xdy (U):=\{X \xbe \xdy:U \xcc X\}.$ If $A \xbe \xdy
(X \xcv Y),$ then $A \xbe \xdy (X)$ and $A \xbe \xdy (Y),$ so
$ \wt{X \xcv Y}$ $ \xcd $ $ \wt{X} \xcv \wt{Y}.$
If $A \xbe \xdy (X)$ and $B \xbe \xdy (Y),$ then $A \xcv B \xbe \xdy (X
\xcv Y),$ so $ \wt{X \xcv Y}$ $ \xcc $ $ \wt{X} \xcv \wt{Y}.$

$(Cl \xcs )$ Let $X',Y' \xbe \xdy,$ $X \xcc X',$ $Y \xcc Y',$ then $X
\xcs Y \xcc X' \xcs Y',$ so $ \wt{X \xcs Y} \xcc \wt{X} \xcs \wt{Y}.$
For the converse, set $X:=M_{ \xdl }-\{m\},$ $Y:=\{m\}$ in Example \ref{Example
Not-Def}.

$ \xCf (Cl-)$ Let $A-B \xcc X \xbe \xdy,$ $B \xcc Y \xbe \xdy,$ so $A
\xcc X \xcv Y \xbe \xdy.$ Let $x \xce \wt{B}$ $ \xch $ $ \xcE Y \xbe \xdy
(B \xcc Y,$ $x \xce Y),$
$x \xce \wt{A-B}$ $ \xch $ $ \xcE X \xbe \xdy (A-B \xcc X,$ $x \xce X),$
so $x \xce X \xcv Y,$ $A \xcc X \xcv Y,$ so $x \xce \wt{A}.$ Thus, $x \xce
\wt{B},$ $x \xce \wt{A-B}$ $ \xch $
$x \xce \wt{A},$ or $x \xbe \wt{A}- \wt{B}$ $ \xch $ $x \xbe \wt{A-B}.$

$(Cl \xcs +)$ $ \wt{A} \xcs X \xcd \wt{A \xcs X}$ by $(Cl \xcs ).$
For `` $ \xcc $ '': Let $A \xcs X \xcc A' \xbe \xdy,$ then by closure
under $( \xcv ),$
$A \xcc A' \xcv \xdC X \xbe \xdy,$ $(A' \xcv \xdC X) \xcs X \xcc A'.$ So
$ \wt{A} \xcs X \xcc \wt{A \xcs X}.$

$(Cl-+)$ $ \wt{A-X}= \wt{A \xcs \xdC X}= \wt{A} \xcs \xdC X= \wt{A}-X$ by
$(Cl \xcs +).$

(3.2) Set $A:=M_{ \xdl },$ $B:=\{m\}$ for $m \xbe M_{ \xdl }$ arbitrary, $
\xdl $ infinite. So $A= \wt{A},$ $B= \wt{B},$ but
$ \wt{A-B}=A \xEd A-B.$

$ \xcz $
\\[3ex]
\index{Definition Log-Cond}

\bd

$\hspace{0.01em}$


\label{Definition Log-Cond}

We introduce here formally a list of properties of set functions on the
algebraic side, and their corresponding logical rules on the other side.

Recall that $ \ol{T}:=\{ \xbf:T \xcl \xbf \},$ $ \ol{ \ol{T} }:=\{ \xbf
:T \xcn \xbf \},$
where $ \xcl $ is classical consequence, and $ \xcn $ any other
consequence.

We show, wherever adequate, in parallel the formula version
in the left column, the theory version
in the middle column, and the semantical or algebraic
counterpart in the
right column. The algebraic counterpart gives conditions for a
function $f:\xdy\xcp\xdp (U)$, where $U$ is some set, and
$\xdy\xcc\xdp (U)$.

When the formula version is not commonly used, we omit it,
as we normally work only with the theory version.

Intuitively, $A$ and $B$ in the right hand side column stand for
$M(\xbf)$ for some formula $\xbf$, whereas $X$, $Y$ stand for
$M(T)$ for some theory $T$.

{\footnotesize

\begin{tabular}{|c|c|c|}

\hline

\multicolumn{3}{|c|}{Basics} \xEP

\hline

$(AND)$
\xEH
$(AND)$
\xEH
Closure under
\xEP

$ \xbf \xcn \xbq,  \xbf \xcn \xbq '   \xch $
\xEH
$ T \xcn \xbq, T \xcn \xbq '   \xch $
\xEH
finite
\xEP

$ \xbf \xcn \xbq \xcu \xbq ' $
\xEH
$ T \xcn \xbq \xcu \xbq ' $
\xEH
intersection
\xEP

\hline

$(OR)$ \xEH $(OR)$ \xEH $( \xbm OR)$ \xEP

$ \xbf \xcn \xbq,  \xbf ' \xcn \xbq   \xch $ \xEH
$ \ol{\ol{T}} \xcs \ol{\ol{T'}} \xcc \ol{\ol{T \xco T'}} $ \xEH
$f(X \xcv Y) \xcc f(X) \xcv f(Y)$
\xEP

$ \xbf \xco \xbf ' \xcn \xbq $ \xEH
\xEH
\xEP

\hline

$(wOR)$
\xEH
$(wOR)$
\xEH
$( \xbm wOR)$
\xEP

$ \xbf \xcn \xbq,$ $ \xbf ' \xcl \xbq $ $ \xch $
\xEH
$ \ol{ \ol{T} } \xcs \ol{T' }$ $ \xcc $ $ \ol{ \ol{T \xco T' } }$
\xEH
$f(X \xcv Y) \xcc f(X) \xcv Y$
\xEP

$ \xbf \xco \xbf ' \xcn \xbq $
\xEH
\xEH
\xEP

\hline

$(disjOR)$
\xEH
$(disjOR)$
\xEH
$( \xbm disjOR)$
\xEP

$ \xbf \xcl \xCN \xbf ',$ $ \xbf \xcn \xbq,$
\xEH
$\xCN Con(T \xcv T') \xch$
\xEH
$X \xcs Y= \xCQ $ $ \xch $
\xEP

$ \xbf ' \xcn \xbq $ $ \xch $ $ \xbf \xco \xbf ' \xcn \xbq $
\xEH
$ \ol{ \ol{T} } \xcs \ol{ \ol{T' } } \xcc \ol{ \ol{T \xco T' } }$
\xEH
$f(X \xcv Y) \xcc f(X) \xcv f(Y)$
\xEP

\hline

$(LLE)$
\xEH
$(LLE)$
\xEH
\xEP

Left Logical Equivalence
\xEH
\xEH
\xEP

$ \xcl \xbf \xcr \xbf ',  \xbf \xcn \xbq   \xch $
\xEH
$ \ol{T}= \ol{T' }  \xch   \ol{\ol{T}} = \ol{\ol{T'}}$
\xEH
trivially true
\xEP

$ \xbf ' \xcn \xbq $ \xEH \xEH \xEP

\hline

$(RW)$ Right Weakening
\xEH
$(RW)$
\xEH
upward closure
\xEP

$ \xbf \xcn \xbq,  \xcl \xbq \xcp \xbq '   \xch $
\xEH
$ T \xcn \xbq,  \xcl \xbq \xcp \xbq '   \xch $
\xEH
\xEP

$ \xbf \xcn \xbq ' $
\xEH
$T \xcn \xbq ' $
\xEH
\xEP

\hline

$(CCL)$ Classical Closure \xEH $(CCL)$ \xEH \xEP

\xEH
$ \ol{ \ol{T} }$ is classically
\xEH
trivially true
\xEP

\xEH closed \xEH \xEP

\hline

$(SC)$ Supraclassicality \xEH $(SC)$ \xEH $( \xbm \xcc )$ \xEP

$ \xbf \xcl \xbq $ $ \xch $ $ \xbf \xcn \xbq $ \xEH $ \ol{T} \xcc \ol{
\ol{T} }$ \xEH $f(X) \xcc X$ \xEP

\cline{1-1}

$(REF)$ Reflexivity \xEH \xEH \xEP
$ \xbD,\xba \xcn \xba $ \xEH \xEH \xEP

\hline

$(CP)$ \xEH $(CP)$ \xEH $( \xbm \xCQ )$ \xEP

Consistency Preservation \xEH \xEH \xEP

$ \xbf \xcn \xcT $ $ \xch $ $ \xbf \xcl \xcT $ \xEH $T \xcn \xcT $ $ \xch
$ $T \xcl \xcT $ \xEH $f(X)= \xCQ $ $ \xch $ $X= \xCQ $ \xEP

\hline

\xEH
\xEH $( \xbm \xCQ fin)$
\xEP

\xEH
\xEH $X \xEd \xCQ $ $ \xch $ $f(X) \xEd \xCQ $
\xEP

\xEH
\xEH for finite $X$
\xEP

\hline

\xEH $(PR)$ \xEH $( \xbm PR)$ \xEP

$ \ol{ \ol{ \xbf \xcu \xbf ' } }$ $ \xcc $ $ \ol{ \ol{ \ol{ \xbf } } \xcv
\{ \xbf ' \}}$ \xEH
$ \ol{ \ol{T \xcv T' } }$ $ \xcc $ $ \ol{ \ol{ \ol{T} } \xcv T' }$ \xEH
$X \xcc Y$ $ \xch $
\xEP

\xEH \xEH $f(Y) \xcs X \xcc f(X)$
\xEP

\cline{3-3}

\xEH
\xEH
$(\xbm PR ')$
\xEP

\xEH
\xEH
$f(X) \xcs Y \xcc f(X \xcs Y)$
\xEP

\hline

$(CUT)$ \xEH $(CUT)$ \xEH $ (\xbm CUT) $ \xEP
$ \xbD \xcn \xba; \xbD, \xba \xcn \xbb \xch $ \xEH
$T \xcc \ol{T' } \xcc \ol{ \ol{T} }  \xch $ \xEH
$f(X) \xcc Y \xcc X  \xch $ \xEP
$ \xbD \xcn \xbb $ \xEH
$ \ol{ \ol{T'} } \xcc \ol{ \ol{T} }$ \xEH
$f(X) \xcc f(Y)$
\xEP

\hline

\end{tabular}

}

{\footnotesize

\begin{tabular}{|c|c|c|}

\hline

\multicolumn{3}{|c|}{Cumulativity} \xEP

\hline

$(CM)$ Cautious Monotony \xEH $(CM)$ \xEH $ (\xbm CM) $ \xEP

$ \xbf \xcn \xbq,  \xbf \xcn \xbq '   \xch $ \xEH
$T \xcc \ol{T' } \xcc \ol{ \ol{T} }  \xch $ \xEH
$f(X) \xcc Y \xcc X  \xch $
\xEP

$ \xbf \xcu \xbq \xcn \xbq ' $ \xEH
$ \ol{ \ol{T} } \xcc \ol{ \ol{T' } }$ \xEH
$f(Y) \xcc f(X)$
\xEP

\cline{1-1}

\cline{3-3}

or $(ResM)$ Restricted Monotony \xEH \xEH $(\xbm ResM)$ \xEP
$ \xbD \xcn \xba, \xbb \xch \xbD,\xba \xcn \xbb $ \xEH \xEH
$ f(X) \xcc A \xcs B \xch f(X \xcs A) \xcc B $ \xEP

\hline

$(CUM)$ Cumulativity \xEH $(CUM)$ \xEH $( \xbm CUM)$ \xEP

$ \xbf \xcn \xbq   \xch $ \xEH
$T \xcc \ol{T' } \xcc \ol{ \ol{T} }  \xch $ \xEH
$f(X) \xcc Y \xcc X  \xch $
\xEP

$( \xbf \xcn \xbq '   \xcj   \xbf \xcu \xbq \xcn \xbq ' )$ \xEH
$ \ol{ \ol{T} }= \ol{ \ol{T' } }$ \xEH
$f(Y)=f(X)$ \xEP

\hline

\xEH
$ (\xcc \xcd) $
\xEH
$ (\xbm \xcc \xcd) $
\xEP
\xEH
$T \xcc \ol{\ol{T'}}, T' \xcc \ol{\ol{T}} \xch $
\xEH
$ f(X) \xcc Y, f(Y) \xcc X \xch $
\xEP
\xEH
$ \ol{\ol{T'}} = \ol{\ol{T}}$
\xEH
$ f(X)=f(Y) $
\xEP

\hline

\multicolumn{3}{|c|}{Rationality} \xEP

\hline

$(RatM)$ Rational Monotony \xEH $(RatM)$ \xEH $( \xbm RatM)$ \xEP

$ \xbf \xcn \xbq,  \xbf \xcN \xCN \xbq '   \xch $ \xEH
$Con(T \xcv \ol{\ol{T'}})$, $T \xcl T'$ $ \xch $ \xEH
$X \xcc Y, X \xcs f(Y) \xEd \xCQ   \xch $
\xEP

$ \xbf \xcu \xbq ' \xcn \xbq $ \xEH
$ \ol{\ol{T}} \xcd \ol{\ol{\ol{T'}} \xcv T} $ \xEH
$f(X) \xcc f(Y) \xcs X$ \xEP

\hline

\xEH $(RatM=)$ \xEH $( \xbm =)$ \xEP

\xEH
$Con(T \xcv \ol{\ol{T'}})$, $T \xcl T'$ $ \xch $ \xEH
$X \xcc Y, X \xcs f(Y) \xEd \xCQ   \xch $
\xEP

\xEH
$ \ol{\ol{T}} = \ol{\ol{\ol{T'}} \xcv T} $ \xEH
$f(X) = f(Y) \xcs X$ \xEP

\hline

\xEH
$(Log=' )$
\xEH $( \xbm =' )$
\xEP

\xEH
$Con( \ol{ \ol{T' } } \xcv T)$ $ \xch $
\xEH $f(Y) \xcs X \xEd \xCQ $ $ \xch $
\xEP

\xEH
$ \ol{ \ol{T \xcv T' } }= \ol{ \ol{ \ol{T' } } \xcv T}$
\xEH $f(Y \xcs X)=f(Y) \xcs X$
\xEP

\hline

\xEH
$(Log \xFO )$
\xEH $( \xbm \xFO )$
\xEP

\xEH
$ \ol{ \ol{T \xco T' } }$ is one of
\xEH $f(X \xcv Y)$ is one of
\xEP

\xEH
$\ol{\ol{T}},$ or $\ol{\ol{T'}},$ or $\ol{\ol{T}} \xcs \ol{\ol{T'}}$ (by (CCL))
\xEH $f(X),$ $f(Y)$ or $f(X) \xcv f(Y)$
\xEP

\hline

\xEH
$(Log \xcv )$
\xEH $( \xbm \xcv )$
\xEP

\xEH
$Con( \ol{ \ol{T' } } \xcv T),$ $ \xCN Con( \ol{ \ol{T' } }
\xcv \ol{ \ol{T} })$ $ \xch $
\xEH $f(Y) \xcs (X-f(X)) \xEd \xCQ $ $ \xch $
\xEP

\xEH
$ \xCN Con( \ol{ \ol{T \xco T' } } \xcv T' )$
\xEH $f(X \xcv Y) \xcs Y= \xCQ$
\xEP

\hline

\xEH
$(Log \xcv ' )$
\xEH $( \xbm \xcv ' )$
\xEP

\xEH
$Con( \ol{ \ol{T' } } \xcv T),$ $ \xCN Con( \ol{ \ol{T' }
} \xcv \ol{ \ol{T} })$ $ \xch $
\xEH $f(Y) \xcs (X-f(X)) \xEd \xCQ $ $ \xch $
\xEP

\xEH
$ \ol{ \ol{T \xco T' } }= \ol{ \ol{T} }$
\xEH $f(X \xcv Y)=f(X)$
\xEP

\hline

\xEH
\xEH $( \xbm \xbe )$
\xEP

\xEH
\xEH $a \xbe X-f(X)$ $ \xch $
\xEP

\xEH
\xEH $ \xcE b \xbe X.a \xce f(\{a,b\})$
\xEP

\hline

\end{tabular}

}

$(PR)$ is also called infinite conditionalization - we choose the name for
its central role for preferential structures $(PR)$ or $( \xbm PR).$

The system of rules $(AND)$ $(OR)$ $(LLE)$ $(RW)$ $(SC)$ $(CP)$ $(CM)$ $(CUM)$
is also called system $P$ (for preferential), adding $(RatM)$ gives the system
$R$ (for rationality or rankedness).

Roughly: Smooth preferential structures generate logics satisfying system
$P$, ranked structures logics satisfying system $R$.

A logic satisfying $(REF)$, $(ResM)$, and $(CUT)$ is called a consequence
relation.

$(LLE)$ and$(CCL)$ will hold automatically, whenever we work with model sets.

$(AND)$ is obviously closely related to filters, and corresponds to closure
under finite intersections. $(RW)$ corresponds to upward closure of filters.

More precisely, validity of both depend on the definition, and the
direction we consider.

Given $f$ and $(\xbm \xcc )$, $f(X)\xcc X$ generates a pricipal filter:
$\{X'\xcc X:f(X)\xcc X'\}$, with
the definition: If $X=M(T)$, then $T\xcn \xbf$  iff $f(X)\xcc M(\xbf )$.
Validity of $(AND)$ and
$(RW)$ are then trivial.

Conversely, we can define for $X=M(T)$

$\xdx:=\{X'\xcc X: \xcE \xbf (X'=X\xcs M(\xbf )$ and $T\xcn \xbf )\}$.

$(AND)$ then makes $\xdx$  closed under
finite intersections, $(RW)$ makes $\xdx$  upward
closed. This is in the infinite case usually not yet a filter, as not all
subsets of $X$ need to be definable this way.
In this case, we complete $\xdx$  by
adding all $X''$ such that there is $X'\xcc X''\xcc X$, $X'\xbe\xdx$.

Alternatively, we can define

$\xdx:=\{X'\xcc X: \xcS\{X \xcs M(\xbf ): T\xcn \xbf \} \xcc X' \}$.

$(SC)$ corresponds to the choice of a subset.

$(CP)$ is somewhat delicate, as it presupposes that the chosen model set is
non-empty. This might fail in the presence of ever better choices, without
ideal ones; the problem is addressed by the limit versions.

$(PR)$ is an infinitary version of one half of the deduction theorem: Let $T$
stand for $\xbf$, $T'$ for $\xbq$, and $\xbf \xcu \xbq \xcn \xbs$,
so $\xbf \xcn \xbq \xcp \xbs$, but $(\xbq \xcp \xbs )\xcu \xbq \xcl \xbs$.

$(CUM)$ (whose most interesting half in our context is $(CM)$) may best be seen
as
normal use of lemmas: We have worked hard and found some lemmas. Now
we can take a rest, and come back again with our new lemmas. Adding them to the
axioms will neither add new theorems, nor prevent old ones to hold.

\index{Fact Mu-Base}

\ed

\bfa

$\hspace{0.01em}$


\label{Fact Mu-Base}

This table is to be read as follows: If the left hand side holds for some
function $f: \xdy \xcp \xdp (U),$ and the auxiliary properties noted in
the middle also
hold for $f$ or $ \xdy,$ then the right hand side will hold, too - and
conversely.

{\small

\begin{tabular}{|c|c|c|c|}

\hline

\multicolumn{4}{|c|}{Basics} \xEP

\hline

(1.1)
\xEH
$(\xbm PR)$
\xEH
$\xch$ $(\xcs)+(\xbm \xcc)$
\xEH
$(\xbm PR')$
\xEP

\cline{1-1}

\cline{3-3}

(1.2)
\xEH
\xEH
$\xci$
\xEH
\xEP

\hline

(2.1)
\xEH
$(\xbm PR)$
\xEH
$\xch$ $(\xbm \xcc)$
\xEH
$(\xbm OR)$
\xEP

\cline{1-1}

\cline{3-3}

(2.2)
\xEH
\xEH
$\xci$ $(\xbm \xcc)$ + closure
\xEH
\xEP

\xEH
\xEH
under set difference
\xEH
\xEP

\hline

(3)
\xEH
$(\xbm PR)$
\xEH
$\xch$
\xEH
$( \xbm CUT)$
\xEP

\hline

(4)
\xEH
$(\xbm \xcc )+(\xbm \xcc \xcd )+(\xbm CUM)+$
\xEH
$\xcH$
\xEH
$( \xbm PR)$
\xEP

\xEH
$(\xbm RatM)+(\xcs )$
\xEH
\xEH
\xEP

\hline

\multicolumn{4}{|c|}{Cumulativity} \xEP

\hline

(5.1)
\xEH
$(\xbm CM)$
\xEH
$\xch$ $(\xcs)+(\xbm \xcc)$
\xEH
$(\xbm ResM)$
\xEP

\cline{1-1}

\cline{3-3}

(5.2)
\xEH
\xEH
$\xci$ (infin.)
\xEH
\xEP

\hline

(6)
\xEH
$(\xbm CM)+(\xbm CUT)$
\xEH
$\xcj$
\xEH
$(\xbm CUM)$
\xEP

\hline

(7)
\xEH
$( \xbm \xcc )+( \xbm \xcc \xcd )$
\xEH
$\xch$
\xEH
$( \xbm CUM)$
\xEP

\hline

(8)
\xEH
$( \xbm \xcc )+( \xbm CUM)+( \xcs )$
\xEH
$\xch$
\xEH
$( \xbm \xcc \xcd )$
\xEP

\hline

(9)
\xEH
$( \xbm \xcc )+( \xbm CUM)$
\xEH
$\xcH$
\xEH
$( \xbm \xcc \xcd )$
\xEP

\hline

\multicolumn{4}{|c|}{Rationality} \xEP

\hline

(10)
\xEH
$( \xbm RatM )+( \xbm PR )$
\xEH
$\xch$
\xEH
$( \xbm =)$
\xEP

\hline

(11)
\xEH
$( \xbm =)$
\xEH
$ \xch $
\xEH
$( \xbm PR),$
\xEP

\hline

(12.1)
\xEH
$( \xbm =)$
\xEH
$ \xch $ $(\xcs)+( \xbm \xcc )$
\xEH
$( \xbm =' ),$
\xEP
\cline{1-1}
\cline{3-3}
(12.2)
\xEH
\xEH
$ \xci $
\xEH
\xEP

\hline

(13)
\xEH
$( \xbm \xcc ),$ $( \xbm =)$
\xEH
$ \xch $ $(\xcv)$
\xEH
$( \xbm \xcv ),$
\xEP

\hline

(14)
\xEH
$( \xbm \xcc ),$ $( \xbm \xCQ ),$ $( \xbm =)$
\xEH
$ \xch $ $(\xcv)$
\xEH
$( \xbm \xFO ),$ $( \xbm \xcv ' ),$ $( \xbm CUM),$
\xEP

\hline

(15)
\xEH
$( \xbm \xcc )+( \xbm \xFO )$
\xEH
$ \xch $ $\xdy$ closed under set difference
\xEH
$( \xbm =),$
\xEP

\hline

(16)
\xEH
$( \xbm \xFO )+( \xbm \xbe )+( \xbm PR)+$
\xEH
$ \xch $ $(\xcv)$ + $\xdy$ contains singletons
\xEH
$( \xbm =),$
\xEP
\xEH
$( \xbm \xcc )$
\xEH
\xEH
\xEP

\hline

(17)
\xEH
$( \xbm CUM)+( \xbm =)$
\xEH
$ \xch $ $(\xcv)$ + $\xdy$ contains singletons
\xEH
$( \xbm \xbe ),$
\xEP

\hline

(18)
\xEH
$( \xbm CUM)+( \xbm =)+( \xbm \xcc )$
\xEH
$ \xch $ $(\xcv)$
\xEH
$( \xbm \xFO ),$
\xEP

\hline

(19)
\xEH
$( \xbm PR)+( \xbm CUM)+( \xbm \xFO )$
\xEH
$ \xch $ sufficient, e.g. true in $\xdD_{\xdl}$
\xEH
$( \xbm =)$.
\xEP

\hline

(20)
\xEH
$( \xbm \xcc )+( \xbm PR)+( \xbm =)$
\xEH
$ \xcH $
\xEH
$( \xbm \xFO ),$
\xEP

\hline

(21)
\xEH
$( \xbm \xcc )+( \xbm PR)+( \xbm \xFO )$
\xEH
$ \xcH $ (without closure
\xEH
$( \xbm =)$
\xEP
\xEH
\xEH
under set difference),
\xEH
\xEP

\hline

(22)
\xEH
$( \xbm \xcc )+( \xbm PR)+( \xbm \xFO )+$
\xEH
$ \xcH $
\xEH
$( \xbm \xbe )$
\xEP
\xEH
$( \xbm =)+( \xbm \xcv )$
\xEH
\xEH
(thus not representability
\xEP
\xEH
\xEH
\xEH
by ranked structures)
\xEP

\hline

\end{tabular}

}

\index{Fact Mu-Base Proof}

\efa

\subparagraph{
Proof
}

$\hspace{0.01em}$


All sets are to be in $ \xdy.$

(1.1) $( \xbm PR)+( \xcs )+( \xbm \xcc )$ $ \xch $ $( \xbm PR' ):$

By $X \xcs Y \xcc X$ and $( \xbm PR),$ $f(X) \xcs X \xcs Y \xcc f(X \xcs
Y).$ By $( \xbm \xcc )$ $f(X) \xcs Y=f(X) \xcs X \xcs Y.$

(1.2) $( \xbm PR' ) \xch ( \xbm PR):$

Let $X \xcc Y,$ so $X=X \xcs Y,$ so by $( \xbm PR' )$ $f(Y) \xcs X \xcc
f(X \xcs Y)=f(X).$

(2.1) $( \xbm PR)+( \xbm \xcc )$ $ \xch $ $( \xbm OR):$

$f(X \xcv Y) \xcc X \xcv Y$ by $( \xbm \xcc ),$ so $f(X \xcv Y)$ $=$ $(f(X
\xcv Y) \xcs X) \xcv (f(X \xcv Y) \xcs Y)$ $ \xcc $ $f(X) \xcv f(Y).$

(2.2) $( \xbm OR)$ $+$ $( \xbm \xcc )$ $+$ closure under set difference $
\xch $ $( \xbm PR):$

Let $X \xcc Y,$ $X':=Y-X$. $f(Y) \xcc f(X) \xcv f(X' )$ by $( \xbm OR),$
so $f(Y) \xcs X$ $ \xcc $
$(f(X) \xcs X) \xcv (f(X' ) \xcs X)$ $=_{( \xbm \xcc )}$ $f(X) \xcv \xCQ $
$=$ $f(X).$

(3) $( \xbm PR)$ $ \xch $ $( \xbm CUT):$

$f(X) \xcc Y \xcc X$ $ \xch $ $f(X) \xcc f(X) \xcs Y \xcc f(Y)$ by $( \xbm
PR).$

(4) $( \xbm \xcc )+( \xbm \xcc \xcd )+( \xbm CUM)+( \xbm RatM)+( \xcs )$ $
\xcH $ $( \xbm PR):$

This is shown in Example \ref{Example Need-Pr}.

(5.1) $( \xbm CM)+( \xcs )+( \xbm \xcc )$ $ \xch $ $( \xbm ResM):$

Let $f(X) \xcc A \xcs B,$ so $f(X) \xcc A,$ so by $( \xbm \xcc )$ $f(X)
\xcc A \xcs X \xcc X,$
so by $( \xbm CM)$ $f(A \xcs X) \xcc f(X) \xcc B.$

(5.2) $( \xbm ResM) \xch ( \xbm CM):$

We consider here the infinitary version, where all sets can be model sets
of infinite theories.
Let $f(X) \xcc Y \xcc X,$ so $f(X) \xcc Y \xcs f(X),$ so by $( \xbm ResM)$
$f(Y)=f(X \xcs Y) \xcc f(X).$

(6) $( \xbm CM)+( \xbm CUT)$ $ \xcj $ $( \xbm CUM):$

Trivial.

(7) $( \xbm \xcc )+( \xbm \xcc \xcd )$ $ \xch $ $( \xbm CUM):$

Suppose $f(D) \xcc E \xcc D.$ So by $( \xbm \xcc )$ $f(E) \xcc E \xcc D,$
so by $( \xbm \xcc \xcd )$ $f(D)=f(E).$

(8) $( \xbm \xcc )+( \xbm CUM)+( \xcs )$ $ \xch $ $( \xbm \xcc \xcd ):$

Let $f(D) \xcc E,$ $f(E) \xcc D,$ so by $( \xbm \xcc )$ $f(D) \xcc D \xcs
E \xcc D,$ $f(E) \xcc D \xcs E \xcc E.$ As $f(D \xcs E)$
is defined, so $f(D)=f(D \xcs E)=f(E)$ by $( \xbm CUM).$

(9) $( \xbm \xcc )+( \xbm CUM)$ $ \xcH $ $( \xbm \xcc \xcd ):$

This is shown in Example \ref{Example Mu-Cum-Cd}.

(10) $( \xbm RatM)+( \xbm PR)$ $ \xch $ $( \xbm =):$

Trivial.

(11) $( \xbm =)$ entails $( \xbm PR):$

Trivial.

(12.1) $( \xbm =) \xcp ( \xbm =' ):$

Let $f(Y) \xcs X \xEd \xCQ,$ we have to show $f(X \xcs Y)=f(Y) \xcs X.$
By $( \xbm \xcc )$ $f(Y) \xcc Y,$ so $f(Y) \xcs X=f(Y) \xcs (X \xcs Y),$
so by $( \xbm =)$ $f(Y) \xcs X$ $=$
$f(Y) \xcs (X \xcs Y)$ $=$ $f(X \xcs Y).$

(12.2) $( \xbm =' ) \xcp ( \xbm =):$

Let $X \xcc Y,$ $f(Y) \xcs X \xEd \xCQ,$ then $f(X)=f(Y \xcs X)=f(Y) \xcs
X.$

(13) $( \xbm \xcc ),$ $( \xbm =)$ $ \xcp $ $( \xbm \xcv ):$

If not, $f(X \xcv Y) \xcs Y \xEd \xCQ,$ but $f(Y) \xcs (X-f(X)) \xEd \xCQ
.$ By (11), $( \xbm PR)$ holds,
so $f(X \xcv Y) \xcs X \xcc f(X),$ so $ \xCQ $ $ \xEd $ $f(Y) \xcs
(X-f(X))$ $ \xcc $ $f(Y) \xcs (X-f(X \xcv Y)),$ so
$f(Y)-f(X \xcv Y) \xEd \xCQ,$ so by $( \xbm \xcc )$ $f(Y) \xcc Y$ and
$f(Y) \xEd f(X \xcv Y) \xcs Y.$
But by $( \xbm =)$ $f(Y)=f(X \xcv Y) \xcs Y,$ a contradiction.

(14)

$( \xbm \xcc ),$ $( \xbm \xCQ ),$ $( \xbm =)$ $ \xch $ $( \xbm \xFO ):$

If $X$ or $Y$ or both are empty, then this is trivial.
Assume then $X \xcv Y \xEd \xCQ,$ so by $( \xbm \xCQ )$ $f(X \xcv Y) \xEd
\xCQ.$
By $( \xbm \xcc )$ $f(X \xcv Y) \xcc X \xcv Y,$ so $f(X \xcv Y) \xcs X=
\xCQ $ and $f(X \xcv Y) \xcs Y= \xCQ $ together are
impossible.
Case 1, $f(X \xcv Y) \xcs X \xEd \xCQ $ and $f(X \xcv Y) \xcs Y \xEd \xCQ
:$ By $( \xbm =)$ $f(X \xcv Y) \xcs X=f(X)$ and
$f(X \xcv Y) \xcs Y=f(Y),$ so by $( \xbm \xcc )$ $f(X \xcv Y)=f(X) \xcv
f(Y).$
Case 2, $f(X \xcv Y) \xcs X \xEd \xCQ $ and $f(X \xcv Y) \xcs Y= \xCQ:$
So by $( \xbm =)$ $f(X \xcv Y)=f(X \xcv Y) \xcs X=f(X).$
Case 3, $f(X \xcv Y) \xcs X= \xCQ $ and $f(X \xcv Y) \xcs Y \xEd \xCQ:$
Symmetrical.

$( \xbm \xcc ),$ $( \xbm \xCQ ),$ $( \xbm =)$ $ \xch $ $( \xbm \xcv ' ):$

If $X \xcv Y= \xCQ,$ then $f(X \xcv Y)=f(X)= \xCQ $ by $( \xbm \xcc ).$
So suppose $X \xcv Y \xEd \xCQ.$ By
(13), $f(X \xcv Y) \xcs Y= \xCQ,$ so $f(X \xcv Y) \xcc X$ by $( \xbm \xcc
).$ By $( \xbm \xCQ ),$ $f(X \xcv Y) \xEd \xCQ,$ so
$f(X \xcv Y) \xcs X \xEd \xCQ,$ and $f(X \xcv Y)=f(X)$ by $( \xbm =).$

$( \xbm \xcc ),$ $( \xbm \xCQ ),$ $( \xbm =)$ $ \xch $ $( \xbm CUM):$

If $Y= \xCQ,$ this is trivial by $( \xbm \xcc ).$ If $Y \xEd \xCQ,$ then
by $( \xbm \xCQ )$ - which is
crucial here - $f(Y) \xEd \xCQ,$ so by $f(Y) \xcc X$ $f(Y) \xcs X \xEd
\xCQ,$ so by $( \xbm =)$
$f(Y)=f(Y) \xcs X=f(X).$

(15) $( \xbm \xcc )+( \xbm \xFO )$ $ \xcp $ $( \xbm =):$

Let $X \xcc Y,$ and consider $Y=X \xcv (Y-$X). Then $f(Y)=f(X) \xFO
f(Y-$X). As
$f(Y-X) \xcs X= \xCQ,$ $f(Y) \xcs X \xcc f(X).$ If $f(Y) \xcs X \xEd \xCQ
,$ then by the same argument
$f(X)$ is involved, so $f(X) \xcc f(Y).$

(16) $( \xbm \xFO )+( \xbm \xbe )+( \xbm PR)+( \xbm \xcc )$ $ \xcp $ $(
\xbm =):$

Suppose $X \xcc Y,$ $x \xbe f(Y) \xcs X,$ we have to show $f(Y) \xcs
X=f(X).$ `` $ \xcc $ '' is trivial
by $( \xbm PR).$ `` $ \xcd $ '': Assume $a \xce f(Y)$ (by $( \xbm \xcc )),$
but $a \xbe f(X).$ By $( \xbm \xbe )$ $ \xcE b \xbe Y.a \xce f(\{a,b\}).$
As $a \xbe f(X),$ by $( \xbm PR),$ $a \xbe f(\{a,x\}).$ By $( \xbm \xFO
),$ $f(\{a,b,x\})$ $=$ $f(\{a,x\}) \xFO f(\{b\}).$
As $a \xce f(\{a,b,x\}),$ $f(\{a,b,x\})$ $=$ $f(\{b\}),$ so $x \xce
f(\{a,b,x\}),$ contradicting $( \xbm PR),$
as $a,b,x \xbe Y.$

(17) $( \xbm CUM)+( \xbm =)$ $ \xcp $ $( \xbm \xbe ):$

Let $a \xbe X-f(X).$ If $f(X)= \xCQ,$ then $f(\{a\})= \xCQ $ by $( \xbm
CUM).$ If not: Let
$b \xbe f(X),$ then $a \xce f(\{a,b\})$ by $( \xbm =).$

(18) $( \xbm CUM)+( \xbm =)+( \xbm \xcc )$ $ \xcp $ $( \xbm \xFO ):$

By $( \xbm CUM),$ $f(X \xcv Y) \xcc X \xcc X \xcv Y$ $ \xcp $ $f(X)=f(X
\xcv Y),$ and $f(X \xcv Y) \xcc Y \xcc X \xcv Y$ $ \xcp $
$f(Y)=f(X \xcv Y).$ Thus, if $( \xbm \xFO )$ were to fail, $f(X \xcv Y)
\xcC X,$ $f(X \xcv Y) \xcC Y,$ but then
by $( \xbm \xcc )$ $f(X \xcv Y) \xcs X \xEd \xCQ,$ so $f(X)=f(X \xcv Y)
\xcs X,$ and $f(X \xcv Y) \xcs Y \xEd \xCQ,$ so
$f(Y)=f(X \xcv Y) \xcs Y$ by $( \xbm =).$ Thus, $f(X \xcv Y)$ $=$ $(f(X
\xcv Y) \xcs X) \xcv (f(X \xcv Y) \xcs Y)$ $=$
$f(X) \xcv f(Y).$

(19) $( \xbm PR)+( \xbm CUM)+( \xbm \xFO )$ $ \xcp $ $( \xbm =):$

Suppose $( \xbm =)$ does not hold. So, by $( \xbm PR),$ there are $X,Y,y$
s.t. $X \xcc Y,$ $X \xcs f(Y) \xEd \xCQ,$
$y \xbe Y-f(Y),$ $y \xbe f(X).$ Let $a \xbe X \xcs f(Y).$ If $f(Y)=\{a\},$
then by $( \xbm CUM)$ $f(Y)=f(X),$ so
there must be $b \xbe f(Y),$ $b \xEd a.$ Take now $Y',$ $Y'' $ s.t. $Y=Y'
\xcv Y'',$ $a \xbe Y',$ $a \xce Y'',$ $b \xbe Y'',$
$b \xce Y',$ $y \xbe Y' \xcs Y''.$ Assume now $( \xbm \xFO )$ to hold,
we show a contradiction.
If $y \xce f(Y'' ),$ then by $( \xbm PR)$ $y \xce f(Y'' \xcv \{a\}).$ But
$f(Y'' \xcv \{a\})$ $=$ $f(Y'' ) \xFO f(\{a,y\}),$
so $f(Y'' \xcv \{a\})=f(Y'' ),$ contradicting $a \xbe f(Y).$ If $y \xbe
f(Y'' ),$ then by $f(Y)$ $=$
$f(Y' ) \xFO f(Y'' ),$ $f(Y)=f(Y' ),$ $contradiction$ as $b \xce f(Y' ).$

(20) $( \xbm \xcc )+( \xbm PR)+( \xbm =)$ $ \xcH $ $( \xbm \xFO ):$

See Example \ref{Example Mu-Barbar}.

(21) $( \xbm \xcc )+( \xbm PR)+( \xbm \xFO )$ $ \xcH $ $( \xbm =):$

See Example \ref{Example Mu-Equal}.

(22) $( \xbm \xcc )+( \xbm PR)+( \xbm \xFO )+( \xbm =)+( \xbm \xcv )$ $
\xcH $ $( \xbm \xbe ):$

See Example \ref{Example Mu-Epsilon}.

Thus, by Fact \ref{Fact Rank-Hold}, the conditions do not assure
representability by ranked structures.

$ \xcz $
\\[3ex]
\index{Example Mu-Cum-Cd}

\be

$\hspace{0.01em}$


\label{Example Mu-Cum-Cd}

We show here $( \xbm \xcc )+( \xbm CUM)$ $ \xcH $ $( \xbm \xcc \xcd ).$

Consider $X:=\{a,b,c\},$ $Y:=\{a,b,d\},$ $f(X):=\{a\},$ $f(Y):=\{a,b\},$ $
\xdy:=\{X,Y\}.$
(If $f(\{a,b\})$ were defined, we would have $f(X)=f(\{a,b\})=f(Y),$
$contradiction.)$

Obviously, $( \xbm \xcc )$ and $( \xbm CUM)$ hold, but not $( \xbm \xcc
\xcd ).$

$ \xcz $
\\[3ex]
\index{Example Need-Pr}

\ee

\be

$\hspace{0.01em}$


\label{Example Need-Pr}

We show here $( \xbm \xcc )+( \xbm \xcc \xcd )+( \xbm CUM)+( \xbm RatM)+(
\xcs )$ $ \xcH $ $( \xbm PR).$

Let $U:=\{a,b,c\}.$ Let $ \xdy = \xdp (U).$ So $( \xcs )$ is trivially
satisfied.
Set $f(X):=X$ for all $X \xcc U$ except for $f(\{a,b\})=\{b\}.$ Obviously,
this cannot be represented by a preferential structure and $( \xbm PR)$ is
false
for $U$ and $\{a,b\}.$ But it satisfies $( \xbm \xcc ),$ $( \xbm CUM),$ $(
\xbm RatM).$ $( \xbm \xcc )$ is trivial.
$( \xbm CUM):$ Let $f(X) \xcc Y \xcc X.$ If $f(X)=X,$ we are done.
Consider $f(\{a,b\})=\{b\}.$ If
$\{b\} \xcc Y \xcc \{a,b\},$ then $f(Y)=\{b\},$ so we are done again. It
is shown in
Fact \ref{Fact Mu-Base}, (8) that $( \xbm \xcc \xcd )$ follows.
$( \xbm RatM):$ Suppose $X \xcc Y,$ $X \xcs f(Y) \xEd \xCQ,$ we have to
show $f(X) \xcc f(Y) \xcs X.$ If $f(Y)=Y,$ the result holds by $X \xcc Y,$
so it does if $X=Y.$
The only remaining case is $Y=\{a,b\},$ $X=\{b\},$ and the result holds
again.

$ \xcz $
\\[3ex]
\index{Example Mu-Barbar}

\ee

\be

$\hspace{0.01em}$


\label{Example Mu-Barbar}

The example shows that $( \xbm \xcc )+( \xbm PR)+( \xbm =)$ $ \xcH $ $(
\xbm \xFO ).$

Consider the following structure without transitivity:
$U:=\{a,b,c,d\},$ $c$ and $d$ have $ \xbo $ many copies in descending
order $c_{1} \xed c_{2}$  \Xl., etc.
$a,b$ have one single copy each. $a \xed b,$ $a \xed d_{1},$ $b \xed a,$
$b \xed c_{1}.$
$( \xbm \xFO )$ does not hold: $f(U)= \xCQ,$ but $f(\{a,c\})=\{a\},$
$f(\{b,d\})=\{b\}.$
$( \xbm PR)$ holds as in all preferential structures.
$( \xbm =)$ holds: If it were to fail, then for some $A \xcc B,$ $f(B)
\xcs A \xEd \xCQ,$ so $f(B) \xEd \xCQ.$
But the only possible cases for $B$ are now: $(a \xbe B,$ $b,d \xce B)$ or
$(b \xbe B,$ $a,c \xce B).$
Thus, $B$ can be $\{a\},$ $\{a,c\},$ $\{b\},$ $\{b,d\}$ with $f(B)=$
$\{a\},$ $\{a\},$ $\{b\},$ $\{b\}.$
If $A=B,$ then the result will hold trivially. Moreover, $ \xCf A$ has to
be $ \xEd \xCQ.$
So the remaining cases of $B$ where it might fail are $B=$ $\{a,c\}$ and
$\{b,d\},$ and
by $f(B) \xcs A \xEd \xCQ,$ the only cases of $ \xCf A$ where it might
fail, are $A=$ $\{a\}$ or $\{b\}$
respectively.
So the only cases remaining are: $B=\{a,c\},$ $A=\{a\}$ and $B=\{b,d\},$
$A=\{b\}.$
In the first case, $f(A)=f(B)=\{a\},$ in the second $f(A)=f(B)=\{b\},$ but
$( \xbm =)$
holds in both.

$ \xcz $
\\[3ex]
\index{Example Mu-Equal}

\ee

\be

$\hspace{0.01em}$


\label{Example Mu-Equal}

The example shows that $( \xbm \xcc )+( \xbm PR)+( \xbm \xFO )$ $ \xcH $
$( \xbm =).$

Work in the set of theory definable model sets of an infinite
propositional
language. Note that this is not closed under set difference, and closure
properties will play a crucial role in the argumentation.
Let $U:=\{y,a,x_{i< \xbo }\},$ where $x_{i} \xcp a$ in the standard
topology. For the order,
arrange s.t. $y$ is minimized by any set iff this set contains a cofinal
subsequence of
the $x_{i},$ this can be done by the standard construction. Moreover, let
the $x_{i}$
all kill themselves, i.e. with $ \xbo $ many copies $x^{1}_{i} \xed
x^{2}_{i} \xed $  \Xl. There are no other
elements in the relation. Note that if $a \xce \xbm (X),$ then $a \xce X,$
and $X$ cannot contain
a cofinal subsequence of the $x_{i},$ as $X$ is closed in the standard
topology.
(A short argument: suppose $X$ contains such a subsequence, but $a \xce
X.$ Then the
theory of a $Th(a)$ is inconsistent with $Th(X),$ so already a finite
subset of
$Th(a)$ is inconsistent with $Th(X),$ but such a finite subset will
finally hold
in a cofinal sequence converging to a.)
Likewise, if $y \xbe \xbm (X),$ then $X$ cannot contain a cofinal
subsequence of the $x_{i}.$

Obviously, $( \xbm \xcc )$ and $( \xbm PR)$ hold, but $( \xbm =)$ does not
hold: Set $B:=U,$ $A:=\{a,y\}.$
Then $ \xbm (B)=\{a\},$ $ \xbm (A)=\{a,y\},$ contradicting $( \xbm =).$

It remains to show that $( \xbm \xFO )$ holds.

$ \xbm (X)$ can only be $ \xCQ,$ $\{a\},$ $\{y\},$ $\{a,y\}.$ As $ \xbm
(A \xcv B) \xcc \xbm (A) \xcv \xbm (B)$ by $( \xbm PR),$

Case 1, $ \xbm (A \xcv B)=\{a,y\}$ is settled.

Note that if $y \xbe X- \xbm (X),$ then $X$ will contain a cofinal
subsequence, and thus
$a \xbe \xbm (X).$

Case 2: $ \xbm (A \xcv B)=\{a\}.$

Case 2.1: $ \xbm (A)=\{a\}$ - we are done.

Case 2.2: $ \xbm (A)=\{y\}:$ $ \xCf A$ does not contain $ \xCf a,$ nor a
cofinal subsequence.
If $ \xbm (B)= \xCQ,$ then $a \xce B,$ so $a \xce A \xcv B,$ a
contradiction.
If $ \xbm (B)=\{a\},$ we are done.
If $y \xbe \xbm (B),$ then $y \xbe B,$ but $B$ does not contain a cofinal
subsequence, so
$A \xcv B$ does not either, so $y \xbe \xbm (A \xcv B),$ $contradiction.$

Case 2.3: $ \xbm (A)= \xCQ:$ $ \xCf A$ cannot contain a cofinal
subsequence.
If $ \xbm (B)=\{a\},$ we are done.
$a \xbe \xbm (B)$ does have to hold, so $ \xbm (B)=\{a,y\}$ is the only
remaining possibility.
But then $B$ does not contain a cofinal subsequence, and neither does $A
\xcv B,$ so
$y \xbe \xbm (A \xcv B),$ $contradiction.$

Case 2.4: $ \xbm (A)=\{a,y\}:$ $ \xCf A$ does not contain a cofinal
subsequence.
If $ \xbm (B)=\{a\},$ we are done.
If $ \xbm (B)= \xCQ,$ $B$ does not contain a cofinal subsequence (as $a
\xce B),$ so neither
does $A \xcv B,$ so $y \xbe \xbm (A \xcv B),$ $contradiction.$
If $y \xbe \xbm (B),$ $B$ does not contain a cofinal subsequence, and we
are done again.

Case 3: $ \xbm (A \xcv B)=\{y\}:$
To obtain a contradiction, we need $a \xbe \xbm (A)$ or $a \xbe \xbm (B).$
But in both cases
$a \xbe \xbm (A \xcv B).$

Case 4: $ \xbm (A \xcv B)= \xCQ:$
Thus, $A \xcv B$ contains no cofinal subsequence. If, e.g. $y \xbe \xbm
(A),$ then $y \xbe \xbm (A \xcv B),$
if $a \xbe \xbm (A),$ then $a \xbe \xbm (A \xcv B),$ so $ \xbm (A)= \xCQ
.$

$ \xcz $
\\[3ex]
\index{Example Mu-Epsilon}

\ee

\be

$\hspace{0.01em}$


\label{Example Mu-Epsilon}

The example show that $( \xbm \xcc )+( \xbm PR)+( \xbm \xFO )+( \xbm =)+(
\xbm \xcv )$ $ \xcH $ $( \xbm \xbe ).$

Let $U:=\{y,x_{i< \xbo }\},$ $x_{i}$ a sequence, each $x_{i}$ kills
itself, $x^{1}_{i} \xed x^{2}_{i} \xed  \Xl $
and $y$ is killed by all cofinal subsequences of the $x_{i}.$ Then for any
$X \xcc U$
$ \xbm (X)= \xCQ $ or $ \xbm (X)=\{y\}.$

$( \xbm \xcc )$ and $( \xbm PR)$ hold obviously.

$( \xbm \xFO ):$ Let $A \xcv B$ be given. If $y \xce X,$ then for all $Y
\xcc X$ $ \xbm (Y)= \xCQ.$
So, if $y \xce A \xcv B,$ we are done. If $y \xbe A \xcv B,$ if $ \xbm (A
\xcv B)= \xCQ,$ one of $A,B$ must contain
a cofinal sequence, it will have $ \xbm = \xCQ.$ If not, then $ \xbm (A
\xcv B)=\{y\},$ and this will
also hold for the one $y$ is in.

$( \xbm =):$ Let $A \xcc B,$ $ \xbm (B) \xcs A \xEd \xCQ,$ show $ \xbm
(A)= \xbm (B) \xcs A.$ But now $ \xbm (B)=\{y\},$ $y \xbe A,$
so $B$ does not contain a cofinal subsequence, neither does A, so $ \xbm
(A)=\{y\}.$

$( \xbm \xcv ):$ $(A- \xbm (A)) \xcs \xbm (A' ) \xEd \xCQ,$ so $ \xbm (A'
)=\{y\},$ so $ \xbm (A \xcv A' )= \xCQ,$ as $y \xbe A- \xbm (A).$

But $( \xbm \xbe )$ does not hold: $y \xbe U- \xbm (U),$ but there is no
$x$ s.t. $y \xce \xbm (\{x,y\}).$

$ \xcz $
\\[3ex]
\index{Fact Mwor}

\ee

\bfa

$\hspace{0.01em}$


\label{Fact Mwor}

$( \xbm wOR)+( \xbm \xcc )$ $ \xch $ $f(X \xcv Y) \xcc f(X) \xcv f(Y) \xcv
(X \xcs Y)$
\index{Fact Mwor Proof}

\efa

\subparagraph{
Proof
}

$\hspace{0.01em}$


$f(X \xcv Y) \xcc f(X) \xcv Y,$ $f(X \xcv Y) \xcc X \xcv f(Y),$ so $f(X
\xcv Y)$ $ \xcc $ $(f(X) \xcv Y) \xcs (X \xcv f(Y))$ $=$
$f(X) \xcv f(Y) \xcv (X \xcs Y)$ $ \xcz $
\\[3ex]
\index{Proposition Alg-Log}

\bp

$\hspace{0.01em}$


\label{Proposition Alg-Log}

The following table is to be read as follows:

Let a logic $ \xcn $ satisfies $ \xCf (LLE)$ and $ \xCf (CCL),$ and define
a function $f: \xdD_{ \xdl } \xcp \xdD_{ \xdl }$
by $f(M(T)):=M( \ol{ \ol{T} }).$ Then $f$ is well defined, satisfies $(
\xbm dp),$ and $ \ol{ \ol{T} }=Th(f(M(T))).$

If $ \xcn $ satisfies a rule in the left hand side,
then - provided the additional properties noted in the middle for $ \xch $
hold, too -
$f$ will satisfy the property in the right hand side.

Conversely, if $f: \xdy \xcp \xdp (M_{ \xdl })$ is a function, with $
\xdD_{ \xdl } \xcc \xdy,$ and we define a logic
$ \xcn $ by $ \ol{ \ol{T} }:=Th(f(M(T))),$ then $ \xcn $ satisfies $ \xCf
(LLE)$ and $ \xCf (CCL).$
If $f$ satisfies $( \xbm dp),$ then $f(M(T))=M( \ol{ \ol{T} }).$

If $f$ satisfies a property in the right hand side,
then - provided the additional properties noted in the middle for $ \xci $
hold, too -
$ \xcn $ will satisfy the property in the left hand side.

If ``formula'' is noted in the table, this means that, if one of the
theories
(the one named the same way in Definition \ref{Definition Log-Cond})
is equivalent to a formula, we can renounce on $( \xbm dp).$

{\small

\begin{tabular}{|c|c|c|c|}

\hline

\multicolumn{4}{|c|}{Basics} \xEP

\hline

(1.1) \xEH $(OR)$ \xEH $\xch$ \xEH $(\xbm OR)$ \xEP

\cline{1-1}

\cline{3-3}

(1.2) \xEH \xEH $\xci$ \xEH \xEP

\hline

(2.1) \xEH $(disjOR)$ \xEH $\xch$ \xEH $(\xbm disjOR)$ \xEP

\cline{1-1}

\cline{3-3}

(2.2) \xEH \xEH $\xci$ \xEH \xEP

\hline

(3.1) \xEH $(wOR)$ \xEH $\xch$ \xEH $(\xbm wOR)$ \xEP

\cline{1-1}

\cline{3-3}

(3.2) \xEH \xEH $\xci$ \xEH \xEP

\hline

(4.1) \xEH $(SC)$ \xEH $\xch$ \xEH $(\xbm \xcc)$ \xEP

\cline{1-1}

\cline{3-3}

(4.2) \xEH \xEH $\xci$ \xEH \xEP

\hline

(5.1) \xEH $(CP)$ \xEH $\xch$ \xEH $(\xbm \xCQ)$ \xEP

\cline{1-1}

\cline{3-3}

(5.2) \xEH \xEH $\xci$ \xEH \xEP

\hline

(6.1) \xEH $(PR)$ \xEH $\xch$ \xEH $(\xbm PR)$ \xEP

\cline{1-1}

\cline{3-3}

(6.2) \xEH \xEH $\xci$ $(\xbm dp)+(\xbm \xcc)$ \xEH \xEP

\cline{1-1}

\cline{3-3}

(6.3) \xEH \xEH $\xcI$ without $(\xbm dp)$ \xEH \xEP

\cline{1-1}

\cline{3-3}

(6.4) \xEH \xEH $\xci$ $(\xbm \xcc)$ \xEH \xEP

\xEH \xEH $T'$ a formula \xEH \xEP

\hline

(6.5) \xEH $(PR)$ \xEH $\xci$ \xEH $(\xbm PR')$ \xEP

\xEH \xEH $T'$ a formula \xEH \xEP

\hline

(7.1) \xEH $(CUT)$ \xEH $\xch$ \xEH $(\xbm CUT)$ \xEP

\cline{1-1}

\cline{3-3}

(7.2) \xEH \xEH $\xci$ \xEH \xEP

\hline

\multicolumn{4}{|c|}{Cumulativity} \xEP

\hline

(8.1) \xEH $(CM)$ \xEH $\xch$ \xEH $(\xbm CM)$ \xEP

\cline{1-1}

\cline{3-3}

(8.2) \xEH \xEH $\xci$ \xEH \xEP

\hline

(9.1) \xEH $(ResM)$ \xEH $\xch$ \xEH $(\xbm ResM)$ \xEP

\cline{1-1}

\cline{3-3}

(9.2) \xEH \xEH $\xci$ \xEH \xEP

\hline

(10.1) \xEH $(\xcc \xcd)$ \xEH $\xch$ \xEH $(\xbm \xcc \xcd)$ \xEP

\cline{1-1}

\cline{3-3}

(10.2) \xEH \xEH $\xci$ \xEH \xEP

\hline

(11.1) \xEH $(CUM)$ \xEH $\xch$ \xEH $(\xbm CUM)$ \xEP

\cline{1-1}

\cline{3-3}

(11.2) \xEH \xEH $\xci$ \xEH \xEP

\hline

\multicolumn{4}{|c|}{Rationality} \xEP

\hline

(12.1) \xEH $(RatM)$ \xEH $\xch$ \xEH $(\xbm RatM)$ \xEP

\cline{1-1}

\cline{3-3}

(12.2) \xEH \xEH $\xci$ $(\xbm dp)$ \xEH \xEP

\cline{1-1}

\cline{3-3}

(12.3) \xEH \xEH $\xcI$ without $(\xbm dp)$ \xEH \xEP

\cline{1-1}

\cline{3-3}

(12.4) \xEH \xEH $\xci$ \xEH \xEP

\xEH \xEH $T$ a formula \xEH \xEP

\hline

(13.1) \xEH $(RatM=)$ \xEH $\xch$ \xEH $(\xbm =)$ \xEP

\cline{1-1}

\cline{3-3}

(13.2) \xEH \xEH $\xci$ $(\xbm dp)$ \xEH \xEP

\cline{1-1}

\cline{3-3}

(13.3) \xEH \xEH $\xcI$ without $(\xbm dp)$ \xEH \xEP

\cline{1-1}

\cline{3-3}

(13.4) \xEH \xEH $\xci$ \xEH \xEP

\xEH \xEH $T$ a formula \xEH \xEP

\hline

(14.1) \xEH $(Log = ')$ \xEH $\xch$ \xEH $(\xbm = ')$ \xEP

\cline{1-1}
\cline{3-3}

(14.2) \xEH \xEH $\xci$ $(\xbm dp)$ \xEH \xEP

\cline{1-1}
\cline{3-3}

(14.3) \xEH \xEH $\xcI$ without $(\xbm dp)$ \xEH \xEP

\cline{1-1}
\cline{3-3}

(14.4) \xEH \xEH $\xci$ $T$ a formula \xEH \xEP

\hline

(15.1) \xEH $(Log \xFO )$ \xEH $\xch$ \xEH $(\xbm \xFO )$ \xEP

\cline{1-1}
\cline{3-3}

(15.2) \xEH \xEH $\xci$ \xEH \xEP

\hline

(16.1)
\xEH
$(Log \xcv )$
\xEH
$\xch$ $(\xbm \xcc)+(\xbm =)$
\xEH
$(\xbm \xcv )$
\xEP

\cline{1-1}
\cline{3-3}

(16.2) \xEH \xEH $\xci$ $(\xbm dp)$ \xEH \xEP

\cline{1-1}
\cline{3-3}

(16.3) \xEH \xEH $\xcI$ without $(\xbm dp)$ \xEH \xEP

\hline

(17.1)
\xEH
$(Log \xcv ')$
\xEH
$\xch$ $(\xbm \xcc)+(\xbm =)$
\xEH
$(\xbm \xcv ')$
\xEP

\cline{1-1}
\cline{3-3}

(17.2) \xEH \xEH $\xci$ $(\xbm dp)$ \xEH \xEP

\cline{1-1}
\cline{3-3}

(17.3) \xEH \xEH $\xcI$ without $(\xbm dp)$ \xEH \xEP

\hline

\end{tabular}

}

\index{Proposition Alg-Log Proof}

\ep

\subparagraph{
Proof
}

$\hspace{0.01em}$


Set $f(T):=f(M(T)),$ note that $f(T \xcv T' ):=f(M(T \xcv T' ))=f(M(T)
\xcs M(T' )).$

We show first the general framework.

Let $ \xcn $ satisfy $ \xCf (LLE)$ and $ \xCf (CCL).$ Let $f: \xdD_{ \xdl
} \xcp \xdD_{ \xdl }$ be defined by $f(M(T)):=M( \ol{ \ol{T} }).$
If $M(T)=M(T' ),$ then $ \ol{T}= \ol{T' },$ so by $ \xCf (LLE)$ $ \ol{
\ol{T} }= \ol{ \ol{T' } },$ so $f(M(T))=f(M(T' )),$ so $f$ is
well defined and satisfies $( \xbm dp).$ By $ \xCf (CCL)$ $Th(M( \ol{
\ol{T} }))= \ol{ \ol{T} }.$

Let $f$ be given, and $ \xcn $ be defined by $ \ol{ \ol{T}
}:=Th(f(M(T))).$ Obviously, $ \xcn $ satisfies
$ \xCf (LLE)$ and $ \xCf (CCL)$ (and thus $ \xCf (RW)).$ If $f$ satisfies
$( \xbm dp),$ then $f(M(T))=M(T' )$ for
some $T',$ and $f(M(T))=M(Th(f(M(T))))=M( \ol{ \ol{T} })$ by Fact \ref{Fact
Dp-Base}. (We will use
Fact \ref{Fact Dp-Base} now without further mentioning.)

Next we show the following fact:

(a) If $f$ satisfies $( \xbm dp),$ or $T' $ is equivalent to a formula,
then
$Th(f(T) \xcs M(T' ))= \ol{ \ol{ \ol{T} } \xcv T' }.$

Case 1, $f$ satisfies $( \xbm dp).$ $Th(f(M(T)) \xcs M(T' ))$ $=$ $Th(M(
\ol{ \ol{T} }) \xcs M(T' )$ $=$ $ \ol{ \ol{ \ol{T} } \xcv T' }$
by Fact \ref{Fact Log-Form} (5).

Case 2, $T' $ is equivalent to $ \xbf '.$ $Th(f(M(T)) \xcs M( \xbf ' ))$
$=$ $ \ol{Th(f(M(T))) \xcv \{ \xbf ' \}}$ $=$
$ \ol{ \ol{ \ol{T} } \xcv \{ \xbf ' \}}$ by Fact \ref{Fact Log-Form}
(3).

We now prove the individual properties.

(1.1) $ \xCf (OR)$ $ \xch $ $( \xbm OR)$

Let $X=M(T),$ $Y=M(T' ).$ $f(X \xcv Y)$ $=$ $f(M(T) \xcv M(T' ))$ $=$
$f(M(T \xco T' ))$ $:=$ $M( \ol{ \ol{T \xco T' } })$ $ \xcc_{(OR)}$
$M( \ol{ \ol{T} } \xcs \ol{ \ol{T' } })$ $=_{(CCL)}$ $M( \ol{ \ol{T} })
\xcv M( \ol{ \ol{T' } })$ $=:$ $f(X) \xcv f(Y).$

(1.2) $( \xbm OR)$ $ \xch $ $ \xCf (OR)$

$ \ol{ \ol{T \xco T' } }$ $:=$ $Th(f(M(T \xco T' )))$ $=$ $Th(f(M(T) \xcv
M(T' )))$ $ \xcd_{( \xbm OR)}$ $Th(f(M(T)) \xcv f(M(T' )))$ $=$
(by Fact \ref{Fact Th-Union}) $Th(f(M(T))) \xcs Th(f(M(T' )))$ $=:$ $
\ol{ \ol{T} } \xcs \ol{ \ol{T' } }.$

(2) By $ \xCN Con(T,T' ) \xcj M(T) \xcs M(T' )= \xCQ,$ we can use
directly the proofs for 1.

(3.1) $ \xCf (wOR)$ $ \xch $ $( \xbm wOR)$

Let $X=M(T),$ $Y=M(T' ).$ $f(X \xcv Y)$ $=$ $f(M(T) \xcv M(T' ))$ $=$
$f(M(T \xco T' ))$ $:=$ $M( \ol{ \ol{T \xco T' } })$ $ \xcc_{(wOR)}$
$M( \ol{ \ol{T} } \xcs \ol{T' })$ $=_{(CCL)}$ $M( \ol{ \ol{T} }) \xcv M(
\ol{T' })$ $=:$ $f(X) \xcv Y.$

(3.2) $( \xbm wOR)$ $ \xch $ $ \xCf (wOR)$

$ \ol{ \ol{T \xco T' } }$ $:=$ $Th(f(M(T \xco T' )))$ $=$ $Th(f(M(T) \xcv
M(T' )))$ $ \xcd_{( \xbm wOR)}$ $Th(f(M(T)) \xcv M(T' ))$ $=$
(by Fact \ref{Fact Th-Union}) $Th(f(M(T))) \xcs Th(M(T' ))$ $=:$ $
\ol{ \ol{T} } \xcs \ol{T' }.$

(4.1) $ \xCf (SC)$ $ \xch $ $( \xbm \xcc )$

Trivial.

(4.2) $( \xbm \xcc )$ $ \xch $ $ \xCf (SC)$

Trivial.

(5.1) $ \xCf (CP)$ $ \xch $ $( \xbm \xCQ )$

Trivial.

(5.2) $( \xbm \xCQ )$ $ \xch $ $ \xCf (CP)$

Trivial.

(6.1) $ \xCf (PR)$ $ \xch $ $( \xbm PR):$

Suppose $X:=M(T),$ $Y:=M(T' ),$ $X \xcc Y,$ we have to show $f(Y) \xcs X
\xcc f(X).$
By prerequisite, $ \ol{T' } \xcc \ol{T},$ so $ \ol{T \xcv T' }= \ol{T},$
so $ \ol{ \ol{T \xcv T' } }= \ol{ \ol{T} }$ by $ \xCf (LLE).$ By $ \xCf
(PR)$ $ \ol{ \ol{T \xcv T' } } \xcc \ol{ \ol{ \ol{T'
} } \xcv T},$
so $f(Y) \xcs X=f(T' ) \xcs M(T)=M( \ol{ \ol{T' } } \xcv T) \xcc M( \ol{
\ol{T \xcv T' } })=M( \ol{ \ol{T} })=f(X).$

(6.2) $( \xbm PR)+( \xbm dp)+( \xbm \xcc )$ $ \xch $ $ \xCf (PR):$

$f(T) \xcs M(T' )=_{( \xbm \xcc )}f(T) \xcs M(T) \xcs M(T' )=f(T) \xcs M(T
\xcv T' ) \xcc_{( \xbm PR)}f(T \xcv T' ),$ so
$ \ol{ \ol{T \xcv T' } }=Th(f(T \xcv T' )) \xcc Th(f(T) \xcs M(T' ))= \ol{
\ol{ \ol{T} } \xcv T' }$ by (a) above and $( \xbm dp).$

(6.3) $( \xbm PR)$ $ \xcH $ $ \xCf (PR)$ without $( \xbm dp):$

$( \xbm PR)$ holds in all preferential structures
(see Definition \ref{Definition Pref-Str})
by Fact \ref{Fact Pref-Sound}.
Example \ref{Example Pref-Dp} shows that $ \xCf (DP)$ may fail in the
resulting logic.

(6.4) $( \xbm PR)+( \xbm \xcc )$ $ \xch $ $ \xCf (PR)$ if $T' $ is
classically equivalent to a formula:

It was shown in the proof of (6.2) that $f(T) \xcs M( \xbf ' ) \xcc f(T
\xcv \{ \xbf ' \}),$ so
$ \ol{ \ol{T \xcv \{ \xbf ' \}} }=Th(f(T \xcv \{ \xbf ' \})) \xcc Th(f(T)
\xcs M( \xbf ' ))= \ol{ \ol{ \ol{T} } \xcv \{ \xbf ' \}}$ by (a) above.

(6.5) $( \xbm PR' )$ $ \xch $ $ \xCf (PR),$ if $T' $ is classically
equivalent to a formula:

$f(M(T)) \xcs M( \xbf ' )$ $ \xcc_{( \xbm PR' )}$ $f(M(T) \xcs M( \xbf '
))$ $=$ $f(M(T \xcv \{ \xbf ' \})).$ So again
$ \ol{ \ol{T \xcv \{ \xbf ' \}} }=Th(f(T \xcv \{ \xbf ' \})) \xcc Th(f(T)
\xcs M( \xbf ' ))= \ol{ \ol{ \ol{T} } \xcv \{ \xbf ' \}}$ by (a) above.

(7.1) $ \xCf (CUT)$ $ \xch $ $( \xbm CUT)$

So let $X=M(T),$ $Y=M(T' ),$ and $f(T):=M( \ol{ \ol{T} }) \xcc M(T' ) \xcc
M(T)$ $ \xcp $ $ \ol{T} \xcc \ol{T' } \xcc \ol{ \ol{T} }=_{ \xCf (LLE)}
\ol{ \ol{( \ol{T})} }$ $ \xcp $
(by $ \xCf (CUT))$
$ \ol{ \ol{T} }= \ol{ \ol{( \ol{T})} } \xcd \ol{ \ol{( \ol{T' })} }= \ol{
\ol{T' } }$ $ \xcp $ $f(T)=M( \ol{ \ol{T} }) \xcc M( \ol{ \ol{T' } })=f(T'
),$ $thus$ $f(X) \xcc f(Y).$

(7.2) $( \xbm CUT)$ $ \xch $ $ \xCf (CUT)$

Let $T$ $ \xcc $ $ \ol{T' }$ $ \xcc $ $ \ol{ \ol{T} }.$ Thus $f(T) \xcc M(
\ol{ \ol{T} })$ $ \xcc $ $M(T' )$ $ \xcc $ $M(T),$ so by $( \xbm CUT)$
$f(T) \xcc f(T' ),$
so $ \ol{ \ol{T} }$ $=$ $Th(f(T))$ $ \xcd $ $Th(f(T' ))$ $=$ $ \ol{ \ol{T'
} }.$

(8.1) $ \xCf (CM)$ $ \xch $ $( \xbm CM)$

So let $X=M(T),$ $Y=M(T' ),$ and $f(T):=M( \ol{ \ol{T} }) \xcc M(T' ) \xcc
M(T)$ $ \xcp $ $ \ol{T} \xcc \ol{T' } \xcc \ol{ \ol{T} }=_{ \xCf (LLE)}
\ol{ \ol{( \ol{T})} }$ $ \xcp $
(by $ \xCf (LLE),$ $ \xCf (CM))$
$ \ol{ \ol{T} }= \ol{ \ol{( \ol{T})} } \xcc \ol{ \ol{( \ol{T' })} }= \ol{
\ol{T' } }$ $ \xcp $ $f(T)=M( \ol{ \ol{T} }) \xcd M( \ol{ \ol{T' } })=f(T'
),$ $thus$ $f(X) \xcd f(Y).$

(8.2) $( \xbm CM)$ $ \xch $ $ \xCf (CM)$

Let $T$ $ \xcc $ $ \ol{T' }$ $ \xcc $ $ \ol{ \ol{T} }.$ Thus by $( \xbm
CM)$ and $f(T) \xcc M( \ol{ \ol{T} })$ $ \xcc $ $M(T' )$ $ \xcc $ $M(T),$
so $f(T) \xcd f(T' )$
by $( \xbm CM),$ so $ \ol{ \ol{T} }$ $=$ $Th(f(T))$ $ \xcc $ $Th(f(T' ))$
$=$ $ \ol{ \ol{T' } }.$

(9.1) $ \xCf (ResM)$ $ \xch $ $( \xbm ResM)$

Let $f(X):=M( \ol{ \ol{ \xbD } }),$ $A:=M( \xba ),$ $B:=M( \xbb ).$ So
$f(X) \xcc A \xcs B$ $ \xch $ $ \xbD \xcn \xba, \xbb $ $ \xch_{ \xCf
(ResM)}$
$ \xbD, \xba \xcn \xbb $ $ \xch $ $M( \ol{ \ol{ \xbD, \xba } }) \xcc M(
\xbb )$ $ \xch $ $f(X \xcs A) \xcc B.$

(9.2) $( \xbm ResM)$ $ \xch $ $ \xCf (ResM)$

Let $f(X):=M( \ol{ \ol{ \xbD } }),$ $A:=M( \xba ),$ $B:=M( \xbb ).$ So $
\xbD \xcn \xba, \xbb $ $ \xch $ $f(X) \xcc A \xcs B$ $ \xch_{( \xbm
ResM)}$
$f(X \xcs A) \xcc B$ $ \xch $ $ \xbD, \xba \xcn \xbb.$

(10.1) $( \xcc \xcd )$ $ \xch $ $( \xbm \xcc \xcd )$

Let $f(T) \xcc M(T' ),$ $f(T' ) \xcc M(T).$
So $Th(M(T' )) \xcc Th(f(T)),$ $Th(M(T)) \xcc Th(f(T' )),$ so $T' \xcc
\ol{T' } \xcc \ol{ \ol{T} },$ $T \xcc \ol{T} \xcc \ol{ \ol{T' } },$
so by $( \xcc \xcd )$ $ \ol{ \ol{T} }= \ol{ \ol{T' } },$ so $f(T):=M( \ol{
\ol{T} })=M( \ol{ \ol{T' } })=:f(T' ).$

(10.2) $( \xbm \xcc \xcd )$ $ \xch $ $( \xcc \xcd )$

Let $T \xcc \ol{ \ol{T' } }$ and $T' \xcc \ol{ \ol{T} }.$ So by $ \xCf
(CCL)$ $Th(M(T))= \ol{T} \xcc \ol{ \ol{T' } }=Th(f(T' )).$
But $Th(M(T)) \xcc Th(X) \xch X \xcc M(T):$ $X \xcc M(Th(X)) \xcc
M(Th(M(T)))=M(T).$
So $f(T' ) \xcc M(T),$ likewise $f(T) \xcc M(T' ),$ so by $( \xbm \xcc
\xcd )$ $f(T)=f(T' ),$ so $ \ol{ \ol{T} }= \ol{ \ol{T' } }.$

(11.1) $ \xCf (CUM)$ $ \xch $ $( \xbm CUM):$

So let $X=M(T),$ $Y=M(T' ),$ and $f(T):=M( \ol{ \ol{T} }) \xcc M(T' ) \xcc
M(T)$ $ \xcp $ $ \ol{T} \xcc \ol{T' } \xcc \ol{ \ol{T} }=_{ \xCf (LLE)}
\ol{ \ol{( \ol{T})} }$ $ \xcp $
$ \ol{ \ol{T} }= \ol{ \ol{( \ol{T})} }= \ol{ \ol{( \ol{T' })} }= \ol{
\ol{T' } }$ $ \xcp $ $f(T)=M( \ol{ \ol{T} })=M( \ol{ \ol{T' } })=f(T' ),$
$thus$ $f(X)=f(Y).$

(11.2) $( \xbm CUM)$ $ \xch $ $ \xCf (CUM)$:

Let $T$ $ \xcc $ $ \ol{T' }$ $ \xcc $ $ \ol{ \ol{T} }.$ Thus by $( \xbm
CUM)$ and $f(T) \xcc M( \ol{ \ol{T} })$ $ \xcc $ $M(T' )$ $ \xcc $ $M(T),$
so $f(T)=f(T' ),$
so $ \ol{ \ol{T} }$ $=$ $Th(f(T))$ $=$ $Th(f(T' ))$ $=$ $ \ol{ \ol{T' }
}.$

(12.1) $ \xCf (RatM)$ $ \xch $ $( \xbm RatM)$

Let $X=M(T),$ $Y=M(T' ),$ and $X \xcc Y,$ $X \xcs f(Y) \xEd \xCQ,$ so $T
\xcl T' $ and $M(T) \xcs f(M(T' )) \xEd \xCQ,$
so $Con(T, \ol{ \ol{T' } }),$ so $ \ol{ \ol{ \ol{T' } } \xcv T} \xcc \ol{
\ol{T} }$ by $ \xCf (RatM),$ so $f(X)=f(M(T))=M( \ol{ \ol{T} }) \xcc M(
\ol{ \ol{T' } } \xcv T)=$
$M( \ol{ \ol{T' } }) \xcs M(T)=f(Y) \xcs X.$

(12.2) $( \xbm RatM)+( \xbm dp)$ $ \xch $ $ \xCf (RatM):$

Let $X=M(T),$ $Y=M(T' ),$ $T \xcl T',$ $Con(T, \ol{ \ol{T' } }),$ so $X
\xcc Y$ and by $( \xbm dp)$ $X \xcs f(Y) \xEd \xCQ,$
so by $( \xbm RatM)$ $f(X) \xcc f(Y) \xcs X,$ so
$ \ol{ \ol{T \xcv T' } }=Th(f(T \xcv T' )) \xcc Th(f(T) \xcs M(T' ))= \ol{
\ol{ \ol{T} } \xcv T' }$ by (a) above and $( \xbm dp).$

(12.3) $( \xbm RatM)$ $ \xcH $ $ \xCf (RatM)$ without $( \xbm dp):$

$( \xbm RatM)$ holds in all ranked preferential structures
(see Definition \ref{Definition Rank-Rel})
by Fact \ref{Fact Rank-Hold}. Example \ref{Example Rank-Dp} (2)
shows
that $ \xCf (RatM)$ may fail in the resulting logic.

(12.4) $( \xbm RatM)$ $ \xch $ $ \xCf (RatM)$ if $T$ is classically
equivalent to a formula:

$ \xbf \xcl T' $ $ \xch $ $M( \xbf ) \xcc M(T' ).$ $Con( \xbf, \ol{
\ol{T' } })$ $ \xcj $ $M( \ol{ \ol{T' } }) \xcs M( \xbf ) \xEd \xCQ $ $
\xcj $ $f(T' ) \xcs M( \xbf )= \xCQ $
by Fact \ref{Fact Log-Form} (4). Thus $f(M( \xbf )) \xcc f(M(T' ))
\xcs M( \xbf )$ by $( \xbm RatM).$
Thus by (a) above $ \ol{ \ol{ \ol{T' } } \xcv \{ \xbf \}} \xcc \ol{ \ol{
\xbf } }.$

(13.1) $ \xCf (RatM=)$ $ \xch $ $( \xbm =)$

Let $X=M(T),$ $Y=M(T' ),$ and $X \xcc Y,$ $X \xcs f(Y) \xEd \xCQ,$ so $T
\xcl T' $ and $M(T) \xcs f(M(T' )) \xEd \xCQ,$
so $Con(T, \ol{ \ol{T' } }),$ so $ \ol{ \ol{ \ol{T' } } \xcv T}= \ol{
\ol{T} }$ by $ \xCf (RatM=),$ so $f(X)=f(M(T))=M( \ol{ \ol{T} })=M( \ol{
\ol{T' } } \xcv T)=$
$M( \ol{ \ol{T' } }) \xcs M(T)=f(Y) \xcs X.$

(13.2) $( \xbm =)+( \xbm dp)$ $ \xch $ $ \xCf (RatM=)$

Let $X=M(T),$ $Y=M(T' ),$ $T \xcl T',$ $Con(T, \ol{ \ol{T' } }),$ so $X
\xcc Y$ and by $( \xbm dp)$ $X \xcs f(Y) \xEd \xCQ,$
so by $( \xbm =)$ $f(X)=f(Y) \xcs X.$ So $ \ol{ \ol{ \ol{T' } } \xcv T}=
\ol{ \ol{T} }$ (a) above and $( \xbm dp).$

(13.3) $( \xbm =)$ $ \xcH $ $ \xCf (RatM=)$ without $( \xbm dp):$

$( \xbm =)$ holds in all ranked preferential structures
(see Definition \ref{Definition Rank-Rel})
by Fact \ref{Fact Rank-Hold}. Example \ref{Example Rank-Dp} (1)
shows
that $ \xCf (RatM=)$ may fail in the resulting logic.

(13.4) $( \xbm =)$ $ \xch $ $ \xCf (RatM=)$ if $T$ is classically
equivalent to a formula:

The proof is almost identical to the one for (12.4). Again, the
prerequisites of $( \xbm =)$ are satisfied, so $f(M( \xbf ))=f(M(T' ))
\xcs M( \xbf ).$
Thus, $ \ol{ \ol{ \ol{T' } } \xcv \{ \xbf \}}= \ol{ \ol{ \xbf } }$ by (a)
above.

Of the last four, we show (14), (15), (17), the proof for (16) is similar
to the
one for (17).

(14.1) $(Log=' )$ $ \xch $ $( \xbm =' ):$

$f(M(T' )) \xcs M(T) \xEd \xCQ $ $ \xch $ $Con( \ol{ \ol{T' } } \xcv T)$ $
\xch_{(Log=' )}$ $ \ol{ \ol{T \xcv T' } }= \ol{ \ol{ \ol{T' } } \xcv T}$ $
\xch $
$f(M(T \xcv T' ))=f(M(T' )) \xcs M(T).$

(14.2) $( \xbm =' )+( \xbm dp)$ $ \xch $ $(Log=' ):$

$Con( \ol{ \ol{T' } } \xcv T)$ $ \xch_{( \xbm dp)}$ $f(M(T' )) \xcs M(T)
\xEd \xCQ $ $ \xch $ $f(M(T' \xcv T))=f(M(T' ) \xcs M(T))$ $=_{( \xbm ='
)}$
$f(M(T' )) \xcs M(T),$ so $ \ol{ \ol{T' \xcv T} }$ $=$ $ \ol{ \ol{ \ol{T'
} } \xcv T}$ by (a) above and $( \xbm dp).$

(14.3) $( \xbm =' )$ $ \xcH $ $(Log=' )$ without $( \xbm dp):$

By Fact \ref{Fact Rank-Hold} $( \xbm =' )$ holds in ranked
structures.
Consider Example \ref{Example Rank-Dp} (2). There, $Con(T, \ol{ \ol{T' }
}),$ $T=T \xcv T',$ and
it was shown that $ \ol{ \ol{ \ol{T' } } \xcv T}$ $ \xcC $ $ \ol{ \ol{T}
}$ $=$ $ \ol{ \ol{T \xcv T' } }$

(14.4) $( \xbm =' )$ $ \xch $ $(Log=' )$ if $T$ is classically equivalent
to a formula:

$Con( \ol{ \ol{T' } } \xcv \{ \xbf \})$ $ \xch $ $ \xCQ \xEd M( \ol{
\ol{T' } }) \xcs M( \xbf )$ $ \xch $ $f(T' ) \xcs M( \xbf ) \xEd \xCQ $ by
Fact \ref{Fact Log-Form} (4). So $f(M(T' \xcv \{ \xbf \}))$ $=$
$f(M(T' ) \xcs M( \xbf ))$ $=$
$f(M(T' )) \xcs M( \xbf )$ by $( \xbm =' ),$ so $ \ol{ \ol{T' \xcv \{ \xbf
\}} }= \ol{ \ol{ \ol{T' } } \xcv \{ \xbf \}}$ by (a) above.

(15.1) $(Log \xFO )$ $ \xch $ $( \xbm \xFO ):$

Trivial.

(15.2) $( \xbm \xFO )$ $ \xch $ $(Log \xFO ):$

Trivial.

(16) $(Log \xcv )$ $ \xcj $ $( \xbm \xcv ):$ Analogous to the proof of
(17).

(17.1) $(Log \xcv ' )+( \xbm \xcc )+( \xbm =)$ $ \xch $ $( \xbm \xcv ' ):$

$f(M(T' )) \xcs (M(T)-f(M(T))) \xEd \xCQ $ $ \xch $ (by $( \xbm \xcc ),$
$( \xbm =),$
Fact \ref{Fact Rank-Auxil}) $f(M(T' )) \xcs M(T) \xEd \xCQ,$ $f(M(T'
)) \xcs f(M(T))= \xCQ $ $ \xch $
$Con( \ol{ \ol{T' } },T),$ $ \xCN Con( \ol{ \ol{T' } }, \ol{ \ol{T} })$ $
\xch $ $ \ol{ \ol{T \xco T' } }= \ol{ \ol{T} }$ $ \xch $ $f(M(T))=f(M(T
\xco T' ))=f(M(T) \xcv M(T' )).$

(17.2) $( \xbm \xcv ' )+( \xbm dp)$ $ \xch $ $(Log \xcv ' ):$

$Con( \ol{ \ol{T' } } \xcv T),$ $ \xCN Con( \ol{ \ol{T' } } \xcv \ol{
\ol{T} })$ $ \xch_{( \xbm dp)}$ $f(T' ) \xcs M(T) \xEd \xCQ,$ $f(T' )
\xcs f(T)= \xCQ $ $ \xch $
$f(M(T' )) \xcs (M(T)-f(M(T))) \xEd \xCQ $ $ \xch $ $f(M(T))$ $=$ $f(M(T)
\xcv M(T' ))$ $=$ $f(M(T \xco T' )).$
So $ \ol{ \ol{T} }= \ol{ \ol{T \xco T' } }.$

(17.3) and (16.3) are solved by Example \ref{Example Rank-Dp} (3).

$ \xcz $
\\[3ex]
\index{Example Rank-Dp}

\be

$\hspace{0.01em}$


\label{Example Rank-Dp}

(1) $( \xbm =)$ without $( \xbm dp)$ does not imply $(RatM=):$

Take $\{p_{i}:i \xbe \xbo \}$ and put $m:=m_{ \xcU p_{i}},$ the model
which makes all $p_{i}$ true, in the top
layer, all the other in the bottom layer. Let $m' \xEd m,$ $T':= \xCQ,$
$T:=Th(m,m' ).$ Then
Then $ \ol{ \ol{T' } }=T',$ so $Con( \ol{ \ol{T' } },T),$ $ \ol{ \ol{T}
}=Th(m' ),$ $ \ol{ \ol{ \ol{T' } } \xcv T}=T.$

So $(RatM=)$ fails, but $( \xbm =)$ holds in all ranked structures.

(2) $( \xbm RatM)$ without $( \xbm dp)$ does not imply (RatM):

Take $\{p_{i}:i \xbe \xbo \}$ and let $m:=m_{ \xcU p_{i}},$ the model
which makes all $p_{i}$ true.

Let $X:=M( \xCN p_{0}) \xcv \{m\}$ be the top layer, put the rest of $M_{
\xdl }$ in the bottom layer.
Let $Y:=M_{ \xdl }.$ The structure is ranked, as shown in Fact \ref{Fact
Rank-Hold},
$( \xbm RatM)$ holds.

Let $T':= \xCQ,$ $T:=Th(X).$ We have to show that $Con(T, \ol{ \ol{T' }
}),$ $T \xcl T',$ but
$ \ol{ \ol{ \ol{T' } } \xcv T} \xcC \ol{ \ol{T} }.$ $ \ol{ \ol{T' } }$ $=$
$Th(M(p_{0})-\{m\})$ $=$ $ \ol{p_{0}}.$ $T$ $=$ $ \ol{\{ \xCN p_{0}\} \xco
Th(m)},$ $ \ol{ \ol{T} }=T.$ $So$ $Con(T, \ol{
\ol{T' } }).$
$M( \ol{ \ol{T' } })=M(p_{0}),$ $M(T)=X,$ $M( \ol{ \ol{T' } } \xcv T)=M(
\ol{ \ol{T' } }) \xcs M(T)=\{m\},$ $m \xcm p_{1},$ so $p_{1} \xbe \ol{
\ol{ \ol{T' } } \xcv T},$ but
$X \xcM p_{1}.$

(3) This example shows that we need $( \xbm dp)$ to go from $( \xbm \xcv
)$ to
$(Log \xcv )$ and from $( \xbm \xcv ' )$ to $(Log \xcv ' ).$

Let $v( \xdl ):=\{p,q\} \xcv \{p_{i}:i< \xbo \}.$ Let $m$ make all
variables true.

Put all models of $ \xCN p,$ and $m,$ in the upper layer, all other models
in the
lower layer. This is ranked, so by Fact \ref{Fact Rank-Hold} $( \xbm
\xcv )$
and $( \xbm \xcv ' )$ hold.
Set $X:=M( \xCN q) \xcv \{m\},$ $X':=M(q)-\{m\},$ $T:=Th(X)= \xCN q \xco
Th(m),$ $T':=Th(X' )= \ol{q}.$
Then $ \ol{ \ol{T} }= \ol{p \xcu \xCN q},$ $ \ol{ \ol{T' } }= \ol{p \xcu
q}.$ We have $Con( \ol{ \ol{T' } },T),$ $ \xCN Con( \ol{ \ol{T' } }, \ol{
\ol{T} }).$
But $ \ol{ \ol{T \xco T' } }= \ol{p} \xEd \ol{ \ol{T} }= \ol{p \xcu \xCN
q}$ and $Con( \ol{ \ol{T \xco T' } },T' ),$ so $(Log \xcv )$ and $(Log
\xcv ' )$ fail.

$ \xcz $
\\[3ex]
\index{Fact Cut-Pr}

\ee

\bfa

$\hspace{0.01em}$


\label{Fact Cut-Pr}

$ \xCf (CUT)$ $ \xcH $ $ \xCf (PR)$
\index{Fact Cut-Pr Proof}

\efa

\subparagraph{
Proof
}

$\hspace{0.01em}$


We give two proofs:

(1) If $ \xCf (CUT)$ $ \xch $ $ \xCf (PR),$ then by $( \xbm PR)$ $ \xch $
(by Fact \ref{Fact Mu-Base} (3))
$( \xbm CUT)$ $ \xch $ (by Proposition \ref{Proposition Alg-Log} (7.2) $
\xCf (CUT)$ $ \xch $ $ \xCf (PR)$
we would have a proof of $( \xbm PR)$ $ \xch $ $ \xCf (PR)$ without $(
\xbm dp),$ which is impossible,
as shown by Example \ref{Example Pref-Dp}.

(2) Reconsider Example \ref{Example Need-Pr}, and say $a \xcm p \xcu q,$
$b \xcm p \xcu \xCN q,$ $c \xcm \xCN p \xcu q.$
It is shown there that $( \xbm CUM)$ holds, so $( \xbm CUT)$ holds, so by
Proposition \ref{Proposition Alg-Log} (7.2) $ \xCf (CUT)$ holds, if we
define $ \ol{ \ol{T} }:=Th(f(M(T)).$
Set $T:=\{p \xco ( \xCN p \xcu q)\},$ $T':=\{p\},$ then
$ \ol{ \ol{T \xcv T' } }= \ol{ \ol{T' } }= \ol{\{p \xcu \xCN q\}},$ $ \ol{
\ol{T} }= \ol{T},$ $ \ol{T \xcv T' }= \ol{T' }= \ol{\{p\}},$ $so$ $ \xCf
(PR)$ $fails.$

$ \xcz $
\\[3ex]
\newpage

\section{
Preferential structures
}
\label{Preferential structures}


\subsection{
General and smooth preferential structures
}
\label{General and smooth preferential structures}


\subsubsection{
Definitions and basics
}
\label{General preferential structures - Definitions}

\index{Definition Pref-Str}

\bd

$\hspace{0.01em}$


\label{Definition Pref-Str}

Fix $U \xEd \xCQ,$ and consider arbitrary $X.$
Note that this $X$ has not necessarily anything to do with $U,$ or $ \xdu
$ below.
Thus, the functions $ \xbm_{ \xdm }$ below are in principle functions from
$V$ to $V$ - where $V$
is the set theoretical universe we work in.

(A) Preferential models or structures.

(1) The version without copies:

A pair $ \xdm:=<U, \xeb >$ with $U$ an arbitrary set, and $ \xeb $ an
arbitrary binary relation
is called a preferential model or structure.

(2) The version with copies:

A pair $ \xdm:=< \xdu, \xeb >$ with $ \xdu $ an arbitrary set of pairs,
and $ \xeb $ an arbitrary binary
relation is called a preferential model or structure.

If $<x,i> \xbe \xdu,$ then $x$ is intended to be an element of $U,$ and
$i$ the index of the
copy.

We sometimes also need copies of the relation $ \xeb,$ we will then
replace $ \xeb $
by one or several arrows $ \xba $ attacking non-minimal elements, e.g. $x
\xeb y$ will
be written $ \xba:x \xcp y,$ $<x,i> \xeb <y,i>$ will be written $ \xba
:<x,i> \xcp <y,i>,$ and
finally we might have $< \xba,k>:x \xcp y$ and $< \xba,k>:<x,i> \xcp
<y,i>,$ etc.

(B) Minimal elements, the functions $ \xbm_{ \xdm }$

(1) The version without copies:

Let $ \xdm:=<U, \xeb >,$ and define

$ \xbm_{ \xdm }(X)$ $:=$ $\{x \xbe X:$ $x \xbe U$ $ \xcu $ $ \xCN \xcE x'
\xbe X \xcs U.x' \xeb x\}.$

$ \xbm_{ \xdm }(X)$ is called the set of minimal elements of $X$ (in $
\xdm ).$

(2) The version with copies:

Let $ \xdm:=< \xdu, \xeb >$ be as above. Define

$ \xbm_{ \xdm }(X)$ $:=$ $\{x \xbe X:$ $ \xcE <x,i> \xbe \xdu. \xCN \xcE
<x',i' > \xbe \xdu (x' \xbe X$ $ \xcu $ $<x',i' >' \xeb <x,i>)\}.$

Again, by abuse of language, we say that $ \xbm_{ \xdm }(X)$ is the set of
minimal elements
of $X$ in the structure. If the context is clear, we will also write just
$ \xbm.$

We sometimes say that $<x,i>$ ``kills'' or ``minimizes'' $<y,j>$ if
$<x,i> \xeb <y,j>.$ By abuse of language we also say a set $X$ kills or
minimizes a set
$Y$ if for all $<y,j> \xbe \xdu,$ $y \xbe Y$ there is $<x,i> \xbe \xdu,$
$x \xbe X$ s.t. $<x,i> \xeb <y,j>.$

$ \xdm $ is also called injective or 1-copy, iff there is always at most
one copy
$<x,i>$ for each $x.$ Note that the existence of copies corresponds to a
non-injective labelling function - as is often used in nonclassical
logic, e.g. modal logic.

We say that $ \xdm $ is transitive, irreflexive, etc., iff $ \xeb $ is.

Note that $ \xbm (X)$ might well be empty, even if $X$ is not.
\index{Definition Pref-Log}

\ed

\bd

$\hspace{0.01em}$


\label{Definition Pref-Log}

We define the consequence relation of a preferential structure for a
given propositional language $ \xdl.$

(A)

(1) If $m$ is a classical model of a language $ \xdl,$ we say by abuse
of language

$<m,i> \xcm \xbf $ iff $m \xcm \xbf,$

and if $X$ is a set of such pairs, that

$X \xcm \xbf $ iff for all $<m,i> \xbe X$ $m \xcm \xbf.$

(2) If $ \xdm $ is a preferential structure, and $X$ is a set of $ \xdl
-$models for a
classical propositional language $ \xdl,$ or a set of pairs $<m,i>,$
where the $m$ are
such models, we call $ \xdm $ a classical preferential structure or model.

(B)

Validity in a preferential structure, or the semantical consequence
relation
defined by such a structure:

Let $ \xdm $ be as above.

We define:

$T \xcm_{ \xdm } \xbf $ iff $ \xbm_{ \xdm }(M(T)) \xcm \xbf,$ i.e. $
\xbm_{ \xdm }(M(T)) \xcc M( \xbf ).$

$ \xdm $ will be called definability preserving iff for all $X \xbe \xdD_{
\xdl }$ $ \xbm_{ \xdm }(X) \xbe \xdD_{ \xdl }.$

As $ \xbm_{ \xdm }$ is defined on $ \xdD_{ \xdl },$ but need by no means
always result in some new
definable set, this is (and reveals itself as a quite strong) additional
property.
\index{Example NeedCopies}

\ed

\be

$\hspace{0.01em}$


\label{Example NeedCopies}

This simple example illustrates the
importance of copies. Such examples seem to have appeared for the first
time
in print in  \cite{KLM90}, but can probably be
attibuted to folklore.

Consider the propositional language $ \xdl $ of two propositional
variables $p,q$, and
the classical preferential model $ \xdm $ defined by

$m \xcm p \xcu q,$ $m' \xcm p \xcu q,$ $m_{2} \xcm \xCN p \xcu q,$ $m_{3}
\xcm \xCN p \xcu \xCN q,$ with $m_{2} \xeb m$, $m_{3} \xeb m' $, and
let $ \xcm_{ \xdm }$ be its consequence relation. (m and $m' $ are
logically identical.)

Obviously, $Th(m) \xco \{ \xCN p\} \xcm_{ \xdm } \xCN p$, but there is no
complete theory $T' $ s.t.
$Th(m) \xco T' \xcm_{ \xdm } \xCN p$. (If there were one, $T' $ would
correspond to $m$, $m_{2},$ $m_{3},$
or the missing $m_{4} \xcm p \xcu \xCN q$, but we need two models to kill
all copies of $m.)$
On the other hand, if there were just one copy of $m,$ then one other
model,
i.e. a complete theory would suffice. More formally, if we admit at most
one
copy of each model in a structure $ \xdm,$ $m \xcM T,$ and $Th(m) \xco T
\xcm_{ \xdm } \xbf $ for some $ \xbf $ s.t.
$m \xcm \xCN \xbf $ - i.e. $m$ is not minimal in the models of $Th(m) \xco
T$ - then there is a
complete $T' $ with $T' \xcl T$ and $Th(m) \xco T' \xcm_{ \xdm } \xbf $,
i.e. there is $m'' $ with $m'' \xcm T' $ and
$m'' \xeb m.$ $ \xcz $
\\[3ex]
\index{Definition Smooth}

\ee

\bd

$\hspace{0.01em}$


\label{Definition Smooth}

Let $ \xdy \xcc \xdp (U).$ (In applications to logic, $ \xdy $ will be $
\xdD_{ \xdl }.)$

A preferential structure $ \xdm $ is called $ \xdy -$smooth iff in every
$X \xbe \xdy $ every element
$x \xbe X$ is either minimal in $X$ or above an element, which is minimal
in $X.$ More
precisely:

(1) The version without copies:

If $x \xbe X \xbe \xdy,$ then either $x \xbe \xbm (X)$ or there is $x'
\xbe \xbm (X).x' \xeb x.$

(2) The version with copies:

If $x \xbe X \xbe \xdy,$ and $<x,i> \xbe \xdu,$ then either there is no
$<x',i' > \xbe \xdu,$ $x' \xbe X,$
$<x',i' > \xeb <x,i>$ or there is $<x',i' > \xbe \xdu,$ $<x',i' > \xeb
<x,i>,$ $x' \xbe X,$ s.t. there is
no $<x'',i'' > \xbe \xdu,$ $x'' \xbe X,$ with $<x'',i'' > \xeb <x',i'
>.$

When considering the models of a language $ \xdl,$ $ \xdm $ will be
called smooth iff
it is $ \xdD_{ \xdl }-$smooth; $ \xdD_{ \xdl }$ is the default.

Obviously, the richer the set $ \xdy $ is, the stronger the condition $
\xdy -$smoothness
will be.
\newpage

\subsubsection{
Representation
}
\label{General preferential structures - Representation}

\index{Table Pref-Representation-With-Ref}

\ed

The following table summarizes representation by general or smooth
preferential structures. The implications on the right are shown in
Proposition \ref{Proposition Alg-Log} (going via the $ \xbm -$functions),
those on the left
are shown in the respective representation theorems.
\label{Table Pref-Representation-With-Ref}

{\scriptsize

\begin{tabular}{|c|c|c|c|c|}

\hline

$\xbm-$ function
\xEH
\xEH
Pref.Structure
\xEH
\xEH
Logic
\xEP

\hline

$(\xbm \xcc)+(\xbm PR)$
\xEH
$\xci$
\xEH
general
\xEH
$\xch$ $(\xbm dp)$
\xEH
$(LLE)+(RW)+(SC)+(PR)$
\xEP

\xEH
Fact \ref{Fact Pref-Sound}
\xEH
\xEH
\xEH
\xEP

\cline{2-2}
\cline{4-4}

\xEH
$\xch$
\xEH
\xEH
$\xci$
\xEH
\xEP

\xEH
Proposition \ref{Proposition Pref-Complete}
\xEH
\xEH
\xEH
\xEP

\cline{2-2}
\cline{4-4}

\xEH
\xEH
\xEH
$\xcH$ without $(\xbm dp)$
\xEH
\xEP

\xEH
\xEH
\xEH
Example \ref{Example Pref-Dp}
\xEH
\xEP

\hline

$(\xbm \xcc)+(\xbm PR)$
\xEH
$\xci$
\xEH
transitive
\xEH
$\xch$ $(\xbm dp)$
\xEH
$(LLE)+(RW)+(SC)+(PR)$
\xEP

\xEH
Fact \ref{Fact Pref-Sound}
\xEH
\xEH
\xEH
\xEP

\cline{2-2}
\cline{4-4}

\xEH
$\xch$
\xEH
\xEH
$\xci$
\xEH
\xEP

\xEH
Proposition \ref{Proposition Pref-Complete-Trans}
\xEH
\xEH
\xEH
\xEP

\cline{2-2}
\cline{4-4}

\xEH
\xEH
\xEH
$\xcH$ without $(\xbm dp)$
\xEH
\xEP

\xEH
\xEH
\xEH
Example \ref{Example Pref-Dp}
\xEH
\xEP

\hline

$(\xbm \xcc)+(\xbm PR)+(\xbm CUM)$
\xEH
$\xci$
\xEH
smooth
\xEH
$\xch$ $(\xbm dp)$
\xEH
$(LLE)+(RW)+(SC)+(PR)+$
\xEP

\xEH
Fact \ref{Fact Smooth-Sound}
\xEH
\xEH
\xEH
$(CUM)$
\xEP

\cline{2-2}
\cline{4-4}

\xEH
$\xch$ $(\xcv)$
\xEH
\xEH
$\xci$ $(\xcv)$
\xEH
\xEP

\xEH
Proposition \ref{Proposition Smooth-Complete}
\xEH
\xEH
\xEH
\xEP

\cline{2-2}
\cline{4-4}

\xEH
\xEH
\xEH
$\xcH$ without $(\xbm dp)$
\xEH
\xEP

\xEH
\xEH
\xEH
Example \ref{Example Pref-Dp}
\xEH
\xEP

\hline

$(\xbm \xcc)+(\xbm PR)+(\xbm CUM)$
\xEH
$\xci$
\xEH
smooth+transitive
\xEH
$\xch$ $(\xbm dp)$
\xEH
$(LLE)+(RW)+(SC)+(PR)+$
\xEP

\xEH
Fact \ref{Fact Smooth-Sound}
\xEH
\xEH
\xEH
$(CUM)$
\xEP

\cline{2-2}
\cline{4-4}

\xEH
$\xch$ $(\xcv)$
\xEH
\xEH
$\xci$ $(\xcv)$
\xEH
\xEP

\xEH
Proposition \ref{Proposition Smooth-Complete-Trans}
\xEH
\xEH
\xEH
\xEP

\cline{2-2}
\cline{4-4}

\xEH
\xEH
\xEH
$\xcH$ without $(\xbm dp)$
\xEH
\xEP

\xEH
\xEH
\xEH
Example \ref{Example Pref-Dp}
\xEH
\xEP

\hline

\end{tabular}

}

\index{Fact Pref-Sound}

\bfa

$\hspace{0.01em}$


\label{Fact Pref-Sound}

$( \xbm \xcc )$ and $( \xbm PR)$ hold in all preferential structures.
\index{Fact Pref-Sound Proof}

\efa

\subparagraph{
Proof
}

$\hspace{0.01em}$


Trivial. The central argument is: if $x,y \xbe X \xcc Y,$ and $x \xeb y$
in $X,$ then also
$x \xeb y$ in $Y.$

$ \xcz $
\\[3ex]
\index{Fact Smooth-Sound}

\bfa

$\hspace{0.01em}$


\label{Fact Smooth-Sound}

$( \xbm \xcc ),$ $( \xbm PR),$ and $( \xbm CUM)$ hold in all smooth
preferential structures.
\index{Fact Smooth-Sound Proof}

\efa

\subparagraph{
Proof
}

$\hspace{0.01em}$


By Fact \ref{Fact Pref-Sound}, we only have to show $( \xbm CUM).$
By Fact \ref{Fact Mu-Base}, $( \xbm CUT)$ follows from $( \xbm PR),$
so it remains to show
$( \xbm CM).$ So suppose $ \xbm (X) \xcc Y \xcc X,$ we have to show $ \xbm
(Y) \xcc \xbm (X).$ Let
$x \xbe X- \xbm (X),$ so there is $x' \xbe X,$ $x' \xeb x,$ by smoothness,
there must be $x'' \xbe \xbm (X),$
$x'' \xeb x,$ so $x'' \xbe Y,$ and $x \xce \xbm (Y).$ The proof for the
case with copies is
analogous.
\index{Example Pref-Dp}

\be

$\hspace{0.01em}$


\label{Example Pref-Dp}

This example was first given in [Sch92]. It shows
that condition $ \xCf (PR)$ may fail in preferential structures which are
not
definability preserving.

Let $v( \xdl ):=\{p_{i}:i \xbe \xbo \},$ $n,n' \xbe M_{ \xdl }$ be defined
by $n \xcm \{p_{i}:i \xbe \xbo \},$
$n' \xcm \{ \xCN p_{0}\} \xcv \{p_{i}:0<i< \xbo \}.$

Let $ \xdm:=<M_{ \xdl }, \xeb >$ where only $n \xeb n',$ i.e. just two
models are
comparable. Note that the structure is transitive and smooth.
Thus, by Fact \ref{Fact Smooth-Sound} $( \xbm \xcc ),$ $( \xbm PR),$
$( \xbm CUM)$ hold.

Let $ \xbm:= \xbm_{ \xdm },$ and $ \xcn $ be defined as usual by $ \xbm
.$

Set $T:= \xCQ,$ $T':=\{p_{i}:0<i< \xbo \}.$ We have $M_{T}=M_{ \xdl },$
$f(M_{T})=M_{ \xdl }-\{n' \},$ $M_{T' }=\{n,n' \},$
$f(M_{T' })=\{n\}.$ So by the result of Example \ref{Example Not-Def},
$f$ is not
definability preserving, and, furthermore,
$ \ol{ \ol{T} }= \ol{T},$ $ \ol{ \ol{T' } }= \ol{\{p_{i}:i< \xbo \}},$
$so$ $p_{0} \xbe \ol{ \ol{T \xcv T' } },$ $but$ $ \ol{ \ol{ \ol{T} } \xcv
T' }= \ol{ \ol{T} \xcv T' }= \ol{T' },$ $so$ $p_{0} \xce
 \ol{ \ol{ \ol{T} } \xcv T' },$
contradicting $ \xCf (PR),$ which holds in all definability preserving
preferential structures $ \xcz $
\\[3ex]
\index{Proposition Pref-Complete}

\ee

\bp

$\hspace{0.01em}$


\label{Proposition Pref-Complete}

Let $ \xbm: \xdy \xcp \xdp (U)$ satisfy $( \xbm \xcc )$ and $( \xbm PR).$
Then there is a preferential
structure $ \xdx $ s.t. $ \xbm = \xbm_{ \xdx }.$
\index{Proposition Pref-Complete-Trans}

\ep

\bp

$\hspace{0.01em}$


\label{Proposition Pref-Complete-Trans}

Let $ \xbm: \xdy \xcp \xdp (U)$ satisfy $( \xbm \xcc )$ and $( \xbm PR).$
Then there is a transitive
preferential structure $ \xdx $ s.t. $ \xbm = \xbm_{ \xdx }.$
\index{Proposition Smooth-Complete}

\ep

\bp

$\hspace{0.01em}$


\label{Proposition Smooth-Complete}

Let $ \xbm: \xdy \xcp \xdp (U)$ satisfy $( \xbm \xcc ),$ $( \xbm PR),$
and $( \xbm CUM),$ and the domain $ \xdy $ $( \xcv ).$

Then there is a $ \xdy -$smooth preferential structure $ \xdx $ s.t. $
\xbm = \xbm_{ \xdx }.$
\index{Proposition Smooth-Complete-Trans}

\ep

\bp

$\hspace{0.01em}$


\label{Proposition Smooth-Complete-Trans}

Let $ \xbm: \xdy \xcp \xdp (U)$ satisfy $( \xbm \xcc ),$ $( \xbm PR),$
and $( \xbm CUM),$ and the domain $ \xdy $ $( \xcv ).$

Then there is a transitive $ \xdy -$smooth preferential structure $ \xdx $
s.t. $ \xbm = \xbm_{ \xdx }.$
\newpage

\subsection{
Ranked structures
}


\subsubsection{
Definitions and basics
}
\label{Ranked - Definitions}

\index{Fact Rank-Base}

\ep

\bfa

$\hspace{0.01em}$


\label{Fact Rank-Base}

Let $ \xeb $ be an irreflexive, binary relation on $X,$ then the following
two conditions
are equivalent:

(1) There is $ \xbO $ and an irreflexive, total, binary relation $ \xeb '
$ on $ \xbO $ and a
function $f:X \xcp \xbO $ s.t. $x \xeb y$ $ \xcr $ $f(x) \xeb ' f(y)$ for
all $x,y \xbe X.$

(2) Let $x,y,z \xbe X$ and $x \xcT y$ wrt. $ \xeb $ (i.e. neither $x \xeb
y$ nor $y \xeb x),$ then $z \xeb x$ $ \xcp $ $z \xeb y$
and $x \xeb z$ $ \xcp $ $y \xeb z.$

$ \xcz $
\\[3ex]
\index{Definition Rank-Rel}

\efa

\bd

$\hspace{0.01em}$


\label{Definition Rank-Rel}

We call an irreflexive, binary relation $ \xeb $ on $X,$ which satisfies
(1)
(equivalently (2)) of Fact \ref{Fact Rank-Base}, ranked.
By abuse of language, we also call a preferential structure $<X, \xeb >$
ranked, iff
$ \xeb $ is.
\index{Fact Rank-Auxil}

\ed

\bfa

$\hspace{0.01em}$


\label{Fact Rank-Auxil}

$M(T)-M(T' )$ is normally not definable.

In the presence of $( \xbm =)$ and $( \xbm \xcc ),$ $f(Y) \xcs (X-f(X))
\xEd \xCQ $ is equivalent to
$f(Y) \xcs X \xEd \xCQ $ and $f(Y) \xcs f(X)= \xCQ.$
\index{Fact Rank-Auxil Proof}

\efa

\subparagraph{
Proof
}

$\hspace{0.01em}$


$f(Y) \xcs (X-f(X))$ $=$ $(f(Y) \xcs X)-(f(Y) \xcs f(X)).$

`` $ \xci $ '': Let $f(Y) \xcs X \xEd \xCQ,$ $f(Y) \xcs f(X)= \xCQ,$ so
$f(Y) \xcs (X-f(X)) \xEd \xCQ.$

`` $ \xch $ '': Suppose $f(Y) \xcs (X-f(X)) \xEd \xCQ,$ so $f(Y) \xcs X
\xEd \xCQ.$ Suppose $f(Y) \xcs f(X) \xEd \xCQ,$ so
by $( \xbm \xcc )$ $f(Y) \xcs X \xcs Y \xEd \xCQ,$ so
by $( \xbm =)$ $f(Y) \xcs X \xcs Y=f(X \xcs Y),$ and $f(X) \xcs X \xcs Y
\xEd \xCQ,$ so by $( \xbm =)$
$f(X) \xcs X \xcs Y=f(X \xcs Y),$ so $f(X) \xcs Y=f(Y) \xcs X$ and $f(Y)
\xcs (X-f(X))= \xCQ.$

$ \xcz $
\\[3ex]
\index{Fact Rank-Trans}

\bfa

$\hspace{0.01em}$


\label{Fact Rank-Trans}

If $ \xeb $ on $X$ is ranked, and free of cycles, then $ \xeb $ is
transitive.
\index{Fact Rank-Trans Proof}

\efa

\subparagraph{
Proof
}

$\hspace{0.01em}$


Let $x \xeb y \xeb z.$ If $x \xcT z,$ then $y \xee z,$ resulting in a
cycle of length 2. If $z \xeb x,$ then
we have a cycle of length 3. So $x \xeb z.$ $ \xcz $
\\[3ex]
\index{Remark RatM=}

\br

$\hspace{0.01em}$


\label{Remark RatM=}

Note that $( \xbm =' )$ is very close to $ \xCf (RatM):$ $ \xCf (RatM)$
says:
$ \xba \xcn \xbb,$ $ \xba \xcN \xCN \xbg $ $ \xch $ $ \xba \xcu \xbg \xcn
\xbb.$ Or, $f(A) \xcc B,$ $f(A) \xcs C \xEd \xCQ $ $ \xcp $
$f(A \xcs C) \xcc B$ for all $A,B,C.$ This is not quite, but almost: $f(A
\xcs C) \xcc f(A) \xcs C$
(it depends how many $B$ there are, if $f(A)$ is some such $B,$ the fit is
perfect).
\index{Fact Rank-Hold}

\er

\bfa

$\hspace{0.01em}$


\label{Fact Rank-Hold}

In all ranked structures, $( \xbm \xcc ),$ $( \xbm =),$ $( \xbm PR),$ $(
\xbm =' ),$ $( \xbm \xFO ),$ $( \xbm \xcv ),$ $( \xbm \xcv ' ),$
$( \xbm \xbe ),$ $( \xbm RatM)$ will hold, if the corresponding closure
conditions are
satisfied.
\index{Fact Rank-Hold Proof}

\efa

\subparagraph{
Proof:
}

$\hspace{0.01em}$


\label{Section Proof:}

$( \xbm \xcc )$ and $( \xbm PR)$ hold in all preferential structures.

$( \xbm =)$ and $( \xbm =' )$ are trivial.

$( \xbm \xcv )$ and $( \xbm \xcv ' ):$ All minimal copies of elements in
$f(Y)$ have the same rank.
If some $y \xbe f(Y)$ has all its minimal copies killed by an element $x
\xbe X,$ by
rankedness, $x$ kills the rest, too.

$( \xbm \xbe ):$ If $f(\{a\})= \xCQ,$ we are done. Take the minimal
copies of a in $\{a\},$ they are
all killed by one element in $X.$

$( \xbm \xFO ):$ Case $f(X)= \xCQ:$ If below every copy of $y \xbe Y$
there is a copy of some $x \xbe X,$
then $f(X \xcv Y)= \xCQ.$ Otherwise $f(X \xcv Y)=f(Y).$ Suppose now $f(X)
\xEd \xCQ,$ $f(Y) \xEd \xCQ,$ then
the minimal ranks decide: if they are equal, $f(X \xcv Y)=f(X) \xcv f(Y),$
etc.

$( \xbm RatM):$ Let $X \xcc Y,$ $y \xbe X \xcs f(Y) \xEd \xCQ,$ $x \xbe
f(X).$ By rankedness, $y \xeb x,$ or
$y \xcT x,$ $y \xeb x$ is impossible, as $y \xbe X,$ so $y \xcT x,$ and $x
\xbe f(Y).$

$ \xcz $
\\[3ex]
\newpage

\subsubsection{
Representation
}
\label{Ranked - Representation}

\index{Definition 1-infin}

\bd

$\hspace{0.01em}$


\label{Definition 1-infin}

Let $ \xdz =< \xdx, \xeb >$ be a preferential structure. Call $ \xdz $
$1- \xca $ over $Z,$
iff for all $x \xbe Z$ there are exactly one or infinitely many copies of
$x,$ i.e.
for all $x \xbe Z$ $\{u \xbe \xdx:$ $u=<x,i>$ for some $i\}$ has
cardinality 1 or $ \xcg \xbo.$
\index{Lemma 1-infin}

\ed

\bl

$\hspace{0.01em}$


\label{Lemma 1-infin}

Let $ \xdz =< \xdx, \xeb >$ be a preferential structure and
$f: \xdy \xcp \xdp (Z)$ with $ \xdy \xcc \xdp (Z)$ be represented by $
\xdz,$ i.e. for $X \xbe \xdy $ $f(X)= \xbm_{ \xdz }(X),$
and $ \xdz $ be ranked and free of cycles. Then there is a structure $
\xdz ' $, $1- \xca $ over
$Z,$ ranked and free of cycles, which also represents $f.$
\index{Lemma 1-infin Proof}

\el

\subparagraph{
Proof
}

$\hspace{0.01em}$


We construct $ \xdz ' =< \xdx ', \xeb ' >.$

Let $A:=\{x \xbe Z$: there is some $<x,i> \xbe \xdx,$ but for all $<x,i>
\xbe \xdx $ there is
$<x,j> \xbe \xdx $ with $<x,j> \xeb <x,i>\},$

let $B:=\{x \xbe Z$: there is some $<x,i> \xbe \xdx,$ s.t. for no $<x,j>
\xbe \xdx $ $<x,j> \xeb <x,i>\},$

let $C:=\{x \xbe Z$: there is no $<x,i> \xbe \xdx \}.$

Let $c_{i}:i< \xbk $ be an enumeration of $C.$ We introduce for each such
$c_{i}$ $ \xbo $ many
copies $<c_{i},n>:n< \xbo $ into $ \xdx ',$ put all $<c_{i},n>$ above all
elements in $ \xdx,$ and order
the $<c_{i},n>$ by $<c_{i},n> \xeb ' <c_{i' },n' >$ $: \xcr $ $(i=i' $ and
$n>n' )$ or $i>i'.$ Thus, all $<c_{i},n>$ are
comparable.

If $a \xbe A,$ then there are infinitely many copies of a in $ \xdx,$ as
$ \xdx $ was
cycle-free, we put them all into $ \xdx '.$
If $b \xbe B,$ we choose exactly one such minimal element $<b,m>$ (i.e.
there
is no $<b,n> \xeb <b,m>)$ into $ \xdx ',$ and omit all other
elements. (For definiteness, assume in all applications $m=0.)$
For all elements from A and $B,$ we take the restriction of the order $
\xeb $ of $ \xdx.$
This is the new structure $ \xdz '.$

Obviously, adding the $<c_{i},n>$ does not introduce cycles, irreflexivity
and
rankedness are preserved. Moreover, any substructure of a cycle-free,
irreflexive,
ranked structure also has these properties, so $ \xdz ' $ is $1- \xca $
over $Z,$ ranked and
free of cycles.

We show that $ \xdz $ and $ \xdz ' $ are equivalent. Let then $X \xcc Z,$
we have to prove
$ \xbm (X)= \xbm ' (X)$ $( \xbm:= \xbm_{ \xdz }$, $ \xbm ':= \xbm_{
\xdz ' }).$

Let $z \xbe X- \xbm (X).$ If $z \xbe C$ or $z \xbe A,$ then $z \xce \xbm '
(X).$ If $z \xbe B,$
let $<z,m>$ be the chosen element. As $z \xce \xbm (X),$ there is $x \xbe
X$ s.t. some $<x,j> \xeb <z,m>.$
$x$ cannot be in $C.$ If $x \xbe A,$ then also $<x,j> \xeb ' <z,m>$. If
$x \xbe B,$ then there is some
$<x,k>$ also in $ \xdx '.$ $<x,j> \xeb <x,k>$ is impossible. If $<x,k>
\xeb <x,j>,$ then $<z,m> \xee <x,k>$
by transitivity. If $<x,k> \xcT <x,j>$, then also $<z,m> \xee <x,k>$ by
rankedness. In any
case, $<z,m> \xee ' <x,k>,$ and thus $z \xce \xbm ' (X).$

Let $z \xbe X- \xbm ' (X).$ If $z \xbe C$ or $z \xbe A,$ then $z \xce \xbm
(X).$ Let $z \xbe B,$ and some $<x,j> \xeb ' <z,m>.$
$x$ cannot be in $C,$ as they were sorted on top, so $<x,j>$ exists in $
\xdx $ too and
$<x,j> \xeb <z,m>.$ But if any other $<z,i>$ is also minimal in $ \xdz $
among the $<z,k>,$
then by rankedness also $<x,j> \xeb <z,i>,$ as $<z,i> \xcT <z,m>,$ so $z
\xce \xbm (X).$ $ \xcz $
\\[3ex]
\index{Proposition Rank-Rep1}

\bp

$\hspace{0.01em}$


\label{Proposition Rank-Rep1}

The first result applies for structures without copies of elements.

(1) Let $ \xdy \xcc \xdp (U)$ be closed under finite unions.
Then $( \xbm \xcc ),$ $( \xbm \xCQ ),$ $( \xbm =)$ characterize ranked
structures for which for all
$X \xbe \xdy $ $X \xEd \xCQ $ $ \xcp $ $ \xbm_{<}(X) \xEd \xCQ $ hold,
i.e. $( \xbm \xcc ),$ $( \xbm \xCQ ),$ $( \xbm =)$ hold in such
structures for $ \xbm_{<},$ and if they hold for some $ \xbm,$ we can
find a ranked relation
$<$ on $U$ s.t. $ \xbm = \xbm_{<}.$ Moreover, the structure can be choosen
$ \xdy -$smooth.

(2) Let $ \xdy \xcc \xdp (U)$ be closed under finite unions, and contain
singletons.
Then $( \xbm \xcc ),$ $( \xbm \xCQ fin),$ $( \xbm =),$ $( \xbm \xbe )$
characterize ranked structures for which
for all finite $X \xbe \xdy $ $X \xEd \xCQ $ $ \xcp $ $ \xbm_{<}(X) \xEd
\xCQ $ hold, i.e. $( \xbm \xcc ),$ $( \xbm \xCQ fin),$ $( \xbm =),$ $(
\xbm \xbe )$
hold in such structures for $ \xbm_{<},$ and if they hold for some $ \xbm
,$ we can find
a ranked relation $<$ on $U$ s.t. $ \xbm = \xbm_{<}.$

Note that the prerequisites of (2) hold in particular in the case
of ranked structures without copies, where all elements of $U$ are present
in the
structure - we need infinite descending chains to have $ \xbm (X)= \xCQ $
for $X \xEd \xCQ.$
\index{Fact Rank-No-Rep}

\ep

\bfa

$\hspace{0.01em}$


\label{Fact Rank-No-Rep}

$( \xbm \xcc )+( \xbm PR)+( \xbm =)+( \xbm \xcv )+( \xbm \xbe )$ do not
imply representation by a ranked
structure.
\index{Fact Rank-No-Rep Proof}

\efa

\subparagraph{
Proof
}

$\hspace{0.01em}$


See Example \ref{Example Rank-Copies}. $ \xcz $
\\[3ex]
\index{Example Rank-Copies}

\be

$\hspace{0.01em}$


\label{Example Rank-Copies}

This example shows that the conditions $( \xbm \xcc )+( \xbm PR)+( \xbm
=)+( \xbm \xcv )+( \xbm \xbe )$
can be satisfied, and still representation by a ranked structure
is impossible.

Consider $ \xbm (\{a,b\})= \xCQ,$ $ \xbm (\{a\})=\{a\},$ $ \xbm
(\{b\})=\{b\}.$ The conditions
$( \xbm \xcc )+( \xbm PR)+( \xbm =)+( \xbm \xcv )+( \xbm \xbe )$
hold trivially. This is representable, e.g. by $a_{1} \xed b_{1} \xed
a_{2} \xed b_{2} \Xl $ without
transitivity. (Note that rankedness implies transitivity,
$a \xec b \xec c,$ but not for $a=c.)$ But this cannot be represented by a
ranked
structure: As $ \xbm (\{a\}) \xEd \xCQ,$ there must be a copy $a_{i}$ of
minimal rank, likewise for
$b$ and some $b_{i}.$ If they have the same rank, $ \xbm
(\{a,b\})=\{a,b\},$ otherwise it will be
$\{a\}$ or $\{b\}.$

$ \xcz $
\\[3ex]
\index{Proposition Rank-Rep2}

\ee

\bp

$\hspace{0.01em}$


\label{Proposition Rank-Rep2}

Let $ \xdy $ be closed under finite unions and contain singletons. Then
$( \xbm \xcc )+( \xbm PR)+( \xbm \xFO )+( \xbm \xcv )+( \xbm \xbe )$
characterize ranked structures,
where elements may appear in several copies.
\newpage

\section{
Theory revision
}


\subsection{
AGM revision
}
\label{Section AGM-revision}

\ep

All material in this Section \ref{Section AGM-revision} is due verbatim
or in essence
to AGM - AGM for Alchourron, Gardenfors, Makinson,
see e.g.  \cite{AGM85}.
\index{Definition AGM}

\bd

$\hspace{0.01em}$


\label{Definition AGM}

We present in parallel the logical
and the semantic (or purely algebraic) side. For the latter, we work in
some
fixed universe $U,$ and the intuition is $U=M_{ \xdl },$ $X=M(K),$ etc.,
so, e.g. $A \xbe K$
becomes $X \xcc B,$ etc.

(For reasons of readability, we omit most caveats about definability.)

$K_{ \xcT }$ will denote the inconsistent theory.

We consider two functions, - and $*,$ taking a deductively closed theory
and a
formula as arguments, and returning a (deductively closed) theory on the
logics
side. The algebraic counterparts work on definable model sets. It is
obvious
that $ \xCf (K-1),$ $(K*1),$ $ \xCf (K-6),$ $(K*6)$ have vacuously true
counterparts on the
semantical side. Note that $K$ $ \xCf (X)$ will never change, everything
is relative
to fixed $K$ $ \xCf (X).$ $K* \xbf $ is the result of revising $K$ with $
\xbf.$ $K- \xbf $ is the result of
subtracting enough from $K$ to be able to add $ \xCN \xbf $ in a
reasonable way, called
contraction.

Moreover,
let $ \xck_{K}$ be a relation on the formulas relative to a deductively
closed theory $K$
on the formulas of $ \xdl,$ and $ \xck_{X}$ a relation on $ \xdp (U)$ or
a suitable subset of $ \xdp (U)$
relative to fixed $X.$ When the context is clear, we simply write $ \xck
.$
$ \xck_{K}$ $( \xck_{X})$ is called a relation of epistemic entrenchment
for $K$ $ \xCf (X).$

The following table presents the ``rationality postulates'' for contraction
(-),
revision $(*)$ and epistemic entrenchment. In AGM tradition, $K$ will be a
deductively closed theory, $ \xbf, \xbq $ formulas. Accordingly, $X$ will
be the set of
models of a theory, $A,B$ the model sets of formulas.

\renewcommand{\arraystretch}{1.2}

{\small
\begin{tabular}{|c|c|c|c|}

\hline

\multicolumn{4}{|c|} {Contraction, $K-\xbf $} \xEP

\hline

$(K-1)$ \xEH $K-\xbf $ is deductively closed \xEH \xEH \xEP

\hline

$(K-2)$ \xEH $K-\xbf $ $ \xcc $ $K$ \xEH $(X \xDN 2)$ \xEH $X \xcc X \xDN A$
\xEP

\hline

$(K-3)$ \xEH $\xbf  \xce K$ $ \xch $ $K-\xbf =K$ \xEH $(X \xDN 3)$ \xEH $X \xcC
A$ $
\xch $ $X \xDN A=X$ \xEP

\hline

$(K-4)$ \xEH $ \xcL \xbf $ $ \xch $ $\xbf  \xce K-\xbf $ \xEH $(X \xDN 4)$ \xEH
$A \xEd
U$ $ \xch $ $X \xDN A \xcC A$ \xEP

\hline

$(K-5)$ \xEH $K \xcc \ol{(K-\xbf ) \xcv \{\xbf \}}$ \xEH $(X \xDN 5)$ \xEH $(X
\xDN
A) \xcs A$ $ \xcc $ $X$ \xEP

\hline

$(K-6)$ \xEH $ \xcl \xbf  \xcr \xbq $ $ \xch $ $K-\xbf =K-\xbq $ \xEH \xEH \xEP

\hline

$(K-7)$ \xEH $(K-\xbf ) \xcs (K-\xbq )  \xcc  $ \xEH
$(X \xDN 7)$ \xEH $X \xDN (A \xcs B)  \xcc  $ \xEP

\xEH $K-(\xbf  \xcu \xbq ) $ \xEH
\xEH $(X \xDN A) \xcv (X \xDN B)$ \xEP

\hline

$(K-8)$ \xEH $\xbf  \xce K-(\xbf  \xcu \xbq )  \xch  $ \xEH
$(X \xDN 8)$ \xEH $X \xDN (A \xcs B) \xcC A  \xch  $ \xEP

\xEH $K-(\xbf  \xcu \xbq ) \xcc K-\xbf $ \xEH
\xEH $X \xDN A \xcc X \xDN (A \xcs B)$ \xEP

\hline
\hline

\multicolumn{4}{|c|} {Revision, $K*\xbf $} \xEP

\hline

$(K*1)$ \xEH $K*\xbf $ is deductively closed \xEH - \xEH \xEP

\hline

$(K*2)$ \xEH $\xbf  \xbe K*\xbf $ \xEH $(X \xfA 2)$ \xEH $X \xfA A \xcc A$
\xEP

\hline

$(K*3)$ \xEH $K*\xbf $ $ \xcc $ $ \ol{K \xcv \{\xbf \}}$ \xEH $(X \xfA 3)$ \xEH
$X
\xcs A \xcc X \xfA A$ \xEP

\hline

$(K*4)$ \xEH $ \xCN \xbf  \xce K  \xch $ \xEH
$(X \xfA 4)$ \xEH $X \xcs A \xEd \xCQ   \xch $ \xEP

\xEH $\ol{K \xcv \{\xbf \}}  \xcc  K*\xbf $ \xEH
\xEH $X \xfA A \xcc X \xcs A$ \xEP

\hline

$(K*5)$ \xEH $K*\xbf =K_{ \xcT }$ $ \xch $ $ \xcl \xCN \xbf $ \xEH $(X \xfA 5)$
\xEH $X \xfA A= \xCQ $ $ \xch $ $A= \xCQ $ \xEP

\hline

$(K*6)$ \xEH $ \xcl \xbf  \xcr \xbq $ $ \xch $ $K*\xbf =K*\xbq $ \xEH - \xEH
\xEP

\hline

$(K*7)$ \xEH $K*(\xbf  \xcu \xbq )  \xcc $ \xEH
$(X \xfA 7)$ \xEH $(X \xfA A) \xcs B  \xcc  $ \xEP

\xEH $\ol{(K*\xbf ) \xcv \{\xbq \}}$ \xEH
\xEH $X \xfA (A \xcs B)$ \xEP

\hline

$(K*8)$ \xEH $ \xCN \xbq  \xce K*\xbf  \xch $ \xEH
$(X \xfA 8)$ \xEH $(X \xfA A) \xcs B \xEd \xCQ \xch $ \xEP

\xEH $\ol{(K*\xbf ) \xcv \{\xbq \}} \xcc K*(\xbf  \xcu \xbq )$ \xEH
\xEH $ X \xfA (A \xcs B) \xcc (X \xfA A) \xcs B$ \xEP

\hline
\hline

\multicolumn{4}{|c|} {Epistemic entrenchment} \xEP

\hline

$(EE1)$ \xEH $ \xck_{K}$ is transitive \xEH
$(EE1)$ \xEH $ \xck_{X}$ is transitive \xEP

\hline

$(EE2)$ \xEH $\xbf  \xcl \xbq   \xch  \xbf  \xck_{K}\xbq $ \xEH
$(EE2)$ \xEH $A \xcc B  \xch  A \xck_{X}B$ \xEP

\hline

$(EE3)$ \xEH $ \xcA  \xbf,\xbq  $ \xEH
$(EE3)$ \xEH $ \xcA A,B $ \xEP
\xEH $ (\xbf  \xck_{K}\xbf  \xcu \xbq $ or $\xbq  \xck_{K}\xbf  \xcu \xbq )$
\xEH
\xEH $ (A \xck_{X}A \xcs B$ or $B \xck_{X}A \xcs B)$ \xEP

\hline

$(EE4)$ \xEH $K \xEd K_{ \xcT }  \xch  $ \xEH
$(EE4)$ \xEH $X \xEd \xCQ   \xch  $ \xEP
\xEH $(\xbf  \xce K$ iff $ \xcA  \xbq.\xbf  \xck_{K}\xbq )$ \xEH
\xEH $(X \xcC A$ iff $ \xcA B.A \xck_{X}B)$ \xEP

\hline

$(EE5)$ \xEH $ \xcA \xbq.\xbq  \xck_{K}\xbf   \xch   \xcl \xbf $ \xEH
$(EE5)$ \xEH $ \xcA B.B \xck_{X}A  \xch  A=U$ \xEP

\hline

\end{tabular}
}
\\

\index{Remark TR-Rank}

\ed

\br

$\hspace{0.01em}$


\label{Remark TR-Rank}

(1) Note that $(X \xfA 7)$ and $(X \xfA 8)$ express a central condition
for ranked
structures, see Section 3.10: If we note $X \xfA.$ by $f_{X}(.),$ we then
have:
$f_{X}(A) \xcs B \xEd \xCQ $ $ \xch $ $f_{X}(A \xcs B)=f_{X}(A) \xcs B.$

(2) It is trivial to see that AGM revision cannot be defined by an
individual
distance (see Definition 2.3.5 below):
Suppose $X \xfA Y$ $:=$ $\{y \xbe Y:$ $ \xcE x_{y} \xbe X( \xcA y' \xbe
Y.d(x_{y},y) \xck d(x_{y},y' ))\}.$
Consider $a,b,c.$ $\{a,b\} \xfA \{b,c\}=\{b\}$ by $(X \xfA 3)$ and $(X
\xfA 4),$ so $d(a,b)<d(a,c).$
But on the other hand $\{a,c\} \xfA \{b,c\}=\{c\},$ so $d(a,b)>d(a,c),$
$contradiction.$
\index{Proposition AGM-Equiv}

\er

\bp

$\hspace{0.01em}$


\label{Proposition AGM-Equiv}

Contraction, revision, and epistemic entrenchment are interdefinable by
the
following equations, i.e., if the defining side has the respective
properties,
so will the defined side.

\renewcommand{\arraystretch}{1.5}

{\scriptsize
\begin{tabular}{|c|c|}

\hline

$K*\xbf:= \ol{(K- \xCN \xbf )} \xcv {\xbf }$  \xEH  $X \xfA A:= (X \xDN  \xdC A)
\xcs A$ \xEP

\hline

$K-\xbf:= K \xcs (K* \xCN \xbf )$  \xEH  $X \xDN A:= X \xcv (X \xfA  \xdC A)$
\xEP

\hline

$K-\xbf:=\{\xbq  \xbe K:$ $(\xbf <_{K}\xbf  \xco \xbq $ or $ \xcl \xbf )\}$ \xEH

$
X \xDN A:=
\left\{
\begin{array}{rcl}
X & iff & A=U, \\
 \xcS \{B: X \xcc B \xcc U, A<_{X}A \xcv B\} & & otherwise \\
\end{array}
\right. $
\xEP

\hline

$
\xbf  \xck_{K}\xbq: \xcr
\left\{
\begin{array}{l}
\xcl \xbf  \xcu \xbq  \\
or \\
\xbf  \xce K-(\xbf  \xcu \xbq ) \\
\end{array}
\right. $
\xEH

$
A \xck_{X}B: \xcr
\left\{
\begin{array}{l}
A,B=U  \\
or \\
X \xDN (A \xcs B) \xcC A \\
\end{array}
\right. $
\xEP

\hline

\end{tabular}
}

\index{Intuit-Entrench}
\paragraph{A remark on intuition}

\ep

The idea of epistemic entrenchment is that $ \xbf $ is more entrenched
than $ \xbq $
(relative to $K)$ iff $M( \xCN \xbq )$ is closer to $M(K)$ than $M( \xCN
\xbf )$ is to $M(K).$ In
shorthand, the more we can twiggle $K$ without reaching $ \xCN \xbf,$ the
more $ \xbf $ is
entrenched. Truth is maximally entrenched - no twiggling whatever will
reach
falsity. The more $ \xbf $ is entrenched,
the more we are certain about it. Seen this way, the properties of
epistemic
entrenchment relations are very natural (and trivial): As only the closest
points of $M( \xCN \xbf )$ count (seen from $M(K)),$ $ \xbf $ or $ \xbq $
will be as entrenched as
$ \xbf \xcu \xbq,$ and there is a logically strongest $ \xbf ' $ which is
as entrenched as $ \xbf $ -
this is just the sphere around $M(K)$ with radius $d(M(K),M( \xCN \xbf
)).$
\newpage

\subsection{
Distance based revision
}


\subsubsection{
Definitions and basics
}
\label{TR - Definitions and basics}

\index{Definition Distance}

\bd

$\hspace{0.01em}$


\label{Definition Distance}

$d:U \xCK U \xcp Z$ is called a pseudo-distance on $U$ iff (d1) holds:

(d1) $Z$ is totally ordered by a relation $<.$

If, in addition, $Z$ has a $<-$smallest element 0, and (d2) holds, we say
that $d$
respects identity:

(d2) $d(a,b)=0$ iff $a=b.$

If, in addition, (d3) holds, then $d$ is called symmetric:

(d3) $d(a,b)=d(b,a).$

(For any $a,b \xbe U.)$

Note that we can force the triangle inequality to hold trivially (if we
can
choose the values in the real numbers): It suffices to choose the values
in
the set $\{0\} \xcv [0.5,1],$ i.e. in the interval from 0.5 to 1, or as 0.
\index{Definition Dist-Indiv-Coll}

\ed

\bd

$\hspace{0.01em}$


\label{Definition Dist-Indiv-Coll}

We define the collective and the individual variant of choosing the
closest
elements in the second operand by two operators,
$ \xfA, \xfB: \xdp (U) \xCK \xdp (U) \xcp \xdp (U):$

Let $d$ be a distance or pseudo-distance.

$X \xfA Y$ $:=$ $\{y \xbe Y:$ $ \xcE x_{y} \xbe X. \xcA x' \xbe X, \xcA y'
\xbe Y(d(x_{y},y) \xck d(x',y' )\}$

(the collective variant, used in theory revision)

and

$X \xfB Y$ $:=$ $\{y \xbe Y:$ $ \xcE x_{y} \xbe X. \xcA y' \xbe
Y(d(x_{y},y) \xck d(x_{y},y' )\}$

(the individual variant, used for counterfactual conditionals and theory
update).

Thus, $A \xfA_{d}B$ is the subset of $B$ consisting of all $b \xbe B$ that
are closest to A.
Note that, if $ \xCf A$ or $B$ is infinite, $A \xfA_{d}B$ may be empty,
even if $ \xCf A$ and $B$ are not
empty. A condition assuring nonemptiness will be imposed when necessary.
\index{Definition Dist-Repr}

\ed

\bd

$\hspace{0.01em}$


\label{Definition Dist-Repr}

An operation $ \xfA: \xdp (U) \xCK \xdp (U) \xcp \xdp (U)$ is
representable iff there is a
pseudo-distance $d:U \xCK U \xcp Z$ such that

$A \xfA B$ $=$ $A \xfA_{d}B$ $:=$ $\{b \xbe B:$ $ \xcE a_{b} \xbe A \xcA
a' \xbe A \xcA b' \xbe B(d(a_{b},b) \xck d(a',b' ))\}.$
\index{Definition TR*d}

\ed

The following is the central definition, it describes the way a revision
$*_{d}$ is
attached to a pseudo-distance $d$ on the set of models.

\bd

$\hspace{0.01em}$


\label{Definition TR*d}

$T*_{d}T' $ $:=$ $Th(M(T) \xfA_{d}M(T' )).$

$*$ is called representable iff there is a pseudo-distance $d$ on the set
of models
s.t. $T*T' =Th(M(T) \xfA_{d}M(T' )).$

\subsubsection{
Representation
}
\label{TR - Representation}

\index{Fact AGM-In-Dist}

\ed

\bfa

$\hspace{0.01em}$


\label{Fact AGM-In-Dist}

A distance based revision satisfies the AGM postulates provided:

(1) it respects identity, i.e. $d(a,a)<d(a,b)$ for all $a \xEd b,$

(2) it satisfies a limit condition: minima exist,

(3) it is definability preserving.

(It is trivial to see that the first two are necessary,
and Example \ref{Example TR-Dp} (2)
below shows the necessity of (3). In particular, (2) and (3) will hold for
finite languages.)
\index{Fact AGM-In-Dist Proof}

\efa

\subparagraph{
Proof:
}

$\hspace{0.01em}$


\label{Section Proof:}

We use $ \xfA $ to abbreviate $ \xfA_{d}.$ As a matter of fact, we show
slightly more, as
we admit also full theories on the right of $*.$

$(K*1),$ $(K*2),$ $(K*6)$ hold by definition, $(K*3)$ and $(K*4)$ as $d$
respects
identity, $(K*5)$ by existence of minima.

It remains to show $(K*7)$ and $(K*8),$ we do them together, and show:
If $T*T' $ is consistent with $T'',$ then $T*(T' \xcv T'' )$ $=$ $
\ol{(T*T' ) \xcv T'' }.$

Note that $M(S \xcv S' )=M(S) \xcs M(S' ),$ and that $M(S*S' )=M(S) \xfA
M(S' ).$ (The latter
is only true if $ \xfA $ is definability preserving.)
By prerequisite, $M(T*T' ) \xcs M(T'' ) \xEd \xCQ,$ so $(M(T) \xfA M(T'
)) \xcs M(T'' ) \xEd \xCQ.$
Let $A:=M(T),$ $B:=M(T' ),$ $C:=M(T'' ).$ `` $ \xcc $ '': Let $b \xbe A
\xfA (B \xcs C).$
By prerequisite, there is $b' \xbe (A \xfA B) \xcs C.$ Thus $d(A,b' ) \xcg
d(A,B \xcs C)=d(A,b).$
As $b \xbe B,$ $b \xbe A \xfA B,$ but $b \xbe C,$ too. `` $ \xcd $ '': Let
$b' \xbe (A \xfA B) \xcs C.$ Thus $d(A,b' )=$
$d(A,B) \xck d(A,B \xcs C),$ so by $b' \xbe B \xcs C$ $b' \xbe A \xfA (B
\xcs C).$
We conclude $M(T) \xfA (M(T' ) \xcs M(T'' ))$ $=$ $(M(T) \xfA M(T' )) \xcs
M(T'' ),$ thus that
$T*(T' \xcv T'' )= \ol{(T*T' ) \xcv T'' }.$

$ \xcz $
\\[3ex]
\index{Definition TR-Umgeb}

\bd

$\hspace{0.01em}$


\label{Definition TR-Umgeb}

For $X,Y \xEd \xCQ,$ set $U_{Y}(X):=\{z:d(X,z) \xck d(X,Y)\}.$
\index{Fact TR-Umgeb}

\ed

\bfa

$\hspace{0.01em}$


\label{Fact Tr-Umgeb}

Let $X,Y,Z \xEd \xCQ.$ Then

(1) $U_{Y}(X) \xcs Z \xEd \xCQ $ iff $(X \xfA (Y \xcv Z)) \xcs Z \xEd \xCQ
,$

(2) $U_{Y}(X) \xcs Z \xEd \xCQ $ iff $ \xdC Z \xck_{X} \xdC Y$ - where $
\xck_{X}$ is epistemic entrenchement relative
to $X.$
\index{Fact TR-Umgeb Proof}

\efa

\subparagraph{
Proof
}

$\hspace{0.01em}$


(1) Trivial.

(2) $ \xdC Z \xck_{X} \xdC Y$ iff $X \xDN ( \xdC Z \xcs \xdC Y) \xcC \xdC
Z.$ $X \xDN ( \xdC Z \xcs \xdC Y)$ $=$ $X \xcv (X \xfA \xdC ( \xdC Z \xcs
\xdC Y))$ $=$
$X \xcv (X \xfA (Z \xcv Y)).$ So $X \xDN ( \xdC Z \xcs \xdC Y) \xcC \xdC
Z$ $ \xcj $ $(X \xcv (X \xfA (Z \xcv Y))) \xcs Z \xEd \xCQ $ $ \xcj $
$X \xcs Z \xEd \xCQ $ or $(X \xfA (Z \xcv Y)) \xcs Z \xEd \xCQ $ $ \xcj $
$d(X,Z) \xck d(X,Y).$

$ \xcz $
\\[3ex]
\index{Conditions TR-Dist}

\bcd

$\hspace{0.01em}$


\label{Conditions TR-Dist}

Let $U \xEd \xCQ,$ $ \xdy \xcc \xdp (U)$ satisfy $( \xcs ),$ $( \xcv ),$
$ \xCQ \xce \xdy.$

Let $A,B,X_{i} \xbe \xdy,$ $ \xfA: \xdy \xCK \xdy \xcp \xdp (U).$

Let $*$ be a revision function defined for
arbitrary consistent theories on both sides. (This is thus a slight
extension of
the AGM framework, as AGM work with formulas only on the right of $*.)$

{\scriptsize

\begin{tabular}{|c|c|c|}

\hline

\xEH
\xEH
$(*Equiv)$
\xEP
\xEH
\xEH
$ \xcm T \xcr S,$ $ \xcm T' \xcr S',$ $\xch$ $T*T' =S*S',$
\xEP

\hline

\xEH
\xEH
$(*CCL)$
\xEP
\xEH
\xEH
$T*T' $ is a consistent, deductively closed theory,
\xEP

\hline

\xEH
$( \xfA Succ)$
\xEH
$(*Succ)$
\xEP
\xEH
$A \xfA B \xcc B$
\xEH
$T' \xcc T*T',$
\xEP

\hline

\xEH
$( \xfA Con)$
\xEH
$(*Con)$
\xEP
\xEH
$A \xcs B \xEd \xCQ $ $ \xch $ $A \xfA B=A \xcs B$
\xEH
$Con(T \xcv T') $ $\xch$ $T*T' = \ol{T \xcv T' },$
\xEP

\hline

Intuitively,
\xEH
$( \xfA Loop)$
\xEH
$(*Loop)$
\xEP
Using symmetry
\xEH
\xEH
\xEP
$d(X_{0},X_{1}) \xck d(X_{1},X_{2}),$
\xEH
$(X_{1} \xfA (X_{0} \xcv X_{2})) \xcs X_{0} \xEd \xCQ,$
\xEH
$Con(T_{0},T_{1}*(T_{0} \xco T_{2})),$
\xEP
$d(X_{1},X_{2}) \xck d(X_{2},X_{3}),$
\xEH
$(X_{2} \xfA (X_{1} \xcv X_{3})) \xcs X_{1} \xEd \xCQ,$
\xEH
$Con(T_{1},T_{2}*(T_{1} \xco T_{3})),$
\xEP
$d(X_{2},X_{3}) \xck d(X_{3},X_{4})$
\xEH
$(X_{3} \xfA (X_{2} \xcv X_{4})) \xcs X_{2} \xEd \xCQ,$
\xEH
$Con(T_{2},T_{3}*(T_{2} \xco T_{4}))$
\xEP
\Xl
\xEH
\Xl
\xEH
\Xl
\xEP
$d(X_{k-1},X_{k}) \xck d(X_{0},X_{k})$
\xEH
$(X_{k} \xfA (X_{k-1} \xcv X_{0})) \xcs X_{k-1} \xEd \xCQ $
\xEH
$Con(T_{k-1},T_{k}*(T_{k-1} \xco T_{0}))$
\xEP
$\xch$
\xEH
$\xch$
\xEH
$\xch$
\xEP
$d(X_{0},X_{1}) \xck d(X_{0},X_{k}),$
\xEH
$(X_{0} \xfA (X_{k} \xcv X_{1})) \xcs X_{1} \xEd \xCQ$
\xEH
$Con(T_{1},T_{0}*(T_{k} \xco T_{1}))$
\xEP

i.e. transitivity, or absence of
\xEH
\xEH
\xEP

loops involving $<$
\xEH
\xEH
\xEP

\hline

\end{tabular}

}

\index{Proposition TR-Alg-Log}

\ecd

\bp

$\hspace{0.01em}$


\label{Proposition TR-Alg-Log}

The following connections between the logical and the algebraic side might
be the most interesting ones. We will consider in all cases also the
variant
with full theories.

Given $*$ which respects logical equivalence, let $M(T) \xfA M(T'
):=M(T*T' ),$
conversely, given $ \xfA,$ let $T*T':=Th(M(T) \xfA M(T' )).$ We then
have:

{\small

\begin{tabular}{|c|c|c|c|}

\hline

(1.1)
\xEH
$(K*7)$
\xEH
$\xch$
\xEH
$(X \xfA 7)$
\xEP

\cline{1-1}
\cline{3-3}

(1.2)
\xEH
\xEH
$\xci$ $(\xbm dp)$
\xEH
\xEP

\cline{1-1}
\cline{3-3}

(1.3)
\xEH
\xEH
$\xci$ B is the model set for some $\xbf$
\xEH
\xEP

\cline{1-1}
\cline{3-3}

(1.4)
\xEH
\xEH
$\xcI$ in general
\xEH
\xEP

\hline

(2.1)
\xEH
$(*Loop)$
\xEH
$\xch$
\xEH
$(\xfA Loop)$
\xEP

\cline{1-1}
\cline{3-3}

(2.2)
\xEH
\xEH
$\xci$ $(\xbm dp)$
\xEH
\xEP

\cline{1-1}
\cline{3-3}

(2.3)
\xEH
\xEH
$\xci$ all $X_i$ are the model sets for some $\xbf_i$
\xEH
\xEP

\cline{1-1}
\cline{3-3}

(2.4)
\xEH
\xEH
$\xcI$ in general
\xEH
\xEP

\hline

\end{tabular}

}

\index{Proposition TR-Alg-Log Proof}

\ep

\subparagraph{
Proof
}

$\hspace{0.01em}$


(1)

We consider the equivalence of $T*(T' \xcv T'' ) \xcc \ol{(T*T' ) \xcv T''
}$ and
$(M(T) \xfA M(T' )) \xcs M(T'' ) \xcc M(T) \xfA (M(T' ) \xcs M(T'' )).$

(1.1)

$(M(T) \xfA M(T' )) \xcs M(T'' )$ $=$ $M(T*T' ) \xcs M(T'' )$ $=$ $M((T*T'
) \xcv T'' )$ $ \xcc_{(K*7)}$
$M(T*(T' \xcv T'' ))$ $=$ $M(T) \xfA M(T' \xcv T'' )$ $=$ $M(T) \xfA (M(T'
) \xcs M(T'' )).$

(1.2)

$T*(T' \xcv T'' )$ $=$ $Th(M(T) \xfA M(T' \xcv T'' ))$ $=$ $Th(M(T) \xfA
(M(T' ) \xcs M(T'' )))$ $ \xcc_{(X \xfA 7)}$
$Th((M(T) \xfA M(T' )) \xcs M(T'' )))$ $=_{( \xbm dp)}$
$Th(M(T) \xfA M(T' )) \xcv T'' $ $=$ $Th(M((T*T' ) \xcv T'' $ $=$ $
\ol{(T*T' ) \xcv T'' }.$

(1.3)

Let $T'' $ be equivalent to $ \xbf ''.$ We can then replace the use of $(
\xbm dp)$
in the proof of (1.2) by Fact \ref{Fact Log-Form} (3).

(1.4)

By Example \ref{Example TR-Dp} (2), $(K*7)$ may fail, though $(X \xfA
7)$ holds.

(2.1) and (2.2):

$Con(T_{0},T_{1}*(T_{0} \xco T_{2}))$ $ \xcj $ $M(T_{0}) \xcs
M(T_{1}*(T_{0} \xco T_{2})) \xEd \xCQ.$

$M(T_{1}*(T_{0} \xco T_{2}))$ $=$ $M(Th(M(T_{1}) \xfA M(T_{0} \xco
T_{2})))$ $=$ $M(Th(M(T_{1}) \xfA (M(T_{0}) \xcv M(T_{2}))))$ $=_{( \xbm
dp)}$
$M(T_{1}) \xfA (M(T_{0}) \xcv (T_{2})),$ so
$Con(T_{0},T_{1}*(T_{0} \xco T_{2}))$ $ \xcj $ $M(T_{0}) \xcs (M(T_{1})
\xfA (M(T_{0}) \xcv (T_{2}))) \xEd \xCQ.$

Thus, all conditions translate one-to-one, and we use $( \xfA Loop)$ and
$(*Loop)$
to go back and forth.

(2.3):

Let $A:=M(Th(M(T_{1}) \xfA (M(T_{0}) \xcv M(T_{2})))),$ $A':=M(T_{1})
\xfA (M(T_{0}) \xcv (T_{2})),$ then we do
not need $A=A',$ it suffices to have $M(T_{0}) \xcs A \xEd \xCQ $ $ \xcj
$ $M(T_{0}) \xcs A' \xEd \xCQ.$ $A= \wt{A' },$ so
we can use Fact \ref{Fact Def-Clos} $(Cl \xcs +),$ if $T_{0}$ is
equivalent to some $ \xbf_{0}.$

This has to hold for all $T_{i},$ so all $T_{i}$ have to be equivalent to
some $ \xbf_{i}.$

(2.4):

By Proposition \ref{Proposition TR-Alg-Repr}, all distance defined $ \xfA $
satisfy
$( \xfA Loop).$ By Example \ref{Example TR-Dp} (1), $(*Loop)$ may fail.

$ \xcz $
\\[3ex]
\index{Example TR-Dp}

The following Example \ref{Example TR-Dp} shows that, in general, a
revision operation
defined on models via a pseudo-distance by $T*T':=Th(M(T) \xfA_{d}M(T'
))$ might not
satisfy $(*Loop)$ or $(K*7),$ unless we require $ \xfA_{d}$ to preserve
definability.

\be

$\hspace{0.01em}$


\label{Example TR-Dp}

Consider an infinite propositional language $ \xdl.$

Let $X$ be an infinite set of models, $m,$ $m_{1},$ $m_{2}$ be models for
$ \xdl.$
Arrange the models of $ \xdl $ in the real plane s.t. all $x \xbe X$ have
the same
distance $<2$ (in the real plane) from $m,$ $m_{2}$ has distance 2 from
$m,$ and $m_{1}$ has
distance 3 from $m.$

Let $T,$ $T_{1},$ $T_{2}$ be complete (consistent) theories, $T' $ a
theory with infinitely
many models, $M(T)=\{m\},$ $M(T_{1})=\{m_{1}\},$ $M(T_{2})=\{m_{2}\}.$ The
two variants diverge now
slightly:

(1) $M(T' )=X \xcv \{m_{1}\}.$ $T,T',T_{2}$ will be pairwise
inconsistent.

(2) $M(T' )=X \xcv \{m_{1},m_{2}\},$ $M(T'' )=\{m_{1},m_{2}\}.$

Assume in both cases $Th(X)=T',$ so $X$ will not be definable by a
theory.

Now for the results:

Then $M(T) \xfA M(T' )=X,$ but $T*T' =Th(X)=T'.$

(1) We easily verify
$Con(T,T_{2}*(T \xco T)),$ $Con(T_{2},T*(T_{2} \xco T_{1})),$
$Con(T,T_{1}*(T \xco T)),$ $Con(T_{1},T*(T_{1} \xco T' )),$
$Con(T,T' *(T \xco T)),$ and conclude by Loop (i.e. $(*Loop))$
$Con(T_{2},T*(T' \xco T_{2})),$ which
is wrong.

(2) So $T*T' $ is consistent with $T'',$ and
$ \ol{(T*T' ) \xcv T'' }=T''.$ But $T' \xcv T'' =T'',$ and $T*(T' \xcv
T'' )=T_{2} \xEd T'',$ contradicting $(K*7).$

$ \xcz $
\\[3ex]
\index{Proposition TR-Alg-Repr}

\ee

\bp

$\hspace{0.01em}$


\label{Proposition TR-Alg-Repr}

Let $U \xEd \xCQ,$ $ \xdy \xcc \xdp (U)$ be closed under finite $ \xcs $
and finite $ \xcv,$ $ \xCQ \xce \xdy.$

(a) $ \xfA $ is representable by a symmetric pseudo-distance $d:U \xCK U
\xcp Z$ iff $ \xfA $
satisfies $( \xfA Succ)$ and $( \xfA Loop)$ in Conditions \ref{Conditions
TR-Dist}.

(b) $ \xfA $ is representable by an identity respecting symmetric
pseudo-distance
$d:U \xCK U \xcp Z$ iff $ \xfA $ satisfies $( \xfA Succ),$ $( \xfA Con),$
and $( \xfA Loop)$
in Conditions \ref{Conditions TR-Dist}.
\index{Proposition TR-Log-Repr}

\ep

\bp

$\hspace{0.01em}$


\label{Proposition TR-Log-Repr}

Let $ \xdl $ be a propositional language.

(a) A revision operation $*$ is representable by a symmetric consistency
and
definability preserving pseudo-distance iff $*$ satisfies $(*Equiv),$
$(*CCL),$ $(*Succ),$ $(*Loop).$

(b) A revision operation $*$ is representable by a symmetric consistency
and
definability preserving, identity respecting pseudo-distance iff $*$
satisfies $(*Equiv),$ $(*CCL),$ $(*Succ),$ $(*Con),$ $(*Loop).$
\index{Example WeakTR}

\ep

\be

$\hspace{0.01em}$


\label{Example WeakTR}

Note that even when the pseudo-distance is a real distance, the
resulting revision operator $ \xfA_{d}$ does not always permit to
reconstruct the
relations of the distances: revision is a coarse instrument to investigate
distances.

Distances with common start (or end, by symmetry) can always be
compared by looking at the result of revision:

$a \xfA_{d}\{b,b' \}=b$ iff $d(a,b)<d(a,b' ),$

$a \xfA_{d}\{b,b' \}=b' $ iff $d(a,b)>d(a,b' ),$

$a \xfA_{d}\{b,b' \}=\{b,b' \}$ iff $d(a,b)=d(a,b' ).$

This is not the case with arbitrary distances $d(x,y)$ and $d(a,b),$
as this example will show.

\ee

We work in the real plane, with the standard distance, the angles have 120
degrees. $a' $ is closer to $y$ than $x$ is to $y,$ a is closer to $b$
than $x$ is to $y,$
but $a' $ is farther away from $b' $ than $x$ is from $y.$ Similarly for
$b,b'.$
But we cannot distinguish the situation $\{a,b,x,y\}$ and the
situation $\{a',b',x,y\}$ through $ \xfA_{d}.$ (See Diagram \ref{Diagram
WeakTR}):

Seen from a, the distances are in that order: $y,b,x.$

Seen from $a',$ the distances are in that order: $y,b',x.$

Seen from $b,$ the distances are in that order: $y,a,x.$

Seen from $b',$ the distances are in that order: $y,a',x.$

Seen from $y,$ the distances are in that order: $a/b,x.$

Seen from $y,$ the distances are in that order: $a' /b',x.$

Seen from $x,$ the distances are in that order: $y,a/b.$

Seen from $x,$ the distances are in that order: $y,a' /b'.$

Thus, any $c \xfA_{d}C$ will be the same in both situations (with a
interchanged with
$a',$ $b$ with $b' ).$ The same holds for any $X \xfA_{d}C$ where $X$ has
two elements.

Thus, any $C \xfA_{d}D$ will be the same in both situations, when we
interchange a with
$a',$ and $b$ with $b'.$ So we cannot determine by $ \xfA_{d}$ whether
$d(x,y)>d(a,b)$ or not.
$ \xcz $
\\[3ex]

\vspace{10mm}

\begin{diagram}

\label{Diagram WeakTR}
\index{Diagram WeakTR}

\unitlength1.0mm
\begin{picture}(110,90)

\newsavebox{\ZWEIeins}
\savebox{\ZWEIeins}(130,85)[bl]
{

\put(5,50){\line(1,0){30}}
\put(35,50){\line(1,2){6}}
\put(35,50){\line(1,-2){6}}

\put(5,50){\circle*{1.5}}
\put(35,50){\circle*{1.5}}

\put(41,62){\circle*{1.5}}
\put(41,38){\circle*{1.5}}

\put(5,47){$x$}
\put(32,47){$y$}
\put(43,61){$a$}
\put(43,37){$b$}

\put(65,50){\line(1,0){35}}
\put(100,50){\line(1,2){12}}
\put(100,50){\line(1,-2){12}}

\put(65,50){\circle*{1.5}}
\put(100,50){\circle*{1.5}}

\put(112,74){\circle*{1.5}}
\put(112,26){\circle*{1.5}}

\put(65,47){$x$}
\put(97,47){$y$}
\put(114,73){$a'$}
\put(114,25){$b'$}

\put(29,10){Indiscernible by revision}

}

\put(0,0){\usebox{\ZWEIeins}}
\end{picture}

\end{diagram}

\vspace{4mm}

\newpage

\section{
Size
}

\index{Definition Filter}

\bd

$\hspace{0.01em}$


\label{Definition Filter}

A filter is an abstract notion of size,
elements of a filter $ \xdf (X)$ on $X$ are called big subsets of $X,$
their complements
are called small, and the rest have medium size. The dual applies to
ideals
$ \xdi (X),$ this is justified by the trivial fact that $\{X-A:A \xbe \xdf
(X)\}$ is an ideal
iff $ \xdf (X)$ is a filter.

In both definitions, the first two conditions (i.e. $ \xCf (FAll),$ $(I
\xCQ ),$ and
$(F \xfB ),$ $(I \xfb ))$ should hold if the notions shall
have anything to do with usual intuition, and there are reasons to
consider
only the weaker, less idealistic, version of the third.

At the same time, we introduce in rough parallel coherence conditions
which
describe what might happen when we change the reference or base set $X.$
$(R \xfB )$
is very natural, $(R \xfb )$ is more daring, and $(R \xfb \xfb )$ even
more so.
$(R \xcv disj)$ is a cautious combination of $(R \xfB )$ and $(R \xcv ),$
as we avoid using
the same big set several times in comparison, so $(R \xcv )$ is used more
cautiously
here. See Remark \ref{Remark Ref-Class} for more details.

Finally, we give a generalized first order quantifier corresponding to a
(weak) filter. The precise connection is formulated in
Definition \ref{Definition Nabla}, Definition \ref{Definition N-Model},
Definition \ref{Definition NablaAxioms}, and
Proposition \ref{Proposition NablaRepr}, respectively their
relativized versions.

Fix now a base set $X \xEd \xCQ.$

A (weak) filter on or over $X$ is a set $ \xdf (X) \xcc \xdp (X),$ s.t.
(FAll), $(F \xfB ),$ $(F \xcs )$
$((FAll),$ $(F \xfB ),$ $(F \xcs ' )$ respectively) hold.

A filter is called a principal filter iff there is $X' \xcc X$ s.t. $ \xdf
=\{A:$ $X' \xcc A \xcc X\}.$

A filter is called an ultrafilter iff for all $X' \xcc X$ $X' \xbe \xdf
(X)$ or $X-X' \xbe \xdf (X).$

A (weak) ideal on or over $X$ is a set $ \xdi (X) \xcc \xdp (X),$ s.t. $(I
\xCQ ),$ $(I \xfb ),$ $(I \xcv )$
$((I \xCQ ),$ $(I \xfb ),$ $(I \xcv ' )$ respectively) hold.

Finally, we set $ \xdm (X):=\{A \xcc X:A \xce \xdi (X),$ $A \xce \xdf
(X)\},$ the ``medium size'' sets, and
$ \xdm^{+}(X):= \xdm (X) \xcv \xdf (X),$
$ \xdm^{+}(X)$ is the set of subsets of $X,$ which are not small, i.e.
have medium or
large size.

For $(R \xfb )$ and $(R \xfb \xfb )$ closure under set difference is
assumed in the following
table.

{\footnotesize

\begin{tabular}{|c|c|c|c|}

\hline

\multicolumn{4}{|c|}{Optimum} \xEP

\hline

$(FAll)$
\xEH
$(I \xCQ )$
\xEH
\xEH
\xEP

$X \xbe \xdf (X)$
\xEH
$ \xCQ \xbe \xdi (X)$
\xEH
\xEH
$\xcA x\xbf(x) \xcp \xeA x\xbf(x)$
\xEP

\hline

\multicolumn{4}{|c|}{Improvement} \xEP

\hline

$(F \xfB )$
\xEH
$(I \xfb )$
\xEH
$(R \xfB )$
\xEH
\xEP

$A \xcc B \xcc X,$
\xEH
$A \xcc B \xcc X,$
\xEH
$X \xcc Y$ $ \xch $ $\xdi (X) \xcc \xdi (Y)$
\xEH
$\xeA x\xbf(x) \xcu $
\xEP

$A \xbe \xdf (X)$ $ \xch $
\xEH
$B \xbe \xdi (X)$ $ \xch $
\xEH
\xEH
$ \xcA x(\xbf(x)\xcp\xbq(x)) \xcp$
\xEP

$B \xbe \xdf (X)$
\xEH
$A \xbe \xdi (X)$
\xEH
\xEH
$\xeA x\xbq(x)$
\xEP

\hline

\multicolumn{4}{|c|}{Adding small sets} \xEP

\hline

$(F \xcs )$
\xEH
$(I \xcv )$
\xEH
$(R \xfb)$
\xEH
\xEP

$A,B \xbe \xdf (X)$ $ \xch $
\xEH
$A,B \xbe \xdi (X)$ $ \xch $
\xEH
$A,B \xbe \xdi (X)$ $ \xch $
\xEH
$\xeA x\xbf(x) \xcu \xeA x\xbq(x) \xcp $
\xEP

$A \xcs B \xbe \xdf (X)$
\xEH
$A \xcv B \xbe \xdi (X)$
\xEH
$A-B \xbe \xdi (X-$B)
\xEH
$\xeA x(\xbf(x) \xcu \xbq(x)) $
\xEP

\xEH
\xEH
or:
\xEH
\xEP

\xEH
\xEH
$A \xbe \xdf (X),$ $B \xbe \xdi (X)$ $ \xch $
\xEH
\xEP

\xEH
\xEH
$A-B \xbe \xdf (X-$B)
\xEH
\xEP

\hline

\multicolumn{4}{|c|}{Cautious addition} \xEP

\hline

$(F \xcs ' )$
\xEH
$(I \xcv ' )$
\xEH
$(R \xcv disj)$
\xEH
\xEP

$A,B \xbe \xdf (X)$ $ \xch $
\xEH
$A,B \xbe \xdi (X)$ $ \xch $
\xEH
$A \xbe \xdi (X),$ $B \xbe \xdi (Y),$ $X \xcs Y= \xCQ $ $ \xch $
\xEH
$\xeA x\xbf(x) \xcp \xCN\xeA x\xCN\xbf(x)$
\xEP

$A \xcs B \xEd \xCQ.$
\xEH
$A \xcv B \xEd X.$
\xEH
$A \xcv B \xbe \xdi (X \xcv Y)$
\xEH
and  $\xeA x\xbf(x) \xcp \xcE x\xbf(x)$
\xEP

\hline

\multicolumn{4}{|c|}{Bold addition} \xEP

\hline

Ultrafilter
\xEH
(Dual of) Ultrafilter
\xEH
$(R \xfb \xfb)$
\xEH
\xEP

\xEH
\xEH
$A \xbe \xdi (X),$ $B \xce \xdf (X)$ $ \xch $
\xEH
$\xCN \xeA x\xbf(x) \xcp \xeA x\xCN\xbf(x)$
\xEP

\xEH
\xEH
$A-B \xbe \xdi (X-B)$
\xEH
\xEP

\xEH
\xEH
or:
\xEH
\xEP

\xEH
\xEH
$A \xbe \xdf (X),$ $B \xce \xdf (X)$ $ \xch $
\xEH
\xEP

\xEH
\xEH
$A-B \xbe \xdf (X-$B)
\xEH
\xEP

\xEH
\xEH
or:
\xEH
\xEP

\xEH
\xEH
$A \xbe \xdm^+ (X),$ $X \xbe \xdm^+ (Y)$ $ \xch $
\xEH
\xEP

\xEH
\xEH
$A \xbe \xdm^+ (Y)$ - Transitivity of $\xdm^+$
\xEH
\xEP

\hline

\end{tabular}

}

\ed

These notions are related to nonmonotonic logics as follows:

We can say that, normally, $ \xbf $ implies $ \xbq $ iff in a big subset
of all $ \xbf -$cases, $ \xbq $
holds. In preferential terms, $ \xbf $ implies $ \xbq $ iff $ \xbq $ holds
in all minimal
$ \xbf -$models. If $ \xbm $ is the model choice function of a
preferential structure, i.e.
$ \xbm ( \xbf )$ is the set of minimal $ \xbf -$models, then $ \xbm ( \xbf
)$ will be a (the smallest)
big subset of the set of $ \xbf -$models, and the filter over the $ \xbf
-$models is the
pricipal filter generated by $ \xbm ( \xbf ).$

Due to the finite intersection property, filters and ideals work well with
logics:
If $ \xbf $ holds normally, as it holds in a big subset, and so does $
\xbf ',$ then
$ \xbf \xcu \xbf ' $ will normally hold, too, as the intersection of two
big subsets is
big again. This is a nice property, but not justified in all situations,
consider e.g. simple counting of a finite subset. (The question has a
name,
``lottery paradox'': normally no single participant wins, but someone wins
in the
end.) This motivates the weak versions.

Normality defined by (weak or not) filters is a local concept: the filter
defined on $X$ and the one defined on $X' $ might be totally independent.

Seen more abstractly, set properties like e.g. $(R \xfB )$
allow the transfer of big (or small) subsets from one to
another base set (and the conclusions drawn on this basis), and we call
them
``coherence properties''. They are very important, not only for working with
a logic which respects them, but also for soundness and completeness
questions,
often they are at the core of such problems.
\index{Remark Ref-Class}

\br

$\hspace{0.01em}$


\label{Remark Ref-Class}

$(R \xfB )$ corresponds to $(I \xfb )$ and $(F \xfB ):$ If $ \xCf A$ is
small in $X \xcc Y,$ then it will a
fortiori be small in the bigger $Y.$

$(R \xfb )$ says that diminishing base sets by a small
amount will keep small subsets small. This goes in the wrong direction, so
we
have to be careful. We cannot diminish arbitrarily, e.g., if $ \xCf A$ is
a small
subset of $B,$ $ \xCf A$ should not be a small subset of $B-(B-A)=A.$ It
still seems quite
safe, if ``small'' is a robust notion, i.e. defined in an abstract way, and
not
anew for each set, and, if ``small'' is sufficiently far from ``big'', as, for
example in a filter.

There is, however, an important conceptual distinction to make here.
Filters
express ``size'' in an abstract way, in the context of
nonmonotonic logics, $ \xba \xcn \xbb $ iff the set of
$ \xba \xcu \xCN \xbb $ is small in $ \xba.$ But here, we were
interested in ``small'' changes in the reference set $X$ (or $ \xba $ in our
example). So
we have two quite different uses of ``size'', one for
nonmonotonic logics, abstractly expressed by
a filter, the other for coherence conditions. It is possible, but not
necessary,
to consider both essentially the same notions. But we should not forget
that we
have two conceptually different uses of size here.

\er

$(R \xfb \xfb )$ is obviously a stronger variant of $(R \xfb ).$

It and its strength is perhaps best understood as transitivity of the
relation $''. \xbe \xdm^{+}(.)''.$

Now, (in comparison to $(R \xfb ))$ $A' $ can be a medium size subset of
$B.$
As a matter of fact, $(R \xfb \xfb )$ is a very big strengthening of $(R
\xfb ):$ Consider a
principal filter $ \xdf:=\{X \xcc B:$ $B' \xcc X\},$ $b \xbe B'.$ Then
$\{b\}$ has at least medium size, so
any small set $A \xcc B$ is smaller than $\{b\}$ - and this is, of course,
just
rankedness. If we only
have $(R \xfb ),$ then we need the whole generating set $B' $ to see that
$ \xCf A$ is small.
This is the strong substitution property of rankedness: any $b$ as above
will
show that $ \xCf A$ is small.

The more we see size as an abstract notion, and the more we see
``small'' different from ``big'' (or ``medium'' ), the more we can go from one
base set
to another and find the same sizes - the more we have coherence when we
reason
with small and big subsets.
$(R \xfb )$ works with iterated use of ``small'', just as do filters, but
not
weak filters. So it is not surprising that weak filters and $(R \xfb )$ do
not
cooperate well: Let $A,B,C$ be small subsets of $X$ - pairwise disjoint,
and
$A \xcv B \xcv C=X,$ this is possible. By $(R \xfb )$ $B$ and $C$ will be
small in $X-A$, so again
by $(R \xfb )$ $C$ will be small in $(X-A)-B=C,$ but this is absurd.

If we think that filters are too strong, but we still want some coherence,
i.e.
abstract size, we can consider $(R \xcv disj):$ If $ \xCf A$ is a
small subset of $B,$ and $A' $ of $B',$ and $B$ and $B' $ are disjoint,
then $A \xcv A' $ is a
small subset of $B \xcv B'.$ It
expresses a uniform approach to size, or distributivity, if you like. It
holds,
e.g. when we consider a set to be small iff it is smaller than a certain
fraction. The important point is here that by disjointness, the big
subsets do
not get ``used up''. (This property generalizes in a straightforward way to
the
infinite case.)
\index{Fact R-down}

\bfa

$\hspace{0.01em}$


\label{Fact R-down}

The two versions of $(R \xfb )$ and the three versions of $(R \xfb \xfb )$
are each equivalent.
For the third version of $(R \xfb \xfb )$ we use $(I \xfb ).$
\index{Fact R-down Proof}

\efa

\subparagraph{
Proof
}

$\hspace{0.01em}$


For $A,B \xcc X,$ $(X-B)-((X-A)-B)=A-$B.

`` $ \xch $ '': Let $A \xbe \xdf (X),$ $B \xbe \xdi (X),$ so $X-A \xbe \xdi
(X),$ so by prerequisite
$(X-A)-B \xbe \xdi (X-$B), so $A-B=(X-B)-((X-A)-B) \xbe \xdf (X-$B).

`` $ \xci $ '': Let $A,B \xbe \xdi (X),$ so $X-A \xbe \xdf (X),$ so by
prerequisite $(X-A)-B \xbe \xdf (X-$B),
so $A-B=(X-B)-((X-A)-B) \xbe \xdi (X-$B).

The proof for $(R \xfb \xfb )$ is the same for the first two cases.

It remains to show equivalence with the last one. We assume closure under
set difference and union.

$(1) \xch (3):$

Suppose $A \xce \xdm^{+}(Y),$ but $X \xbe \xdm^{+}(Y),$ we show $A \xce
\xdm^{+}(X).$ So $A \xbe \xdi (Y),$ $Y-X \xce \xdf (Y),$
so $A=A-(Y-X) \xbe \xdi (Y-(Y-X))= \xdi (X).$

$(3) \xch (1):$

Suppose $A-B \xce \xdi (X-$B), $B \xce \xdf (X),$ we show $A \xce \xdi
(X).$ By prerequisite $A-B \xbe \xdm^{+}(X-$B),
$X-B \xbe \xdm^{+}(X),$ so $A-B \xbe \xdm^{+}(X),$ so by $(I \xfb )$ $A
\xbe \xdm^{+}(X),$ so $A \xce \xdi (X).$

$ \xcz $
\\[3ex]
\index{Fact Ref-Class-Mu}

\bfa

$\hspace{0.01em}$


\label{Fact Ref-Class-Mu}

If $f(X)$ is the smallest $ \xCf A$ s.t. $A \xbe \xdf (X),$ then, given
the property on the
left, the one on the right follows.

Conversely, when we define $ \xdf (X):=\{X':f(X) \xcc X' \xcc X\},$ given
the property on
the right, the one on the left follows. For this direction, we assume
that we can use the full powerset of some base set $U$ - as is the case
for
the model sets of a finite language. This is perhaps not too bold, as
we mainly want to stress here the intuitive connections, without putting
too much weight on definability questions.

{\footnotesize

\begin{tabular}{|c|c|c|c|}

\hline

(1.1)
\xEH
$(R \xfB )$
\xEH
$ \xch $
\xEH
$( \xbm wOR)$
\xEP

\cline{1-1}
\cline{3-3}

(1.2)
\xEH
\xEH
$ \xci $
\xEH
\xEP

\hline

(2.1)
\xEH
$(R \xfB )+(I \xcv )$
\xEH
$ \xch $
\xEH
$( \xbm OR)$
\xEP

\cline{1-1}
\cline{3-3}

(2.2)
\xEH
\xEH
$ \xci $
\xEH
\xEP

\hline

(3.1)
\xEH
$(R \xfB )+(I \xcv )$
\xEH
$ \xch $
\xEH
$( \xbm PR)$
\xEP

\cline{1-1}
\cline{3-3}

(3.2)
\xEH
\xEH
$ \xci $
\xEH
\xEP

\hline

(4.1)
\xEH
$(R \xcv disj )$
\xEH
$ \xch $
\xEH
$( \xbm disjOR)$
\xEP

\cline{1-1}
\cline{3-3}

(4.2)
\xEH
\xEH
$ \xci $
\xEH
\xEP

\hline

(5.1)
\xEH
$(R \xfb)$
\xEH
$ \xch $
\xEH
$( \xbm CM)$
\xEP

\cline{1-1}
\cline{3-3}

(5.2)
\xEH
\xEH
$ \xci $
\xEH
\xEP

\hline

(6.1)
\xEH
$(R \xfb \xfb)$
\xEH
$ \xch $
\xEH
$( \xbm RatM)$
\xEP

\cline{1-1}
\cline{3-3}

(6.2)
\xEH
\xEH
$ \xci $
\xEH
\xEP

\hline

\end{tabular}

}

\index{Fact Ref-Class-Mu Proof}

\efa

\subparagraph{
Proof
}

$\hspace{0.01em}$


(1.1) $(R \xfB )$ $ \xch $ $( \xbm wOR):$

$X-f(X)$ is small in $X,$ so it is small in $X \xcv Y$ by $(R \xfB ),$ so
$A:=X \xcv Y-(X-f(X)) \xbe \xdf (X \xcv Y),$ but $A \xcc f(X) \xcv Y,$ and
$f(X \xcv Y)$ is the smallest element
of $ \xdf (X \xcv Y).$

(1.2) $( \xbm wOR)$ $ \xch $ $(R \xfB ):$

Let $X \xcc Y,$ $X':=Y-$X. Let $A \xbe \xdi (X),$ so $X-A \xbe \xdf (X),$
so $f(X) \xcc X-$A, so
$f(X \xcv X' ) \xcc f(X) \xcv X' \xcc (X-A) \xcv X' $ by prerequisite, so
$(X \xcv X' )-((X-A) \xcv X' )=A \xbe \xdi (X \xcv X' ).$

(2.1) $(R \xfB )+(I \xcv )$ $ \xch $ $( \xbm OR):$

$X-f(X)$ is small in $X,$ $Y-f(Y)$ is small in $Y,$ so both are small in
$X \xcv Y$ by
$(R \xfB ),$ so $A:=(X-f(X)) \xcv (Y-f(Y))$ is small in $X \xcv Y$ by $(I
\xcv ),$ but
$X \xcv Y-(f(X) \xcv f(Y)) \xcc A,$ so $f(X) \xcv f(Y) \xbe \xdf (X \xcv
Y),$ so, as $f(X \xcv Y)$ is the smallest
element of $ \xdf (X \xcv Y),$ $f(X \xcv Y) \xcc f(X) \xcv f(Y).$

(2.2) $( \xbm OR)$ $ \xch $ $(R \xfB )+(I \xcv ):$

Let again $X \xcc Y,$ $X':=Y-$X. Let $A \xbe \xdi (X),$ so $X-A \xbe \xdf
(X),$ so $f(X) \xcc X-$A. $f(X' ) \xcc X',$
so $f(X \xcv X' ) \xcc f(X) \xcv f(X' ) \xcc (X-A) \xcv X' $ by
prerequisite, so
$(X \xcv X' )-((X-A) \xcv X' )=A \xbe \xdi (X \xcv X' ).$

$(I \xcv )$ holds by definition.

(3.1) $(R \xfB )+(I \xcv )$ $ \xch $ $( \xbm PR):$

Let $X \xcc Y.$ $Y-f(Y)$ is the largest element of $ \xdi (Y),$ $X-f(X)
\xbe \xdi (X) \xcc \xdi (Y)$ by
$(R \xfB ),$ so $(X-f(X)) \xcv (Y-f(Y)) \xbe \xdi (Y)$ by $(I \xcv ),$ so
by ``largest'' $X-f(X) \xcc Y-f(Y),$
so $f(Y) \xcs X \xcc f(X).$

(3.2) $( \xbm PR)$ $ \xch $ $(R \xfB )+(I \xcv )$

Let again $X \xcc Y,$ $X':=Y-$X. Let $A \xbe \xdi (X),$ so $X-A \xbe \xdf
(X),$ so $f(X) \xcc X-$A, so
by prerequisite $f(Y) \xcs X \xcc X-$A, so $f(Y) \xcc X' \xcv (X-$A), so
$(X \xcv X' )-(X' \xcv (X-A))=A \xbe \xdi (Y).$

Again, $(I \xcv )$ holds by definition.

(4.1) $(R \xcv disj)$ $ \xch $ $( \xbm disjOR):$

If $X \xcs Y= \xCQ,$ then (1) $A \xbe \xdi (X),B \xbe \xdi (Y) \xch A
\xcv B \xbe \xdi (X \xcv Y)$ and
(2) $A \xbe \xdf (X),B \xbe \xdf (Y) \xch A \xcv B \xbe \xdf (X \xcv Y)$
are equivalent. (By $X \xcs Y= \xCQ,$
$(X-A) \xcv (Y-B)=(X \xcv Y)-(A \xcv B).)$
So $f(X) \xbe \xdf (X),$ $f(Y) \xbe \xdf (Y)$ $ \xch $ (by prerequisite)
$f(X) \xcv f(Y) \xbe \xdf (X \xcv Y).$ $f(X \xcv Y)$
is the smallest element of $ \xdf (X \xcv Y),$ so $f(X \xcv Y) \xcc f(X)
\xcv f(Y).$

(4.2) $( \xbm disjOR)$ $ \xch $ $(R \xcv disj):$

Let $X \xcc Y,$ $X':=Y-$X. Let $A \xbe \xdi (X),$ $A' \xbe \xdi (X' ),$
so $X-A \xbe \xdf (X),$ $X' -A' \xbe \xdf (X' ),$
so $f(X) \xcc X-$A, $f(X' ) \xcc X' -A',$ so $f(X \xcv X' ) \xcc f(X)
\xcv f(X' ) \xcc (X-A) \xcv (X' -A' )$ by
prerequisite, so $(X \xcv X' )-((X-A) \xcv (X' -A' ))=A \xcv A' \xbe \xdi
(X \xcv X' ).$

(5.1) $(R \xfb )$ $ \xch $ $( \xbm CM):$

$f(X) \xcc Y \xcc X$ $ \xch $ $X-Y \xbe \xdi (X),$ $X-f(X) \xbe \xdi (X)$
$ \xch_{(R \xfb )}$ $A:=(X-f(X))-(X-Y) \xbe \xdi (Y)$ $ \xch $
$Y-A=f(X)-(X-Y) \xbe \xdf (Y)$ $ \xch $ $f(Y) \xcc f(X)-(X-Y) \xcc f(X).$

(5.2) $( \xbm CM)$ $ \xch $ $(R \xfb )$

Let $A \xbe \xdf (X),$ $B \xbe \xdi (X),$ so $f(X) \xcc X-B \xcc X,$ so by
prerequisite $f(X-B) \xcc f(X).$
As $A \xbe \xdf (X),$ $f(X) \xcc A,$ so $f(X-B) \xcc f(X) \xcc A \xcs
(X-B)=A-$B, and $A-B \xbe \xdf (X-$B).

(6.1) $(R \xfb \xfb )$ $ \xch $ $( \xbm RatM):$

Let $X \xcc Y,$ $X \xcs f(Y) \xEd \xCQ.$ If $Y-X \xbe \xdf (Y),$ then
$A:=(Y-X) \xcs f(Y) \xbe \xdf (Y),$ but by
$X \xcs f(Y) \xEd \xCQ $ $A \xcb f(Y),$ contradicting ``smallest'' of
$f(Y).$ So $Y-X \xce \xdf (Y),$ and
by $(R \xfb \xfb )$ $X-f(Y)=(Y-f(Y))-(Y-X) \xbe \xdi (X),$ so $X \xcs f(Y)
\xbe \xdf (X),$ so $f(X) \xcc f(Y) \xcs X.$

(6.2) $( \xbm RatM)$ $ \xch $ $(R \xfb \xfb )$

Let $A \xbe \xdf (Y),$ $B \xce \xdf (Y).$ $B \xce \xdf (Y)$ $ \xch $ $Y-B
\xce \xdi (Y)$ $ \xch $ $(Y-B) \xcs f(Y) \xEd \xCQ.$
Set $X:=Y-$B, so $X \xcs f(Y) \xEd \xCQ,$ $X \xcc Y,$ so $f(X) \xcc f(Y)
\xcs X$ by prerequisite.
$f(Y) \xcc A$ $ \xch $ $f(X) \xcc f(Y) \xcs X=f(Y)-B \xcc A-$B.

$ \xcz $
\\[3ex]
\index{Definition Nabla}

\bd

$\hspace{0.01em}$


\label{Definition Nabla}

Augment the language of first order logic by the new quantifier:
If $ \xbf $ and $ \xbq $ are formulas, then so are $ \xeA x \xbf (x),$ $
\xeA x \xbf (x): \xbq (x),$
for any variable $x.$ The:-versions are the restricted variants.
We call any formula of $ \xdl,$ possibly containing $ \xeA $ a $ \xeA -
\xdl -$formula.
\index{Definition N-Model}

\ed

\bd

$\hspace{0.01em}$


\label{Definition N-Model}

$( \xdn -$Model)

Let $ \xdl $ be a first order language, and $M$ be a $ \xdl -$structure.
Let $ \xdn (M)$ be
a weak filter, or $ \xdn -$system - $ \xdn $ for normal - over $M.$
Define $<M, \xdn (M)>$ $ \xcm $ $ \xbf $ for any $ \xeA - \xdl -$formula
inductively as usual, with
one additional induction step:

$<M, \xdn (M)>$ $ \xcm $ $ \xeA x \xbf (x)$ iff there is $A \xbe \xdn (M)$
s.t. $ \xcA a \xbe A$ $(<M, \xdn (M)>$ $ \xcm $ $ \xbf [a]).$
\index{Definition NablaAxioms}

\ed

\bd

$\hspace{0.01em}$


\label{Definition NablaAxioms}

Let any axiomatization of predicate calculus be given. Augment this with
the axiom schemata

(1) $ \xeA x \xbf (x)$ $ \xcu $ $ \xcA x( \xbf (x) \xcp \xbq (x))$ $ \xcp
$ $ \xeA x \xbq (x),$

(2) $ \xeA x \xbf (x)$ $ \xcp $ $ \xCN \xeA x \xCN \xbf (x),$

(3) $ \xcA x \xbf (x)$ $ \xcp $ $ \xeA x \xbf (x)$ and $ \xeA x \xbf (x)$
$ \xcp $ $ \xcE x \xbf (x),$

(4) $ \xeA x \xbf (x)$ $ \xcr $ $ \xeA y \xbf (y)$ if $x$ does not occur
free in $ \xbf (y)$ and $y$ does not
occur free in $ \xbf (x).$

(for all $ \xbf,$ $ \xbq )$.
\index{Proposition NablaRepr}

\ed

\bp

$\hspace{0.01em}$


\label{Proposition NablaRepr}

The axioms given in Definition \ref{Definition NablaAxioms}
are sound and complete for the semantics of Definition \ref{Definition N-Model}
\index{Definition Nabla-System}

\ep

\bd

$\hspace{0.01em}$


\label{Definition Nabla-System}

Call $ \xdn^{+}(M)=< \xdn (N):N \xcc M>$ a $ \xdn^{+}-system$ or system of
weak filters over $M$ iff
for each $N \xcc M$ $ \xdn (N)$ is a weak filter or $ \xdn -$system over
$N.$
(It suffices to consider the definable subsets of $M.)$
\index{Definition N-Model-System}

\ed

\bd

$\hspace{0.01em}$


\label{Definition N-Model-System}

Let $ \xdl $ be a first order language, and $M$ a $ \xdl -$structure. Let
$ \xdn^{+}(M)$ be
a $ \xdn^{+}-system$ over $M.$

Define $<M, \xdn^{+}(M)>$ $ \xcm $ $ \xbf $ for any formula inductively as
usual, with
the additional induction steps:

1. $<M, \xdn^{+}(M)>$ $ \xcm $ $ \xeA x \xbf (x)$ iff there is $A \xbe
\xdn (M)$ s.t. $ \xcA a \xbe A$ $(<M, \xdn^{+}(M)>$ $ \xcm $ $ \xbf [a]),$

2. $<M, \xdn^{+}(M)>$ $ \xcm $ $ \xeA x \xbf (x): \xbq (x)$ iff there is
$A \xbe \xdn (\{x:<M, \xdn^{+}(M)> \xcm \xbf (x)\})$ s.t.
$ \xcA a \xbe A$ $(<M, \xdn^{+}(M)>$ $ \xcm $ $ \xbq [a]).$
\index{Definition NablaAxioms-System}

\ed

\bd

$\hspace{0.01em}$


\label{Definition NablaAxioms-System}

Extend the logic of first order predicate calculus by adding the axiom
schemata

(1) a. $ \xeA x \xbf (x)$ $ \xcr $ $ \xeA x(x=x): \xbf (x),$
$b.$ $ \xcA x( \xbs (x) \xcr \xbt (x))$ $ \xcu $ $ \xeA x \xbs (x): \xbf
(x)$ $ \xcp $ $ \xeA x \xbt (x): \xbf (x),$

(2) $ \xeA x \xbf (x): \xbq (x)$ $ \xcu $ $ \xcA x( \xbf (x) \xcu \xbq (x)
\xcp \xbj (x))$ $ \xcp $ $ \xeA x \xbf (x): \xbj (x),$

(3) $ \xcE x \xbf (x)$ $ \xcu $ $ \xeA x \xbf (x): \xbq (x)$ $ \xcp $ $
\xCN \xeA x \xbf (x): \xCN \xbq (x),$

(4) $ \xcA x( \xbf (x) \xcp \xbq (x))$ $ \xcp $ $ \xeA x \xbf (x): \xbq
(x)$
and $ \xeA x \xbf (x): \xbq (x)$ $ \xcp $ $[ \xcE x \xbf (x)$ $ \xcp $ $
\xcE x( \xbf (x) \xcu \xbq (x))],$

(5) $ \xeA x \xbf (x): \xbq (x)$ $ \xcr $ $ \xeA y \xbf (y): \xbq (y)$
(under the usual caveat for substitution).

(for all $ \xbf,$ $ \xbq,$ $ \xbj,$ $ \xbs,$ $ \xbt )$.
\index{Proposition NablaRepr-System}

\ed

\bp

$\hspace{0.01em}$


\label{Proposition NablaRepr-System}

The axioms of Definition \ref{Definition NablaAxioms-System} are
sound and complete for the $ \xdn^{+}-semantics$
of $ \xeA $ as defined in Definition \ref{Definition N-Model-System}.

\ep


\begin{thebibliography}{xxxxxx}

\addcontentsline{toc}{section}{References}


\bibitem[AGM85]{AGM85}
C.Alchourron, P.Gardenfors, D.Makinson,
``On the Logic of Theory Change: partial meet contraction and
revision functions'', Journal of Symbolic Logic, Vol. 50,
pp. 510-530, 1985

\bibitem[KLM90]{KLM90}
S.Kraus, D.Lehmann, M.Magidor, ``Nonmonotonic reasoning, preferential
models and cumulative logics'', Artificial Intelligence, 44 (1-2),
p.167-207, July 1990

\bibitem[Sch92]{Sch92}
K.Schlechta: ``Some results on classical preferential models'',
Journal of Logic and Computation, Oxford,
Vol.2, No.6 (1992), $p.$ 675-686

\end{thebibliography}
\end{document}